\documentstyle[bezier]{mathnach}
\Volume{187}
\Year{1997}
\Pages{147\,--\,210}
\Received{July 10, 1995}

\def\emline#1#2#3#4#5#6{%
       \put(#1,#2){\special{em:moveto}}
       \put(#4,#5){\special{em:lineto}}}
\def\newpic#1{}
\def\wtilde{\widetilde}
\def\what{\widehat}
\renewcommand{\fr}{\mbox{for \, all}}
\newcommand{\Aol}{A_0^{(l)}}
\newcommand{\Ail}{A_1^{(l)}}
\newcommand{\Bl}{{B}^{(l)}}
\newcommand{\Blis}{\left[ \left( B^{(l)} \right)^{-1} \right]}
\newcommand{\Btl}{\wtilde{B}^{(l)}}
\newcommand{\Dl}{{D}^{(l)}}
\newcommand{\Dtl}{\wtilde{D}^{(l)}}
\newcommand{\re}{{\rm e}}
\newcommand{\rj}{{\rm j}}
\newcommand{\rjt}{{\rm j}^{\dagger}}
\newcommand{\rJ}{{\rm J}}
\newcommand{\rJt}{{\rm J}^{\dagger}}
\newcommand{\rJo}{{\rm J}_{0}}
\newcommand{\rJot}{{\rm J}_{0}^{\dagger}}
\newcommand{\St}{S^{\dagger}}
\newcommand{\bA}{{\bf A}}
\newcommand{\bB}{{\bf B}}
\newcommand{\bBt}{{\bf B}^{\dagger}}

\newcommand{\bI}{{\bf I}}
\newcommand{\bJ}{{\bf J}}
\newcommand{\bJo}{{\bf J}_0}
\newcommand{\bJot}{{\bf J}_0^\dagger}
\newcommand{\bJt}{{\bf J}^{\dagger}}
\newcommand{\bg}{{\bf g}}
\newcommand{\bt}{{\bf t}}
\newcommand{\bR}{{\bf R}}
\newcommand{\bRo}{{\bf R}_{0}}

\newcommand{\bS}{{\bf S}}
\newcommand{\bs}{{\bf s}}
\newcommand{\bv}{{\bf v}}
\newcommand{\cA}{{\cal A}}
\newcommand{\cB}{{\cal B}}
\newcommand{\cD}{{\cal D}}
\newcommand{\cE}{{\cal E}}
\newcommand{\cF}{{\cal F}}
\newcommand{\cG}{{\cal G}}
\newcommand{\cH}{{\cal H}}
\newcommand{\cI}{{\cal I}}
\newcommand{\cJ}{{\cal J}}
\newcommand{\cK}{{\cal K}}
\newcommand{\cN}{{\cal N}}
\newcommand{\cO}{{\cal O}}
\newcommand{\cP}{{\cal P}}
\newcommand{\cQ}{{\cal Q}}
\newcommand{\cR}{{\cal R}}
\newcommand{\cT}{{\cal T}}
\newcommand{\cW}{{\cal W}}
\newcommand{\hbI}{\hat{\bf I}}
\newcommand{\hbIo}{\hat{\bf I}_0}
\newcommand{\hIo}{\hat{I}_0}
\newcommand{\hIi}{\hat{I}_1}
\newcommand{\hk}{\hat{k}}
\newcommand{\hp}{\hat{p}}
\newcommand{\hP}{\widehat{P}}

\newcommand{\sff}{{\sf f}}
\newcommand{\Sls}{\left[ S_l \right]}
\newcommand{\Slis}{\left[ S_l^{-1} \right]}
\newcommand{\tk}{\tilde{k}}
\newcommand{\tL}{\tilde{L}}
\newcommand{\un}{^{(n)}}
\newcommand{\uo}{^{(0)}}
\newcommand{\Xl}{{\bf X}^{(l)}}
\newcommand{\Xol}{{\bf X}_0^{(l)}}
\newcommand{\Xil}{{\bf X}_1^{(l)}}
\newcommand{\an}{\alpha}
\newcommand{\bn}{\beta}
\newcommand{\gn}{\gamma}

\newcommand{\ad}{_{\alpha}}

\newcommand{\aid}{_{\alpha,i}}
\newcommand{\ajd}{_{\alpha,j}}
\newcommand{\adl}{_{\alpha,l}}
\newcommand{\bjd}{_{\beta,j}}
\newcommand{\bkd}{_{\beta,k}}
\newcommand{\bd}{_{\beta}}
\newcommand{\gd}{_{\gamma}}
\newcommand{\gjd}{_{\gamma,j}}
\newcommand{\abd}{_{\alpha\beta}}
\newcommand{\bad}{_{\beta\alpha}}
\newcommand{\au}{^{(\alpha)}}
\newcommand{\bu}{^{(\beta)}}
\newcommand{\gu}{^{(\gamma)}}
\newcommand{\aut}{^{(\alpha)\dagger}}

\newcommand{\aju}{^{(\alpha,j)}}

\newcommand{\Y}{\Upsilon}
\newcommand{\anum}{\alpha \EQ 1,\,2,\,3,\quad j \EQ 1,2,\,\ldots,n_{\alpha}}

\newcommand{\Om}{{\Omega}}
\newcommand{\omo}{\omega_0}
\newcommand{\Omt}{\Omega^\dagger}
\newcommand{\omos}{\omega^{*}_0}
\newcommand{\Phis}{\Phi^{*}}
\newcommand{\Psis}{\Psi^{*}}
\newcommand{\Rt}{{{\R}^3}}
\newcommand{\Rs}{{{\R}^6}}
\newcommand{\Stwo}{{S^2}}

\newcommand{\Ct}{{{\C}^3}}
\newcommand{\Othree}{{\cal O}\big({\C}^3\big)}
\newcommand{\Opthree}{{\cal O}'\big({\C}^3\big)}
\newcommand{\Osix}{{\cal O}\big({\C}^6\big)}
\newcommand{\Opsix}{{\cal O}'\big({\C}^6\big)}
\newcommand{\Pilh}{\Pi_l^{(\rm hol)}}
\newcommand{\Pibj}{\Pi_b^{(\beta,j)}}
\newcommand{\Sum}{\displaystyle\sum\limits}
\newcommand{\Int}{\displaystyle\int\limits}
\newcommand{\Min}[1]{\mathop{{\rm min}}\limits_{#1}}
\newcommand{\Max}[1]{\mathop{{\rm max}}\limits_{#1}}
\newcommand{\Sup}[1]{\mathop{{\rm sup}}\limits_{#1}}

\newcommand{\Bigcap}{\mathop{\bigcap}\limits}
\newcommand{\reduction}[2]{ #1 \big|_{#2}}

\newcommand{\Lim}[1]{\mathop{{\rm lim}}\limits_{#1}}
\newcommand{\diag}{\mathop{\rm diag}}
\newcommand{\Img}{\mathop{\rm Im}}
\newcommand{\Real}{\mathop{\rm Re}}
\newcommand{\Equal}[1]{\mathop{=}\limits_{#1}}
\newcommand{\be}{\begin{equation}}
\newcommand{\ee}{\end{equation}}

\newcommand{\Frac}[2]{{\displaystyle\frac{#1}{#2}}}
\begin{document}
\Title{Representations \vs{0.2cm}\,for \,the \,Three\,--\,Body
\,{\bf T}\,--\,Matrix, \,Scattering \,Matrices
\,and \,Resolvent \,in \,Unphysical \,Energy \,Sheets}
\Shorttitle{Structure of Three\,--\,Body \,Resolvent \,in \,Unphysical
\,Sheets}
\By{{\sc Alexander K.~Motovilov} of Dubna}
\Names{Motovilov}
\Subjclass{Primary 47\hs{0.02cm}A\hs{0.02cm}40, 81\hs{0.02cm}U\hs{0.02cm}10;
Secondary 32\hs{0.02cm}D\hs{0.02cm}15, 35\hs{0.02cm}J\hs{0.02cm}10.}

\maketitle
\begin{abstract}
Explicit representations for the Faddeev components of the
three\,--\,body T\,--\,mat\-rix continued analytically into unphysical
sheets of the energy Riemann surface are formulated.  According
to the representations, the T\,--\,mat\-rix in unphysical sheets is
explicitly expressed in terms of its components taken in the
physical sheet only.  The representations for the T\,--\,mat\-rix are
then used to construct similar representations for the analytic
continuation of the three\,--\,body scattering matrices and the
resolvent. Domains on unphysical sheets are described where the
representations obtained can be applied. A method for finding
three\,--\,body resonances
based on the
Faddeev differential equations   is
proposed.
\end{abstract}
\newsection{Introduction}
\label{SIntro}
The topic of the paper is closely related to the problem of
studying resonances in three\,--\,body quantum systems. The role played
by such resonances is well known, for example, in the physics of
nuclear reactions and in astrophysics.  More generally,  resonance
is one of the most interesting phenomena in quantum scattering
and the problem of definition and studying resonances attracts a
lot of attention both from physicists and mathematicians.
The literature on resonances is enormous and thus
no attempt will be made here to present an exhaustive summary.
For a history of the subject and a review see e.\,g. the
books~\cite{ReedSimonIII} -- \cite{Kukulin}. The main problems
connected with a definition of  resonance are explicitly
emphasized by {\sc B. Simon} in his survey~\cite{SimonChem}. In 
contrast to the usual (real) spectrum, the resonant one is not
a unitary invariant of a (self\,--\,adjoint) operator and, thus,
``no satisfactory definition of a resonance can depend only on
the structure of a single operator on an abstract Hilbert
space''.  Thus, a consideration of resonances is always
connected with a rather concrete system, model, etc. supposing the
presence of an extra structure like a ``free'' or
``unperturbed'' Hamiltonian or geometric features of the system
or model concerned (see \cite{SimonChem}).

The original idea of interpreting resonances in quantum
mechanics as complex poles of the scattering matrix continued
analytically into unphysical sheets of the energy plane goes
back to {\sc G.\,Gamow}~\cite{GGamow}. In this interpretation,
one actually compares the real dynamics of a system with some of
its ``free'' dynamics (which is an extra structure in the sense
of~\cite{SimonChem}).  Here, resonances also manifest themselves
as energy poles of the continued kernels of the wave operators
(the scattering wave functions).  For radially symmetric
potentials, the interpretation of two\,--\,body resonances as
poles of the analytic continuation of the scattering matrix has
been entirely elaborated in terms of the Jost
functions~\cite{Jost} (see e.\,g.~\cite{ReedSimonIII},
\cite{Baz}, \cite{Newton}, \cite{AlfaroRegge}).

Beginning with {\sc E.\,C. Titchmarsh}~\cite{Titchmarsh} the
resonances have also been interpreted  as  poles of the analytic
continuations of the kernel of the Hamiltonian resolvent (or
matrix elements of the resolvent between suitable states,
see  \cite{ReedSimonIII}, \cite{ReedSimonIV}). In different forms
this idea is realized in the
papers~\cite{Schwinger} -- \cite{Hunziker} (see also the Refs.  quoted
in these papers and in the
books~\cite{ReedSimonIII} -- \cite{Exner}).  In particular, such an
interpretation became the basis for a perturbation theory for
eigenvalues embedded in the continuous spectrum of $N$\,--\,body
Hamiltonians and turning into resonances.  This is a well
studied subject now (see e.\,g. \cite{AlbeverioBook},
\cite{Albeverio72}, \cite{Howland}, \cite{Rauch},
\cite{Hunziker}).  Another variant of a perturbation theory for
the two\,--\,body resonances has been elaborated in
~\cite{AlbeverioKrohn} (see also~\cite{AlbeverioBook}) for
the case where the radius of the interaction tends to zero.

In the case where the interaction potentials are analytic functions of
the coordinates, one can  investigate  resonances by the complex
dilation method~\cite{BalslevCombes} (see also
\cite{ReedSimonIV}, \cite{SimonChem}). The complex dilation
makes it possible to rotate the continuous spectrum of the
$N$\,--\,body Hamiltonian in such a way that certain sectors become
accessible for observation in unphysical sheets neighboring
the physical one~\cite{ReedSimonIV}.  Resonances situated in these
sectors turn out to belong to the discrete spectrum of the transformed
Hamiltonian.  It should be noted that the resonances given
by  complex scaling are proved for a wide class of
interactions to be not only poles of the resolvent but also
poles of continued scattering amplitudes \cite{Hagedorn}. A
number of rigorous results were obtained within the framework of this
method (see e.\,g. the papers~\cite{BalslevCombes},
\cite{SimonRes}, \cite{Hagedorn}, \cite{BalslevSkibsted},
\cite{Siedentop1},  \cite{Siedentop2} and the
book~\cite{ReedSimonIV}).  Regarding the detailed of structure of
the $N$\,--\,body scattering matrix and resolvent of the initial
Hamiltonian continued into unphysical sheets, the complex
dilation method gives not too large capacities.

A relation between analytic properties of the scattering matrix in
the complex plane and the space\,--\,time behavior of wave packets has
been studied (see the books~\cite{Exner} and~\cite{BohmQM} for
references).  Also the problem of completeness, normalization and
orthogonalization of the resonance wave functions (``Ga\-mow
vectors''), i.\,e., solutions of the Schr\"odinger equation corresponding
to the resonance energies had been widely discussed (see e.\,g.
\cite{Baz}, \cite{AlbeverioKrohnPrepr}, \cite{Romo},
\cite{GareevBang}).  Also attempts to interpret resonances and
respective Ga\-mow vectors in the framework of  rigged Hilbert
spaces~\cite{Gelfand} have been undertaken for a number of simple
models (see \cite{ParGorSud}, \cite{BohmJMP}, \cite{Antoniou}).

If the support of an interaction is compact, the resonances of the
two\,--\,body system  may be treated in the  framework of the approach created
by {\sc P. Lax} and {\sc R. Phillips}~\cite{LaxP}  (concerning its
further development see \cite{ReedSimonIII}
and~\cite{Adamjan} -- \cite{Fedorov}).  A main advantage of this
approach consists in the possibility of giving an elegant operator
interpretation of the resonances. Precisely, it allows one describe
resonances as the discrete spectrum of a dissipative operator
representing the generator of the compressed evolution (semi)\,group.
Also,  completeness of the resonance wave functions and an expansion
theorem in the translationally invariant subspace~\cite{LaxP} are
naturally proved~\cite{PavlovDAN71}. It should be noted however that
the Lax\,--\,Phillips approach has strong restrictions on the domain of
its applicability related in particular to the dimension of the configuration
space of the system under consideration (the dimension has to be odd,
and thus the $N$\,--\,body problem already with $N=3$ cannot be
treated).  Up to now, the Lax\,--\,Phillips scheme has been realized in those
scattering problems which generate the Riemann surfaces (though
rather complicated) consisting only of two sheets of the complex
energy plane (see~\cite{PavlovFedorov}, \cite{Fedorov}).  In the
multichannel scattering problems with binary channels, this scheme has been
partly realized in \cite{MotLaxP}.

In the present paper we are concerned with the Faddeev
approach~\cite{Faddeev63} to the three\,--\,body problem.  It is
well known that many important conceptual and constructive
results (see \cite{Faddeev63} -- \cite{MF}) concerning the
physical sheet in the three\,--\,body scattering problem have been
obtained on the basis of the Faddeev
equations~\cite{Faddeev63} and their modifications.  In
particular, the structure of the resolvent and scattering operator
was studied in detail, the completeness of the wave operators
was proved and coordinate asymptotics of the scattering wave
functions were investigated for rapidly decreasing as well as Coulomb
interactions$^{1)}$%
\footnote{$^{1)}$In the last decade, the new, more abstract
approaches~\cite{Enss} -- \cite{DerezinskiNbody} (see also the literature
cited in~\cite{DerezinskiNbody}) having no relation to the
Faddeev\,--\,Yakubovsky techniques~\cite{Faddeev63}, \cite{MF} and
\cite{Yakubovsky}, have been developed to prove the existence and
asymptotic completeness of the $N$\,--\,body wave operators.  In
particular, in \cite{DerezinskiNbody} such a proof is given for
arbitrary $N$ in the case where the pair interactions fall off at
infinity like $r^{-\varrho}$, $\varrho>\sqrt{3}-1$, i.\,e., substantially
slower than  Coulomb potentials. Another approach to proving
the absence of the singular continuous spectrum of the
$N$\,--\,body Hamiltonians
including the hard\,--\,core interactions has been worked out in
\cite{BoutetGS}.}%
~\cite{Faddeev63}, \cite{MerkDiss}, \cite{MF}, \cite{EChAYa}.
Analogous results were obtained also for the singular interactions
described by boundary conditions of various
types~\cite{EChAYa} -- \cite{MakMelMot}.  On the basis of the Faddeev
equations, various methods of investigation of concrete physical systems
were developed in \cite{MF}, \cite{EChAYa}, \cite{Schmid},
\cite{Belyaev}.  Regarding the unphysical sheets, the situation with
using these equations is rather different.  Here, when studying a
concrete three\,--\,body problem, one usually restricts oneself usually to
developing a numerical algorithm to search for resonances in the
unphysical sheets neighboring the physical one. A survey of
physical approaches to a study of resonances in  three\,--\,body
nuclear systems based on the Faddeev equations can be found in
\cite{Orlov} and~\cite{Kukulin}.

The present work is devoted to extending the Faddeev
approach~\cite{Faddeev63}, \cite{MF} to study  the
structure of the three\,--\,body T\,--\,matrix,  resolvent and scattering
matrices continued into unphysical sheets.  We restrict
ourselves to the case  where the interaction potentials fall off
in coordinate space not slower than exponentially.  When constructing
a theory of resonances in the two\,--\,body problem with such
interactions one can use the coordinate as well as momentum
representations. It is clear however that the analytic continuation
of the Faddeev integral equations~\cite{Faddeev63}, \cite{MF} into
unphysical sheets turns out to be a very difficult problem if
the equations are written in the configuration space.  The problem is
that there exist noncompact (cylindrical) domains where the pair
(two\,--\,body) potentials are translationally invariant and, therefore,
do not decrease. At the same time the continued kernels  of the
equations increase exponentially and their solutions have to
increase exponentially, too.  Therefore, the integral terms
 diverge and the coordinate space equations do not make sense.
On the other hand, the integral terms given in the momentum space
can be considered as  Cauchy type integrals admitting an
explicit continuation.  So, at least in the sense of
distributions, a continuation of the momentum space Faddeev
equations becomes a solvable problem.  Actually, in the
paper~\cite{OrlovTur} (see also \cite{Orlov}),  such a
continuation into unphysical sheets neighboring  the
physical one has already been realized formally with the $s$\,--\,wave
Faddeev equations corresponding  to the rank~1 separable
(finite-dimensional) pair potentials.  In the present paper, we
construct a continuation of the equations for the Faddeev components
of the three\,--\,body T\,--\,matrix $T(z)$ in the case of sufficiently
arbitrary pair potentials. We do this not only for the
neighboring unphysical sheets but also for all those remote
sheets of the three\,--\,body Riemann surface where it is possible to
guide the spectral parameter (the energy $z$)  around the
two\,--\,body thresholds.

A central result of the paper consists in a substantiation of the
existence of the analytic continuations (in the weak sense) of
the Faddeev components $M\abd(z)$, $\an,\bn=1,2,3,$ of the
operator $T(z)$ and a construction of explicit (i.\,e., given in
terms of the physical sheet only) representations for them in
the unphysical sheets. These representations are found as a result
of exactly solving the Faddeev equations for the matrix
$M(z)=\{ M\abd(z)\}$ continued into unphysical sheets.  Omitting details
[see formula (\ref{Ml3fin})], the
representations read
\be
\label{MlIntro}
M\big|_{\Pi_l} \EQ {M}\big|_{\Pi_0}\,-\,
{\bf Q}^{\dagger}_l\bJt AS^{-1}_l\bJ{\bf Q}_l
\ee
where $\Pi_l$ denotes the unphysical sheets enumerated by a
(multi)\,index $l\neq 0$. The operator ${\bf Q}_l(z)$ and the
``transposed'' one ${\bf Q}_l^\dagger(z)$ are explicitly
constructed from the matrix $M(z)$ taken in the physical sheet
$\Pi_0$. The numerical matrix $A(z)$ is  an entire function of
$z\in{\C}$.  By $S_l(z)$ [see (\ref{Slcut})] we understand a
truncation (depending essentially on $l$) of the total
three\,--\,body scattering matrix $S(z)$. The operators $\bJ(z)$ and
$\bJt(z)$ realize a restriction of the kernels of the operators
${\bf Q}_l(z)$ and ${\bf Q}_l^\dagger(z)$ on energy shells
respectively in the first and last momentum arguments so that the
products ${\bf Q}_l^\dagger\bJt$ and $\bJ{\bf Q}_l$ have
half\,--\,on\,--\,shell kernels.  Note that the structure of the
representations (\ref{MlIntro}) [(\ref{Ml3fin})] for
$\reduction{M(z)}{\Pi_l}$ is quite analogous to that of the
representations found in the author's recent works~\cite{MotTMF}
and~\cite{MotYaF} for the analytic continuation of the T\,--\,matrix in
the multichannel scattering problems with binary channels.
Representations for the analytic continuations of the three\,--\,body
scattering matrices and resolvent follow immediately from the
representations above for $\reduction{M(z)}{\Pi_l}$ [see Equations
(\ref{Slfin}) and~(\ref{R3l}), respectively].  As follows from
the representations (\ref{Ml3fin}), (\ref{Slfin}) and
(\ref{R3l}), the nontrivial (i.\,e.,  differing from the poles at points of
the discrete spectrum  of the three\,--\,body Hamiltonian)
singularities of the T\,--\,matrix, scattering matrices and resolvent
situated in the unphysical sheet $\Pi_l$ are in fact
singularities of the inverse truncated scattering matrix
$S_l^{-1}(z)$.  Therefore the resonances in the sheet $\Pi_l$,
considered as  poles of the T\,--\,matrix, scattering matrix and
resolvent continued on $\Pi_l$,  are actually those values of the
energy $z$ for which the matrix $S_l(z)$ has the eigenvalue zero.
Of course, in  analogy with a similar property of the two\,--\,body
resonances this result can be considered as quite natural and
rather expected.

Some basic results of the present work were announced in the
report~\cite{MotFewBodyCEBAF}.

Let us describe shortly the structure of the paper.
In Sec. 2, some general notations are given.
Sec. 3 contains  information on analytic properties of
the two\,--\,body T\,-- and scattering matrices
which are used in subsequent sections.
Sec. 4 is devoted to a description of properties of the
matrix $M(z)$ and scattering matrices in the physical sheet of the
energy $z$. In particular,  the domains of the
physical sheet where the half\,--\,on\,--\,shell kernels of $M(z)$ as well as
the truncated scattering matrices $S_l(z)$ may be considered as
holomorphic functions of $z$ are described here.  The kernels and matrices
which are included in the explicit representations (\ref{Ml3fin}),
(\ref{Slfin}) and (\ref{R3l}) mentioned above are introduced.
In Sec. 5 we specify the unphysical sheets
included in a part $\Re$ of the three\,--\,body Riemann surface
which we deal with in the paper.
The analytic continuation of the Faddeev equations into unphysical sheets
is made in Sec. 6. \linebreak
Sec. 7 is devoted to proving the validity of the explicit
representations (\ref{Ml3fin}) for the analytic continuation of the
matrix $M(z)$ into unphysical sheets of the surface $\Re$.
In Sec. 8 we derive analogous representations
[given by formulas (\ref{Slfin})]
for the scattering matrices, and in Sec. 9 the
representations [given by formulas (\ref{R3l})] for the resolvent.
In Sec. 10 we discuss the practical meaning of the
results obtained. In particular we give here a sketch of a method
to search for three\,--\,body resonances on the basis of Faddeev
differential equations in the configuration space.
\newsection{Notations}
\label{SNotations}
A three\,--\,body system in  momentum space is considered.  We
enumerate the particles by the index $\alpha=1,2,3$ and write
$k\ad$, $p\ad$ for the scaled relative momenta~\cite{MF}. For
instance
\begin{equation}
\label{Jacobi}
\begin{array}{rcl}
\vs{0.2cm}k_1 \!\!& = & \!\! \ds
\left[\frac{{\rm m}_2+{\rm m}_3}{2{\rm m}_2 {\rm m}_3}\right]^{1/2}
\cdot
\frac{{\rm m}_2 {\rm p}_3 -{\rm m}_3 {\rm p}_2}
{{\rm m}_2+{\rm m}_3}\,,   \\
p_1 \!\!& = & \!\!
\left[    \Frac{{\rm m}_1 +{\rm m}_2 +{\rm m}_3 }
{ 2{\rm m}_1 ({\rm m}_2 + {\rm m}_3) }    \right]^{1/2}\cdot
\Frac{({\rm m}_2+{\rm m}_3){\rm p}_1-{\rm m}_1 ({\rm p}_2 + {\rm p}_3)}
{ {\rm m}_1 +{\rm m}_2 +{\rm m}_3 }
\end{array}
\end{equation}
with ${\rm m}\ad$ the masses and ${\rm p}\ad$ the momenta of the
particles.  The movement of the center of mass of the system is assumed to
be separated.

Expressions for the relative momenta $k\ad$, $p\ad$ with $\an=2,3$
may be obtained from  (\ref{Jacobi}) by cyclic permutation of
indices.  Usually, we combine the relative momenta  $k\ad$, $p\ad$
into six\,--\,component vectors \vs{0.05cm}$P= \lbrace k\ad, p\ad \rbrace$.
A choice of a
certain pair  $\lbrace k\ad, p\ad \rbrace$ fixes a Cartesian
coordinate system in ${\R}^{6}$.  Transition from one pair of the
momenta to another one means a rotation in ${\R}^{6}$,
$ k\ad=c\abd k\bd+s\abd p\bd$,\,\,
$ p\ad=-s\abd k\bd +c\abd p\bd,$\,\,
with coefficients $c\abd$, $s\abd$
depending on the particle masses only~\cite{MF} such that
$ -1< c\abd < 0,$ $s\abd^2=$ $= 1-c\abd^2$,
$c\bad=c\abd$ and $s\bad=-s\abd$, $\bn\neq\an$.

In the momentum representation, the Hamiltonian $H$
of the three\,--\,body system under consideration is defined by
          $$
(Hf)(P)\;=\;P^{2}f(P)\,+\,\Sum_{\an=1}^{3} \,(v\ad f)(P)\,, \quad
P^2 \EQ k\ad^2\,+\,p\ad^2\,, \quad f\in \cH_0 \; \equiv \; L_2\big(\Rs\big)\,,
          $$
with $v\ad$ the pair potentials which are assumed to be
integral operators in $k\ad$ with the kernels  $v\ad(k\ad, k'\ad)$.
Hereafter, by the square of a vector of ${\R}^{N}$, $N=3$ or $N=6$,
we understand the square of its modulus, e.\,g., $P^2=|P|^2$.

For the sake of definiteness we suppose all the potentials
$v\ad$, $\an=1,2,3,$ to be local.  This means that the kernel of
$v\ad$ depends on the difference of the variables $k\ad$ and $k'\ad$
only:  $v\ad(k\ad, k'\ad)$ $=v\ad(k\ad - k'\ad)$.  Actually, we
consider two variants of the potentials $v\ad$.  In the first one,
$v\ad (k)$ are holomorphic functions of the variable $k\in\Ct$
satisfying the estimate
\vs{0.2cm}\be
\label{vpot}
|v\ad (k)|\LE \Frac{c}{(1+|k|)^{\theta_{0}}}
\, e^{a_{0}|{\rm Im} k|} \quad \mbox{for all} \quad k\in{\C}^3
\ee
for some  $c>0,$ $a_0 >0$ and $\theta_{0}\in (3/2,\, 2)$.\,\,
In the second variant, the potentials $v\ad (k)$ are holomorphic
functions of $k$ in the strip
$$
W_{2b}\EQ\{k:\,k\in\Ct,\,|{\rm{Im}}\,k|<2b\}\,, \quad b \GT 0\,,
$$
only and satisfy the condition
\be
\label{vpotb}
|v\ad (k)|\LE \Frac{c}{(1+|k|)^{\theta_{0}}}
 \quad \mbox{for all} \;\; k \;\;{\rm such\;\;that} \;\; |\Img \,k| \LT 2b\,.
\ee
In  both variants the potentials $v\ad$ are supposed to be such
that $v\ad (-k)=$ $\overline{v\ad (k)}$ for any $k\in{\R}^3$.  The
latter condition guarantees  self\,--\,adjointness of the Hamiltonian
$H$ on the set $\cD(H)=\Big\{ f:\, \int_{\R^6} \big(1+P^{2}\big)^{2}|f(P)|^{2} 
\,dP <\infty \Big\}$ (see Theorem~1.1 of \cite{Faddeev63}).

Note that in the first variant, the requirement of holomorphy in all
$\Ct$ and no more than exponential growth  in $|{\rm
Im}\, k|$ (2.2) mean that the potentials $v\ad$ have  compact support in
the coordinate space.  In the second variant, the potentials $v\ad
(k)$, rewritten in the coordinate representation, decrease
exponentially.

By $h\ad$, $(h\ad f)(k\ad)=$ $k\ad ^{2} f(k\ad)+(v\ad f)(k\ad)$, we
denote the Hamiltonian of the pair subsystem $\an$.  The operator
$h\ad$ acts in $L_2 \big(\Rt\big)$.  Due to  condition (\ref{vpot}) or
(\ref{vpotb}), its discrete spectrum  $\sigma_{d}(h\ad)$ is negative
and finite~\cite{ReedSimonIII}.  We enumerate the eigenvalues
$\lambda\ajd\in\sigma_{d}(h\ad),$ $\lambda\ajd < 0$, $ j =
1,2,\,\ldots,n\ad$, $n\ad<\infty$, taking into account their multiplicities:
each eigenvalue being repeated a number of times equal to
its multiplicity. The maximum of these numbers is denoted by
$\lambda_{\rm max}$,\, $\lambda_{\rm max}=\max_{\an,j}\hs{0.02cm}\lambda
\ajd<0.$
The notation $\psi\ajd(k\ad)$ is used for the respective eigenfunctions.

By $\sigma_d(H)$ and $\sigma_c(H)$ we denote respectively the
discrete and (absolutely) continuous components of the spectrum
$\sigma(H)$ of the Hamiltonian $H$. Note that
$\sigma_c(H)=$ \linebreak
$ =[\lambda_{\rm min},+\infty)$ with $\lambda_{\rm min}=
\min_{\an ,j} \hs{0.02cm}\lambda\ajd$.

The notation $H_0$ is used for the operator of kinetic energy, $(H_0
f)(P)=P^{2}f(P)$, while $R_0 (z)$ and $R(z)$ stand for the resolvents
of the operators $H_0$ and $H$, $R_0 (z)=$ \linebreak $=(H_0 -zI)^{-1}$ and
$R(z)=(H -zI)^{-1}$. Here, $I$ is the identity operator in  $\cH_0$.

Let $ M\abd (z) = \delta\abd v\ad -v\ad R(z) v\bd$, $\an,\bn
=1,2,3,$ be the Faddeev components (these components were
introduced in formulae (3.7) in
paper~\cite{Faddeev63} by  {\sc L.\,D. Faddeev}) of the three\,--\,body
$\mbox{T\,--\,matrix}^{2)}$%
\footnote{$^{2)}$ Generally speaking, this operator does not
possess a matrix structure. The traditional term ``matrix'' came
from multi\,--\,channel scattering problems with binary channels
where the operators defined analogously to $T(z)$ are indeed
matrices.}
$ T(z)=V-VR(z)V $ where $V=v_1 +v_1 +v_3 $.  The operators
$M\abd(z)$ satisfy the Faddeev equations (the famous
equations~(3.9) of \cite{Faddeev63})
\be
\label{Fadin}
M\abd (z) \EQ \delta\abd \bt\ad(z) \,-\,
\bt\ad(z)\hs{0.02cm} R_0 (z)\Sum_{\gn\neq\an} M_{\gamma\bn}(z)\,,
\quad \an \EQ 1,\,2,\,3\,.
\ee
Here, the operator ${\bt}\ad(z)$, $\alpha=1,2,3,$ has the kernel
\be
\label{tpair3}
{\bt}\ad(P,P',z) \EQ t\ad\big(k\ad,k'\ad,z-p\ad^2\big)
\hs{0.02cm}\delta(p\ad-p'\ad)\,,
\ee
where $t\ad(k,k',z)$ stands for the kernel of the operator
$t\ad(z)=v\ad-v\ad r\ad(z)v\ad$ with $r\ad(z)=(h\ad-z)^{-1}$.
The operator $t\ad(z)$, called the T\,--\,matrix (or transition
operator)  for the pair subsystem $\alpha$, satisfies in
turn the Lippmann\,--\,Schwinger equation
\be
\label{LScht2}
t\ad(z) \EQ v\ad\,-\,v\ad r_0^{(\an)}(z)\hs{0.02cm}t\ad(z)\,,
\ee
where $r_0^{(\an)}(z)$, $r_0^{(\an)}(z)=\big(h_0^{(\an)}-z\big)^{-1},$ is the
resolvent of the kinetic energy operator $h_0^{(\an)}$ for the
subsystem $\an$, $\big(h_0^{(\an)}f\ad\big)(k\ad)=k\ad^2 f\ad (k\ad),$
$f\ad\in L_2 \big(\Rt\big).$

It is convenient to rewrite the system (\ref{Fadin}) in the matrix
form
\be
\label{MFE}
M(z) \EQ \bt (z) \,-\, \bt (z)\hs{0.02cm} \bRo (z) \Y M(z)\,,
\ee
with $\bt (z)=\diag\hs{0.02cm}\{\bt_1(z),\bt_2(z),\bt_3(z)\}$ and
$\bRo (z)=\diag\hs{0.02cm}\{R_0(z),R_0(z),R_0(z)\}$. By $\Y$
we denote a $3\times 3$ matrix with elements
$\Y\abd=1-\delta\abd.$
$ M(z) $ is a $3 \times 3$ operator matrix
constructed from the components
$M\abd(z)$, $ M=\{ M\abd \},$ $\an,\bn=1,2,3$. The matrices
$M$, $\bt$, $\bRo$
and $\Y$ are considered as operators in the Hilbert space
${\cal G}_0=\bigoplus\nolimits_{\an=1}^{3} L_2\big(\Rs\big)$.

The matrix $M(z)$ obeys also
an alternative variant of the Faddeev equations,
\vs{0.2cm}\be
\label{MFEAl}
M(z) \EQ \bt (z) \,-\, M(z)\hs{0.02cm}\bRo (z)\hs{0.02cm} \Y \hs{0.02cm}
\bt (z) \,.
\ee

We shall also use  the iterated equations (\ref{MFE})
and (\ref{MFEAl}),
\be
\label{MFEmn}
M(z) \EQ \Sum_{k=0}^{m+n+1} \cQ^{(k)}
\,+\, \cQ^{(m)}\hs{0.02cm}\bRo\hs{0.02cm}\Y\hs{0.02cm} M\hs{0.02cm}\Y
\hs{0.02cm}
\bRo \hs{0.02cm} \cQ^{(n)}\,, \quad m,\,n\GE 0\,,
\ee
with
\vs{0.2cm}\be
\label{Qitert}
 \cQ^{(k)}(z) \EQ \bigl( - \bt(z)\hs{0.02cm}\bRo(z)\hs{0.02cm}\Y \bigr)^k
 \bt(z) \EQ \bt(z) \bigl( - \Y\hs{0.02cm}\bRo(z) \hs{0.02cm}\bt(z)  \bigr)^k,
\ee
the iterations of the
absolute term  $\cQ\uo(z)=\bt(z)$ in Eqs.~(\ref{MFE})
and~(\ref{MFEAl}).

The resolvent $R(z)$ is expressed in terms of the matrix $M(z)$
by the formula~\cite{MF}
\vs{0.2cm}\be
\label{RMR}
R(z) \EQ R_0(z)\,-\,R_0(z)\hs{0.02cm} \Om \hs{0.02cm}M(z) \hs{0.02cm}\Omt
\hs{0.02cm}R_0(z)\,,
\ee
where  $\Om: \cG_0\rightarrow\cH_0$ stands for
the matrix\,--\,row $\Om=(1,\,\, 1,\,\, 1)$ and
$\Omt=\Om^{*}=(1,\,\, 1,\,\, 1)^{\dagger}$.
Hereafter, the symbol ``$\dagger$'' means transposition.

Throughout the paper we understand by $\sqrt{z-\lambda}$,
$z\!\in\!{\C}$, $\lambda\!\in\!{\R},$\,  the main branch of the
function $(z-\lambda)^{1/2}$.  By $\hat{q}$ we  usually denote the
unit vector in the direction  $q\!\in\!{\R}^N$,\,\,
$\hat{q}={q}/{|q|}$, and by $S^{N-1}$ the unit sphere in ${\R}^N$,
$\hat{q}\!\in\!  S^{N-1}$.  The inner product in ${\R}^N$ is denoted by
$(\,\cdot\, ,\,\cdot\,)$. The  notation $\langle\,\cdot\, ,\,\cdot\,
\rangle $ is used for inner products in infinite\,--\,dimensional
Hilbert spaces.

Let $\cH^{\aju}=L_2\big(\Rt\big)$ and $\cH^{\au}=\bigoplus
\nolimits_{j=1}^{n\ad}
\cH^{\aju}$.  By $\Psi\ad$ we denote the operator acting from  ${\cal
H}^{\au}$ into ${\cal H}_0$ defined by
              $$
(\Psi\ad f)(P) \EQ
\Sum\limits_{j=1}^{n\ad}\psi\ajd (k\ad)f_j(p\ad) \,, \quad
f\EQ ( f_1, f_2,\,\ldots , f_{n\ad})^{\dag}\,.
            $$
The notation $\Psi\ad^{*}$ is used for the adjoint operator
of  $\Psi\ad$.  By $\Psi$ we denote the block\,--\,diagonal
matrix operator $\Psi =\diag\hs{0.02cm}\{\Psi_1, \Psi_2, \Psi_3\}$ which
acts from  ${\cal H}_1 =\bigoplus\nolimits_{\an=1}^{3} {\cal H}\au$ into ${\cal
G}_0$ and by $\Psi^{*}$ the adjoint operator of $\Psi$.  Analogous
to $\Psi\ad$, $\Psi\ad^{*}$, $\Psi$ and  $\Psi^{*}$ we introduce
operators $\Phi\ad$, $\Phi\ad^{*}$, $\Phi$ and  $\Phi^{*}$, which are
obtained from the former by replacement of the eigenfunctions
$\psi\ajd(k\ad)$ with the ``form factors''
$\phi\ajd(k\ad)=(v\ad\psi\ajd)(k\ad)$, $\an=1,2,3$, $j=1,2,\,\ldots,n\ad$.

The two\,--\,body T\,--\,matrix  $t\ad(z)$ is
known~(see \S\hs{0.02cm}4 of \cite{Faddeev63} or \S\hs{0.02cm}1 of
Chapter~III
of \cite{MF}) to be an analytic operator\,--\,valued function
of the variable $z\in{\C}\backslash[0,+\infty)$ having simple poles at the
points $z\in\sigma_d(h\ad)$. Its kernel admits the
representation
\be
\label{tdistr}
t\ad(k,k',z) \EQ - \Sum_{j=1}^{n\ad} \,\Frac{\phi\ajd(k)\hs{0.02cm}
\overline{\phi\ajd}(k')}{\lambda\ajd-z}
 \,+\,\tilde{t}\ad(k,k',z)\,,
\ee
where $\tilde{t}\ad(k,k',z)$ is a  holomorphic function in  the variable
$z\in{\C}\backslash[0,+\infty)$.  Therefore
\be
\label{tcd}
\bt\ad(z)\;=\;-\,\Phi\hs{0.02cm}\ad\hs{0.02cm}\bg\ad(z)\hs{0.02cm}
\Phi\ad^{*}\,+\,\tilde{\bt}\ad(z)
\ee
where $\tilde{\bt}\ad(z)$ stands for the operator having the kernel
$\tilde{t}\ad\big(k\ad,k'\ad,z-p\ad^{2}\big)\delta(p\ad-p'\ad)$.  At the
same time \vs{0.05cm}$\bg\ad(z)$=$\diag\hs{0.02cm}\{ g_{\alpha,1}(z),\,\ldots,
g_{\an,n\ad}(z)\}$ is a block\,--\,diagonal operator matrix whose
elements  $g\ajd(z)$ are the  operators in ${\cal H}^{(\an,j)}$
with the singular kernels
$ g\ajd(p\ad,p'\ad,z)={\delta(p\ad-p'\ad)}/\big( {\lambda\ajd-z+p\ad^2}\big).$

Below, we consider restrictions of different functions
on the energy (or mass) shell
\be
\label{EnergyS-twobody}
   k \EQ \sqrt{z} \,\hat{k}\,,\quad  \hat{k}\in S^2\,,
\ee
in the two\,--\,body problem and on the energy (or mass) shells
\be
\label{EnergyS3}
   P \EQ \sqrt{z} \,\what{P}\,,\quad  \what{P}\in S^5\,,
\ee
and
\vs{0.2cm}\be \hs{1.0cm}
\label{EnergyS2}
  p\ad \EQ \sqrt{z-\lambda\ajd} \,\hat{p}\ajd\,,\quad
    \hat{p}\ajd\in S^2, \quad \alpha = 1,2,3, \;\; j= 1,2,\,\ldots, n_{\alpha}
        \,,
\ee
in the problem of three particles.  In the last case, the sets
(\ref{EnergyS3}) and (\ref{EnergyS2}) are called respectively
three\,--\,body and two\,--\,body energy shells.

Let ${\cal O}\big({\C}^{N}\big)$ be a linear
space of test functions represented by
the Fourier transform of functions belonging to
$C_{0}^{\infty}\big({\R}^{N}\big)$ (we deal with  $N=3$ or $N=6$ only).
We mean $f\in{\cal O}\big({\C}^{N}\big)$ if
$$
f(q) \EQ\int_{{\R}^N} dx \exp\hs{0.02cm}\{ i\hs{0.02cm}(q_1 x_1+\,\cdots\,+
q_N x_N)\}\,f^{\#}(x) \quad \mbox{for some} \quad  f^{\#}\in C_{0}^{\infty}
\big({\R}^{N}\big)\,.
$$
where, $q=(q_1,q_2,\,\ldots,q_N)\in{\C}^N$, $x=(x_1,x_2,\,\ldots,x_N)
\in{\R}^N$.
Every  $f(q)\in{\cal O}\big({\C}^{N}\big)$
is a holomorphic (entire) function of the variable
$q\in{\C}^{N}$ satisfying the estimates
$$
\left| \frac{\partial^{|m|}}
{\partial q^{m_1}_{1}\cdot \ldots\cdot\partial q^{m_N}_{N}}
f(q) \right| \EQ c\cdot
{\exp\hs{0.02cm}(a|\Img q|)}{(1+|q|)^{-\theta}}\,,
$$
where $a$ is the radius of a ball centered at the origin and containing
the support of the Fourier pre\,--\,image of $f$ in ${\R}^N $,
$
 |m|=m_1+\cdots+m_N \,,
$
and
$
|\Img q | = \sqrt{\sum_{j=1}^{N} |\Img q_j|^2 }.
$
For $\theta$ one can take an arbitrary positive number.
The coefficient $c>0$ depends on $f$, $\theta$ and $ m=(m_1,\ldots,m_N)$.

A sequence of functions $\{f_j\}$ in the space ${\cal O}\big({\C}^{N}\big)$
is said to be convergent to a function $f\in{\cal O}\big({\C}^{N}\big)$ if
for all $\theta>0$ and $m_1,\,\ldots,m_N=0,1,2,\,\ldots$ the equalities
$$
\Lim{j\rightarrow\infty} \Sup{q\in{\C}^N}
(1+|q|)^\theta\exp(-a|\Img q|)\!
\left| \frac{\partial^{|m|}}
{\partial q^{m_1}_{1}\cdot\,\ldots\,\cdot\partial q^{m_N}_{N}}
\big(f_j(q)-f(q)\big) \!\right| \EQ 0
$$
hold for some $a$ not depending on $m$.
Elements of the dual space ${\cal O}'\big({\C}^{N}\big)$, the linear
continuous functionals over ${\cal O}\big({\C}^{N}\big)$,
are  usually called generalized functions or distributions$^{3)}$%
\footnote{$^{3)}$Note that, in fact, we could consider  narrower
classes of distributions over spaces of test functions
holomorphic only in those domains [described in terms of the
energy shells (\ref{EnergyS3}) and (\ref{EnergyS2})] where the
scattering matrices and some half\,--\,on\,--\,shell kernels encountered
below may be continued into the physical energy sheet. However the
results [representations (\ref{Ml3fin}), (\ref{Slfin}) and
(\ref{R3l})] do not depend on such a choice.}
(see e.\,g.~\cite{Gelfand}, \cite{GShilov}).

Let $\rj(z)$ be the operator which restricts
functions $f(k)$, $k\in \Rt$,
on the shell (\ref{EnergyS-twobody}) at \linebreak
$z=E\pm i\,0$, $E>0$,
and continuing them if possible, on a domain
of complex values of the energy  $z$.
On the set ${\cal O}\big(\Ct\big)$, the operator  $\rj(z)$ is defined by
\be
\label{j2}
\bigl(\rj(z)f\bigr)\!\big(\hat{k}\big) \EQ f\big(\sqrt{z}\,\hat{k}\big)\,.
\ee
Its kernel is a holomorphic generalized function
(distribution)~\cite{GShilov},
$\rj\big(\hat{k},k',z\big)=$ \linebreak
$= \delta\big(\sqrt{z}\,\hat{k}-k'\big)$.

By $\rjt(z)$ we denote the transposed operator of $\rj(z)$.
For any $\varphi\in L_2(S^2)$ this operator gives
generalized function (distribution) over  ${\cal O}\big(\Ct\big),$
\be
\label{j2t}
(\rjt(z)\varphi)(k) \EQ
\int_{S^2} d\hat{k'}\,\delta\big(k-\sqrt{z}\,\hat{k'}\big)
\,\varphi\big(\hat{k'}\big) \EQ
\Frac{\delta\big(|k|-\sqrt{z}\,\big)}{z} \,\varphi\big(\hat{k}\big)\,,
\ee
i.\,e.,
\be
\label{j2tint}
(\rjt(z)\varphi,f) \EQ \int_{S^2}d\hat{k}f\big(\sqrt{z}\,\hat{k}\big)
\hs{0.02cm}\varphi\big(\hat{k}\big)\,,
\quad f\in \Othree\,.
\ee

Let $\rJ\ajd(z)$,  $\alpha = 1,2,\,\ldots, ~ j=1,2,\ldots\,$, be the operator
of restriction on the shell
(\ref{EnergyS2}). Its action on  $\Othree$ is defined by
$$
(\rJ\ajd(z)f)(\hat{p}\ad) \EQ f\big(\sqrt{z-\lambda\ajd}\,\,\hat{p}\ad\big)\,,
\quad \an \EQ 1,\,2,\,3, \quad j \EQ 1,\,2,\ldots,n\ad\,.
$$
The  operators $\rJ\ajd(z)$ have the kernels
               $$
\rJ\ajd(\hat{p}\ad,p'\ad,z) \EQ \delta\big(\sqrt{z-\lambda\ajd}\,\hat{p}\ad-
p'\ad\big)\,.
               $$
By $\rJo(z)$ we denote the operator of restriction
on the shell  (\ref{EnergyS3}). On $ \Osix $
this operator is defined by
$
\big(\rJo(z)f\big)\big(\what{P}\big) = f\big(\sqrt{z}\widehat{P}\big).
$
The notations $\rJt\ajd(z)$ and $\rJot(z)$ are used for the
respective transposed operators. Their actions are defined
similarly to (\ref{j2t}), (\ref{j2tint}) by
$$
\big(\rJt\ajd(z)\varphi\big)(p\ad) \EQ
\int_{S^2}d\hat{p}'\ad
\delta\big(p\ad-\sqrt{z-\lambda\ajd}\,\hat{p}'\ad\big)
\hs{0.02cm}\varphi(\hat{p}'\ad)\,,\quad \varphi\in \widehat{\cH}\aju,
$$
$$
\big(\rJot(z)\varphi\big)(P) \EQ
\int_{S^5}d\what{P}'
\delta\big(P-\sqrt{z}\what{P}'\big)
\hs{0.02cm}\varphi\big(\what{P}'\big)\,,\quad \varphi\in \widehat{\cH}_0\,,
$$
where $\widehat{\cH}\aju\equiv L_2\big(S^2\big)$ and $\widehat{\cH}_0\equiv
L_2\big(S^5\big).$ The generalized functions $\rJt\ajd(z)\varphi$ and
$\rJot(z)\varphi$ are elements of the spaces  $\Opthree$ and
$\Opsix $ of distributions over $\Othree $ and $\Osix $,
respectively.

The operators $\rJ\ajd$ and $\rJt\ajd$ are then combined into the
block\,--\,diagonal matrices $\rJ\au(z)=$
$\diag\hs{0.02cm}\{\rJ_{\alpha,1}(z),\,\ldots,\rJ_{\alpha,n\ad}(z)\}$ and
$\rJ^{(\alpha)\dagger}(z)=$ $ \diag\hs{0.02cm}\big\{ \rJt_{\alpha,1}(z),\,
\ldots, \rJt_{\alpha,n\ad}(z)\big\}. $ The latter  are used to
construct the
operators
            $$
\rJ_1 (z) \EQ \diag\hs{0.02cm}\big\{\rJ^{(1)}(z),\rJ^{(2)}(z),
\rJ^{(3)}(z) \big\}\,, \quad \rJt_1 (z) \EQ
\diag\hs{0.02cm}\big\{\rJ^{(1)\dagger}(z),\rJ^{(2)\dagger}(z),
\rJ^{(3)\dagger}(z) \big\}\,.
            $$
The action of $\rJ\au(z)$ and
$\rJ_1(z)$ on elements of the spaces $\cO\au=$
$\times_{\an=1}^{n\ad} \cO\aju$, $\cO\aju\equiv$ $ \equiv \Othree$ and
$\cO_1={\times}_{\an=1}^3 \cO\au$, respectively, can be understood from the
definition of the operators $\rJ\ajd(z)$.  The operators
$\rJ^{(\alpha)\dagger}(z)$ act from
$\widehat{\cH}\au\equiv\oplus_{j=1}^{n\ad} \widehat{\cH}\aju $ to the
space $\cO\au{}'$ of distributions  over $\cO\au$ and
the operator $\rJt_1(z)$  from
$\widehat{\cH}_1
= \oplus_{\an=1}^{3}\widehat{\cH}\aju$ to the space
$\cO'_1$ of distributions over $\cO_1$.

Finally, we use the block\,--\,diagonal operator $3\times
3$\,--\,matrices
            $$
	    \bJo(z) \EQ
\diag \hs{0.02cm}\{ \rJo(z),\rJo(z),\rJo(z)\}\,,
\quad \bJot(z) \EQ \diag\hs{0.02cm}\big\{ \rJot(z),\rJot(z),\rJot(z)\big\}
             $$
constructed from the operators $\rJo(z)$ and $\rJot(z)$,
respectively, as well as the operators
$\bJ(z)=\diag\hs{0.02cm}\{\rJo(z),\rJ_1(z)\}$ and
$\bJt(z)=\diag\hs{0.02cm}\big\{\rJot(z),\rJt_1(z)\big\}.$ The actions of these
operators
is clear from the definitions of the operators $\rJo,$ $\rJ_1$,
$\rJot$ and $\rJt_1$.  In particular, the operator  $\bJt(z)$
acts from  the space $\what{\cG}_0=\oplus_{\an=1}^3 \widehat{\cH}_0$
to the space ${\times}_{\an=1}^3 \cO'\big({\C}^6\big)$.

The identity operators in the spaces $\widehat{\cH}_0$,
$\what{\cG}_0$,  $\widehat{\cH}_1$ and $\what{\cH}_0\oplus\what{\cH}_1$
are denoted by $\hIo$, $\hbIo$, $\hIi$ and $\hbI$, respectively.


\newsection{Analytic continuation of the two\,--\,body
T\,-- and scattering matrices}
\label{St2body}
In this section we recall some analytic properties of pair
T\,--\,matrices  the knowledge of which will be necessary  for
posing the three\,--\,body problem below. Note that these properties are
well known (see e.\,g. ~\cite{ReedSimonIII}, \cite{AlfaroRegge}
and~\cite{Orlov}) for a wide class of the potentials
$v\ad$.  As a matter of fact, we give here only an
explicit representation for the two\,--\,body T\,--\,matrix in
an unphysical sheet which is a particular case of the explicit
representations constructed in the author's work~\cite{MotTMF}
(see Theorem~2 in~\cite{MotTMF} and comments to it) for a
more general situation of analytic continuation of the T\,--\,matrix on
unphysical sheets in the multichannel problem with binary
channels. On the other hand, using this simple example related
to the two\,--\,body problem we can demonstrate the main features of the
scheme which we apply later in the three\,--\,body problem.

Throughout the section we shall consider a fixed two\,--\,body
subsystem of the three\,--\,body system concerned. Therefore, its
index will be omitted in the notations. Statements will be given for
the first variant of potentials~(\ref{vpot}).  If
necessary,  assertions corresponding to the second variant
(\ref{vpotb}) will be written in parentheses.

According to Eq.~(\ref{LScht2}) the kernel $t(k,k',z)$
of the pair T\,--\,matrix satisfies the integral equation
\be
\label{LScht2int}
t(k,k',z) \EQ v(k,k')\,-\,\int_{\Rt} dq \,\Frac{v(k,q)\hs{0.02cm}t(q,k',z)}
{q^2-z}\,.
\ee
All the dependence of $t(k,k',z)$ on $z$ is determined by the
integral term of the equation. The latter, with respect to the
variable $z$, represents a Cauchy type integral.
Integrals of this kind appear in the Faddeev equations
(\ref{Fadin}) as well and they are all of the form
\be
\label{Phi}
\Phi(z) \EQ
\int_{{\R}^{N}} dq\,\Frac{f(q)}{\lambda+q^2-z}
\ee
where $N=3$ or $N=6$ and $\lambda\leq 0$.

Let us describe some properties of the function  $\Phi(z)$, supposing
that $f$ is a holomorphic function of the variable
$q\in {\C}^N$.  We also assume that $f$  satisfies the
estimate $|f(q)|\leq  c\, M(q)$ where  $M(q)>0 $ and
$
\int_{S^{N-1}} d\hat{q} \, M(q) \leq
c\,\frac {\exp\hs{0.02cm}(a|\Img q|)} {(1+|q|)^{\theta} }
$
for some $c>0$ and $\theta\in (N-2,N-1)$.

Let $\Re_{\Phi}$ be the Riemann surface of the function
$$
\zeta(z) \EQ \left\{\!\!
\begin{array}{ll}
(z-\lambda)^{1/2}\,, &  N\,\, {\rm odd}\,,\\
\ln \hs{0.02cm}(z-\lambda)\,,   &  N\,\, {\rm even}\,.
\end{array} \right.
$$
The surface $\Re_{\Phi}$ arises as a result of pasting the
sheets $\Pi_l$ representing copies of the complex plane ${\C}$
cut along the ray $[\lambda,+\infty)$.  The sheet $\Pi_l$ is
identified with a branch of the function $\zeta(z)$. If
$\zeta(z)=(z-\lambda)^{1/2}$, we suppose for the main branch
that $l=0$, otherwise $l=1$.  If $N$ is even, then as  $l$ we
take a number of the function $\ln \hs{0.02cm}(z-\lambda)$ branch,
$
\ln\hs{0.02cm}(z-\lambda) =\ln |z-\lambda |+ i \hs{0.02cm}\varphi
\,+ \,i \hs{0.02cm}2\pi l
$
where $\varphi={\rm arg}\hs{0.02cm} (z-\lambda)$, so that
$z-\lambda=|z-\lambda|\exp\hs{0.02cm}(i\hs{0.02cm}\varphi)$ and
$\varphi\in[0,2\pi)$.
\begin{lemma}\label{LPhi}
The function $\Phi(z)$ is holomorphic in the complex plane
${\C}$ cut along the ray  $[\lambda,+\infty)$ and admits the
analytic continuation on the Riemann surface $\Re_{\Phi}$  given by
\be
\label{Phil}
\reduction{\Phi(z)}{\Pi_l} \EQ
\Phi(z)\,-\,\pi \hs{0.02cm}i \hs{0.02cm}l \big(\sqrt{z-\lambda}\,
\big)^{N-2}
\Int_{S^{N-1}} d\hat{q}f\big(\sqrt{z-\lambda}\hs{0.02cm}\hat{q}\big)
\ee
and the estimates
\begin{eqnarray*}
& \| \Phi(z) \|\LE c \,\| f \|_{\theta}(1+|z|)^{-\nu '},& \\ [0.1cm]
&\Int_{S^{N-1}} d\hat{q} f\big(\sqrt{z-\lambda}\hs{0.02cm}\hat{q}\big) \LE
c \,\| f \|_{\theta}(1+|z|)^{-\theta/2}
\exp\!\big(a\big|\Img\sqrt{z-\lambda}\,\big|\,\big)&
\end{eqnarray*}
with $\|f\|_{\theta} =$ $\mathop{\rm
sup}_{q\in{\C}^N}M^{-1}(q)|f(q)|$ hold for any $\nu'<({\theta-(N-2)})/{2}$.
\end{lemma}

\vs{0.3cm}For a proof see \cite{MotTMF}.

\vs{0.3cm}Using the relations (\ref{Phil}) one can easily obtain
representations for the analytic continuations of the kernels
$r_0(k,k',z)$ and $R_0(P,P',z)$. The Riemann surface $\Re^{(2)}$
of the kernel $r_0(z)$  coincides with $\Re_\Phi$ corresponding
to $\zeta(z)=z^{1/2}$, since, in this case, $N$ is odd ($N=3$).  For
$R_0(z)$, the Riemann surface is logarithmic, $\zeta(z)=\ln\hs{0.02cm}(z)$
($N$ even, $N=6$).  The continuation of the kernels $r_0(k,k',z)$
and $R_0(P,P',z)$ in $z$ is understood in the sense of
generalized functions (distributions) over $\Othree$ and
$\Osix$, respectively.  In the example with $r_0(z)$, we consider
the continuation of the bilinear form
$
\Phi(z)=(r_0(z)f_1,f_2)\equiv
\int_{\Rt} dq\,\frac{f_1(q)f_2(q)}{q^2-z},
\; f_1, f_2\in \Othree.
$
By Lemma~\ref{LPhi} one gets immediately
\begin{eqnarray}
\label{r02}
\reduction{r_0(z)}{\Pi_1} \!\! & = & \!\!
r_0(z)\,+\,{\rm a}_0(z)\hs{0.02cm}\rjt(z)
\hs{0.02cm}\rj(z)\,,
\quad \hs{0.7cm}{\rm a}_0(z)\EQ -\,\pi \hs{0.02cm}i \hs{0.02cm}\sqrt{z}\,,
\\ [0.2cm]
\reduction{R_0(z)}{\Pi_l}\!\! & = & \!\! R_0(z)\,+\,A_0(z)\hs{0.02cm}l\,
\rJot(z)\hs{0.02cm}\rJ_0(z)\,,
\quad A_0(z) \EQ -\,\pi \hs{0.02cm}i \hs{0.02cm}z^2,\\ [0.1cm]
&& \hs{4.5cm} l\EQ \pm 1,\, \pm 2,\,\ldots \, .  \nonumber \label{R06}
\end{eqnarray}
Using  Lemma~\ref{LPhi} one can immerse the
equation (\ref{LScht2int}) in a suitable Banach space and show
then that for the first variant of potentials (\ref{vpot}),
the integral operator $vr_0(z)$ can be continued  analytic  in
$z$ as a compact operator on the whole  Riemann surface
$\Re^{(2)}$.  Regarding the second variant of the potentials
(\ref{vpotb}), the existence of such a continuation into
the unphysical sheet is guaranteed only for the domain
$\Pi_1\cap\cP_b$,
\be
\label{Pib}
\cP_b \EQ \left\{ z:\,\,\,
\Real z>-b^2+\Frac{1}{4b^2}\,(\Img z)^2 \right\}
\ee
bounded by the parabola $\Img\sqrt{z}=b$ inside of which the
function  $v\big(\sqrt{z}\,\big(\hk-\hk'\big)\big)$ is holomorphic in
$z$ for arbitrary $\hk,\,\hk'\in S^2$.  According to
Eq.~(\ref{r02}) the product $\reduction{vr_0(z)}{\Pi_1}$ may be
written as
$
\reduction{vr_0(z)}{\Pi_1}=vr_0(z)+
{\rm a}_0(z)v\hs{0.02cm}\rjt(z)\hs{0.02cm}\rj(z)
$.
Thus, the continued equation (\ref{LScht2int}) for
$t'(z)=\reduction{t(z)}{\Pi_1}$ has the form
\be
\label{tprime}
t'(z) \EQ v-v\hs{0.02cm}r_0(z)\hs{0.02cm}t'(z)\,-\,{\rm a}_0(z)\hs{0.02cm}v\, \rjt(z)\,
\rj(z)\, t'(z)\,.
\ee
Let us transfer the term  $vr_0(z)t'(z)$ to the left\,--\,hand side
of Eq.~(\ref{tprime}) and invert (at $z\not\in\sigma_d(h)$)
the operator $[I+vr_0(z)]^{-1}$ using the relation
$[I+vr_0(z)]^{-1}v=t(z)$.  This leads to the expression
\be
\label{t2in}
t'(z) \EQ t(z)\,-\,{\rm a}_0(z)\hs{0.02cm}t(z)\hs{0.02cm}\rjt(z)\hs{0.02cm}
\rj(z)t'(z)
\ee
for $t'(z)$ in terms of the half\,--\,on\,--\,shell kernel
$t'\big(\sqrt{z}\,\hat{k},k',z\big)$, the first argument of which belongs
to the energy\,--\,shell (\ref{EnergyS-twobody}). Now to find
$\rj\hs{0.02cm}{t'}$ we apply $j$ to the both parts of Eq.~(\ref{t2in})
and then transfer  a term including $\rj(z)\hs{0.02cm}t(z)$ to the
left\,--\,hand side. Then we get
\vs{0.2cm}\be
\label{sjt2}
s(z)\hs{0.02cm}\rj(z)\hs{0.02cm}t'(z) \EQ \rj(z)\hs{0.02cm}t(z)
\ee
where $s(z)$, $z\in\Pi_0$,
is the two\,--\,body scattering matrix (cf., e.\,g., formula~(7.70)
of \cite{MF}),
\vs{0.2cm}\be
\label{s2}
s(z) \EQ \hat{I} \,+ \,{\rm a}_0 (z) \hs{0.02cm}\rj(z)\hs{0.02cm}t(z)
\hs{0.02cm}\rjt(z)\,,
\ee
$s(z): L_2\big(S^2\big)\rightarrow L_2\big(S^2\big)$,
with ${\rm a}_0(z) = -\pi i \sqrt{z}$
and  $\hat{I}$ is the identity operator in  $L_2 \big(S^2\big)$.

The absolute term $(\rj \hs{0.03cm}t)\big(\hat{k},\hat{k'},z\big)$ of
Eq.~(\ref{sjt2}), considered as a function of the first argument
$\hat{k}\in S^2$, is an element of  $L_2\big(S^2\big)$ at $z\not\in
\sigma_d(h)$.  The operator $\rj\hs{0.02cm} t\hs{0.02cm}\rjt$ on the right in
(\ref{s2})
is a compact operator in $L_2\big(S^2\big)$ for $z\not\in\sigma_d(h)$.
Therefore, on the domain of analyticity $\Pi_0\setminus
\sigma_d(h)$
$\bigl(\,\cP_b\cap\Pi_0\setminus\sigma_d(h)\,\bigl)$ of the
operator\,--\,valued  function $(\rj\hs{0.03cm} t\hs{0.03cm} \rjt)(z)$, one
can apply
to Eq.~(\ref{sjt2}) the Fredholm analytic
alternative~\cite{ReedSimonI} (see~\cite{MotTMF}). This means
that, with exception of a countable set $\sigma_{\rm res}$ having
no limit points in ${\C}\setminus\overline{\sigma(h)}$
$\big(\cP_b\setminus\overline{\sigma(h)}\,\big)$, the
 inverse operator $[s(z)]^{-1}$ exists and $\rj\hs{0.02cm} t'(z)=[s(z)]^{-1}
\rj\hs{0.02cm}t(z)$.  As a result we come to the following statement which is
in fact a one\,--\,channel variant of Theorem~2 of \cite{MotTMF}.

\begin{theorem}\label{Tht2body}
The two\,--\,body T\,--\,matrix $t(z)$ admits  analytic continuation in
the variable  $z$ on \vs{0.05cm}the sheet $\Pi_1$ $($on the domain
$\cP_b\cap\Pi_1)$ as a bounded operator in $L_2\big(\Rt\big)$. The result of
the continuation $\reduction{t(z)}{\Pi_1}\!
\Big(\!\reduction{t(z)}{\cP_b\cap\Pi_1}\!\Big)$ is expressed by
T\,-- and S\,--\,matrices taken in the physical sheet as
\vs{0.2cm}\be
\label{3tnf}
\reduction{t(z)}{\Pi_1} \EQ t(z)\,-\,{\rm a}_0(z)\hs{0.02cm} \tau(z)
\ee
where $\tau(z)=\big(t\hs{0.02cm}\rjt s^{-1}\rj \hs{0.02cm}t\big)(z)$.
The kernel $\reduction{t(k,k',z)}{\Pi_1}$
is a holomorphic function of the variables
$k,k'\in \Ct $ and
$z\in\Pi_1\setminus\bigl(\sigma_{\rm res}\cup\sigma_d(h)\bigr)$
$\bigl( k$ and $k'$ belonging to $W_b$ and $z\in\cP_b\cap\Pi_1\setminus
\bigl(\sigma_{\rm res}\cup\sigma_d(h) \bigr)\bigr)$.
\end{theorem}

\vs{0.3cm}On the physical sheet $\Pi_0$, the pair  T\,--\,matrix $t(z)$ admits
the representation (\ref{tdistr}) including the formfactors
$\phi_j$, $j=1,2,\,\ldots, n$. It follows from the Lippmann\,--\,Schwinger
equation for $\phi_j$ that
\be
\label{phij}
\phi_j(k) \EQ -\,\int_\Rt dq\, v(k,q)\,\Frac{1}{q^2-\lambda_j}\,\phi_j(q)\,,
\quad \lambda_j \LT 0\,,
\ee
that the formfactor $\phi_j(k)$ admits an analytic continuation
in $k$ on  $\Ct$ (on $W_{2b}$)  and that, at the same time, it
satisfies the type (\ref{vpot}) estimate where one has to
replace $\theta_0$ by a number $\theta$, $1 < \theta <
\theta_0 $, which can be taken arbitrarily  close to
$\theta_0$~\cite{Faddeev63}.  Hence the eigenfunction
\be
\label{psiphij}
\psi_j(k) \EQ -\,\Frac{\phi_j(k)}{k^2-\lambda_j}
\ee
of $h$ admits also such an analytic continuation on $\Ct$ (on
$W_{2b}$) with the exception of the set
$\{k\in\Ct:\,k^2=\lambda_j\}$ where $\psi_j(k)$ has
singularities $\big($turning for $k=\sqrt{z}\,\hat{k}$, $\hat{k}\in
S^2$, into a pole in $z$ at $z=\lambda_j$).

The regular summand $\tilde{t}(k,k',z)$ of the kernel of $t(z)$
is a holomorphic function in the variables $k,k'\in\Ct$, $z\in\Pi_0$
 ($k,k'\in W_b$, $z\in \cP_b\cap\Pi_0$). It admits the
estimate
$$
|\tilde{t}(k,k',z)| \LT {c}{\,(1+|k-k'|)^{-\theta} }\cdot
\exp\hs{0.02cm}[a(|\!\Img k|+|\!\Img k'|)]
$$
with arbitrary $\theta\in(1,\theta_0)$.

Regarding the operator $\reduction{t(z)}{\Pi_1}$, it follows
from Eq.~(\ref{3tnf}) that the points $z\!\in\!\sigma_d(h)$ become,
 generally speaking, poles of the first order of this operator.
One can easily check however that if the eigenvalue
$\lambda\in\sigma_d(h)$ is simple, then the respective
singularities of  both summands of  (\ref{3tnf}) compensate
each other and the there is no pole of $\reduction{t(z)}{\Pi_1}$
 at  $z=\lambda$.  It follows from the Fredholm analytic
alternative~\cite{ReedSimonI} for Eq.~(\ref{sjt2}) that the
poles of $\reduction{t(z)}{\Pi_1}$ at $z\in \sigma_{\rm res}$
are of a finite order.  It is also easy to show
that if $\cA(\hk)$ is a nontrivial solution of the equation
\be
\label{zero2}
s(z){\cal A} \EQ 0
\ee
at some $z\in\sigma_{\rm res}$, $z\not\in\sigma_d(h)$, then the
Schr\"odinger equation has at this $z$ a nontrivial (resonance)
solution $\psi^{\#}_{\rm res}$. The asymptotics of such a solution
written in configuration space is exponentially increasing,
                     $$
\psi^{\#}_{\rm res}(x) \; \Equal{x\rightarrow\infty} \;
\left({\cal A}(-\hat{x})+o(1)\right)
\Frac{{\rm e}^{-i\sqrt{z}\,|x|}}{\,|x|}\,, \quad
x\in{\R}^3\,.
                     $$
The function $\psi^{\#}_{\rm res}(x)$ is a so\,--\,called Gamow
vector corresponding to the resonance at the energy $z$ (see
e.\,g. \cite{Newton}, \cite{BohmJMP}, \cite{BohmQM}). The
function  ${\cal A}\big(\hat{k}\big)$ gives  sense to the breakup
amplitude of the resonance state.

The formula for the analytic continuation of the scattering matrix
into the unphysical sheet  $\Pi_1$ (on the set $\cP_b\cap\Pi_1$)
follows immediately from Eq.~(\ref{3tnf}) (see
\cite{MotTMF}),
\be
\label{S2P1}
\reduction{s(z)}{\Pi_1} \EQ {\cal E}\hs{0.02cm}[s(z)]^{-1}{\cal E}\,,
\ee
where ${\cal E}$ stands for the inversion in $L_2\big(S^2\big)$, $({\cal
E}f)\big(\hat{k}\big)=f\big(-\hat{k}\big)$.

Utilizing the representation (\ref{3tnf}) one can easily get an
explicit representation in terms of the physical sheet as well
for the analytic continuation on $\Pi_1$ $\left(\mbox{on }
\cP_b\cap\Pi_1 \right)$ of the resolvent $r(z)$:
\be
\label{res2}
\reduction{r(z)}{\Pi_l} \EQ
r\,+\,{\rm a}_0(I-rv)\hs{0.03cm}\rjt s^{-1} \rj (I-vr)\,.
\ee
The continuation is understood again in the sense of generalized
functions (distributions) over  $\Othree$.  This means that one
has to continue the bilinear form
$$
\Phi(z) \EQ \bigl(r(z)f_1,f_2\bigr) \; \equiv \;
\int_{\Rt} dk \int_{\Rt} dk'\hs{0.03cm}
f_2(k)\hs{0.02cm} r(k,k',z) \hs{0.02cm}f_1(k')\,,
\quad f_1,\,f_2\in\Othree\,.
$$

\newsection{The matrix $M(z)$ and the three\,--\,body scattering
            matrices in the physical sheet}
\label{SMSphys}

\vs{-0.8cm}\newsubsection{Faddeev equations}
\label{SFaddReal}
At the beginning, we recall briefly some principal
properties~\cite{Faddeev63}, \cite{MF} of the Faddeev equations
(\ref{MFE}) for the matrix $M(z)$  and  of the kernels
$M\abd(P,P',z)$ themselves for real arguments $P,P'\in\Rs $.  To
formulate these properties we reproduce here the following
definition from ~\cite{Faddeev63}, \S\hs{0.02cm}5.

An operator\,--\,valued function
$\cQ\abd(z):
\cH_0\rightarrow\cH_0,$  $\an,\bn=1,2,3,$ of the variable
$z\in{\C}$, is said to be a function of type $\cD\abd$ if it
admits the representation
\begin{eqnarray}
\nonumber
 \cQ\abd(z)\!\! & = & \!\!\cF\abd(z)\,+\,\Phi\ad\bg\ad(z)\cI\abd(z)\,+\,
 \cJ\abd(z)\bg\bd(z)\Phis\bd \label{QDrepr} \\ [-0.25cm]
 & &  \mbox{  } \\ [-0.25cm]
 & & \!\! + \;\;     \Phi\ad\bg\ad(z)\cK\abd(z)\bg\bd(z)\Phis\bd\,.
 \nonumber
\end{eqnarray}
The operator\,--\,valued functions
$\cF\abd(z)$ and $\cI\abd(z)$ from $\cH_0$ into
 $\cH_0$ and $\cH\au,$ respectively, and
$\cJ\abd(z): \cH\bu\rightarrow\cH_0$
and
$\cK\abd(z): \cH\bu\rightarrow\cH\au$
are called components of the function $\cQ\abd(z)$.
If $\cQ\abd(z)$ is an integral operator, then its kernel is
of the type $\cD\abd$.

Let
$$
\cN(P,\theta) \EQ \Sum_{\an,\bn,\atop \an\!\ne\!\bn} \,
(1+|p\ad|)^{-\theta} (1+|p\bd|)^{-\theta}\,.
$$
A function $\cQ(z)$ of type  $\cD\abd$
is said to be  a function of
class $ \cD_{\an\bn}(\theta,\mu)$ if its
components $\cF\abd$, $\cI\abd$, $\cJ\abd$ and $\cK\abd$
are integral operators and for the kernels  $\cF\abd(P,P',z)$
at $P,\,P', \,\Delta P,\, \Delta P'\in\Rs$,
the estimates
\begin{eqnarray}
\label{Est1}
& & \!\!\left| \cF\abd(P,P',z) \right| \LE \,c\, \cN(P,\theta)
\big(1+{p\bd'}^2 \big)^{-1}\,, \\ [0.2cm]
& & \!\! \left|\cF\abd(P+\Delta P,P'+\Delta P',z+\Delta z)-
\cF(P,P',z)\right|    \nonumber \\ [-0.25cm]
& & \mbox{ } \\ [-0.25cm]
 &  \leq & \!\!  c\,\cN(P,\theta) \big(1+{p\bd'}^2 \big)^{-1}
(|\Delta P|^\mu+|\Delta P'|^\mu+|\Delta z|^\mu) \nonumber \label{Est2}
\end{eqnarray}
are valid for a certain $c>0$ and, at the same time the kernels
$\cI_{\an,j;\bn}\big(p\ad,P',z\big),$ 
$\cJ_{\an;\bn,k}\big(P,p'\bd,z\big)$
and $\cK_{\an,j;\bn,k}\big(p\ad,p'\bd,z\big)$
satisfy the inequalities obtained
from (\ref{Est1}) and (\ref{Est2})
by taking, respectively, $k\ad=0$,\, $k'\bd=0$
or simultaneously $k\ad=0$ and $k'\bd=0$.

Let $\cQ\un(z)$ be the iteration (\ref{Qitert})
of the absolute term of Eq.~(\ref{MFE}).
In contrast to $\cQ^{(0)}(z)=\bt(z)$, the
kernels of the operators $\cQ\un(z)$, beginning already with  $n=1$,
do not include $\delta$\,--\,functions. Moreover, it follows from the
representation (\ref{tcd}) for $\bt\ad(z)$ explicitly manifesting
a contribution of the discrete spectrum of the pair subsystems, that the
matrix elements $\cQ\abd\un(z)$,
$\an,\bn=1,2,3,$ of the operators $\cQ\un(z)$ for $n\geq 1$
are actually functions of type $\cD\abd$.
Their components $\cF\abd\un(z)$, $\cI\abd\un(z)$,
$\cJ\abd\un(z)$ and $\cK\abd\un(z)$
at $z\in{\C}\setminus[\lambda_{\rm min},+\infty)$
are bounded operators depending on $z$ analytically.
In the case of the potentials (\ref{vpot}) and (\ref{vpotb})
at $z\not\in[\lambda_{\rm min},+\infty)$,
the H\"older index $\mu$ of smoothness of their kernels with respect to the
variables $P,P',p\ad$ and $p'\bd$
is equal to 1. If $n\leq 3$,
then, as
$\Img z\rightarrow 0$  with  $\Real z\in[\lambda_{\rm min},+\infty)$,
the kernels $\cF_{\an\bn}\un$,
$\cI_{\an,j;\bn}\un,$
$\cJ_{\an;\bn,k}\un$,
and $\cK_{\an,j;\bn,k}\un$
have the so\,--\,called  minor (three\,--\,particle) singularities
(see \S\hs{0.02cm}5 of ~\cite{Faddeev63}
and \S\hs{0.02cm}2 of Chapter~III of \cite{MF})
which weaken with increasing $n$.
For $n\geq 4$ such singularities do not appear at all,
and these kernels become H\"older functions
in all their variables including
the limit values  $z=E\pm i\hs{0.02cm}0$,  $E\in(\lambda_{\rm min},+\infty).$
A more precise statement~\cite{Faddeev63} is the following:
The operator\,--\,valued functions $\cQ\abd\un(z)$ for $n\geq 4$,
are of type $\cD\abd(\theta,\mu)$ for $0<\theta<\theta_0$, $0<\mu<{1}/{8}$
uniformly with respect to $z$
in any bounded set of the complex plane
${\C}$ cut along the ray $[\lambda_{\rm min},+\infty)$.
One can take  $\theta$, $\theta<\theta_0$, arbitrarily
close to $\theta_0$.
Thus, it is convenient~\cite{Faddeev63}
to choose instead of $M(z)$ the new unknown
$
\cW (z)=$ $ M(z)- $ $ \sum_{n=0}^{3} \cQ\un(z)
$
satisfying the equation
\vs{0.2cm}\be
\label{MFEW}
\cW (z) \EQ \cW \uo(z)\,-\,\bt(z)\hs{0.02cm}\bRo(z)\hs{0.02cm}\Y
\hs{0.02cm}\cW (z)
\ee
analogous to Eq.~(\ref{MFE}) but with another absolute term
$\cW \uo(z)=\cQ^{(4)}(z)$.

The solvability of Eq.~(\ref{MFEW}) is established in the Banach space
$\cB_{\theta\mu}$ whose elements are aggregates
$w=(\rho_1,\rho_2,\rho_3,\sigma_{1,1},\,\ldots,\sigma_{1,n_1},
\sigma_{2,1},\,\ldots,\sigma_{2,n_2},
\sigma_{3,1}, \,\ldots,\sigma_{3,n_3})^\dagger$
consisting of the functions
$\rho\ad(P),$ $P\in\Rs,$  and
$\sigma\ajd(p\ad),$  $p\ad\in\Rt$,  $\alpha = 1\,2,3, $ $j=1,2,\,\ldots,
n_{\alpha}$. The norm  $\|w\|_{\theta\mu}$ is defined in $\cB_{\theta\mu}$ by
\begin{eqnarray*}
\|  w \|_{\theta\mu} \!\! & = & \!\! \Sum_{\an=1}^{3}
\,{\rm sup} \!\left\{ \Frac{1}{\cN(P,\theta)}
\left[
\rho\ad(P)+  \Frac{|\rho\ad(P+\Delta P)-\rho\ad(P)|}{|\Delta P|^\mu}
\right] \right.       \\
 & & \left. \hs{1.3cm}+ \;\;(1+|p\ad|)^{2\theta}\Sum_{j=1}^{n\ad}
\left[         \sigma\ajd(p\ad)+
\Frac{\sigma\ajd(p\ad+\Delta p\ad)-\sigma\ajd(p\ad)}{|\Delta p\ad|^\mu}
\right] \right\}.
\end{eqnarray*}
The operator $-\bt(z)\bRo(z)\Y$ in (\ref{MFEW}) corresponds to the
operator $\bA(z)$ in $\cB_{\theta\mu}$ the action $w'=\bA(z)w$ of which
is defined according to the representation (\ref{tdistr}) in such a
way that
\begin{eqnarray*}
\rho'\!\!   &=& \!\! -\tilde{\bt}\hs{0.02cm}\bRo\hs{0.02cm}\Y
(\rho+\Phi\bg\sigma)\,,      \\[0.1cm]
\sigma'\!\!  &=& \!\!  \Phi\hs{0.02cm}\bg\hs{0.02cm}\Phis\hs{0.02cm}
\bRo\hs{0.02cm}\Y(\rho +\Phi\bg\sigma)
\end{eqnarray*}
where $\bg(z)=\diag\hs{0.02cm}\{\bg_1(z),\bg_2(z),\bg_3(z)\}$, and
$\rho, \,\sigma$ and $\rho',$ $\sigma'$ stand for the components of the
elements $w$ and $w'$, e.\,g. $\rho=(\rho_1,\rho_2,\rho_3)^\dagger,$
$\sigma=(\sigma_{1,1}, $
$\ldots, \sigma_{1,n_1}, \sigma_{2,1},\,\ldots,
\sigma_{2,n_2},\sigma_{3,1},$ \linebreak
$\ldots,\sigma_{3,n_3})^\dagger$.

Eq.~(\ref{MFEW}) is replaced then with the following equations in
$\cB_{\theta\mu}$
\begin{eqnarray}
\label{weq1}
w\bd(P',z)\! \! &=& \!\! w\bd\uo(P',z)\,+\,\bA(z) \hs{0.02cm}w\bd(P',z)\,,\quad
\;\;\bn \EQ 1,\,2,\,3\,,          \\ [0.2cm]
\label{weq2}
w\bkd(p\bd',z)\!\! &=& \!\!w\bkd\uo(p\bd',z)\,+\,\bA(z)\hs{0.02cm}
w\bkd(p\bd',z),\quad
\bn\EQ 1 ,\,2,\,3\,,\\ [0.1cm]
  & &  \hs{4.0cm} k \EQ 1,\,2,\ldots,\,n\bd\,, \nonumber
\end{eqnarray}
where $w\bd\uo(P',z)$  stands for an element of
$\cB_{\theta\mu}$ consisting of the kernels 
\begin{eqnarray*}
\rho\abd\uo(P,P',z)\!\! & = & \!\!\cF\abd\uo(P,P',z)\,, \quad
\hs{0.4cm} (\an \EQ 1,\,2,\,3)\,, \\ [0.1cm]
\sigma\uo_{\an,j;\bn}(p\ad,P',z) \!\! & = & \!\!
\cJ\uo_{\an,j;\bn}(p\ad,P',z)
\quad (\anum)\,,
\end{eqnarray*}
of the components $\cF\uo(z)$ and $\cI\uo(z)$ of the
function $\cW \uo(z)$ being considered at fixed $\bn$, $P'$ and $z$.
Analogously, $w\uo\bkd(p\bd',z)$ is an element of $\cB_{\theta\mu}$
consisting of the kernels
\begin{eqnarray*}
\rho\uo_{\an,j;\bn}(p\ad,P',z) \!\! & = & \!\!
\cJ\uo_{\an,j;\bn}(p\ad,P',z)\, \quad
\;\; \; (\an \EQ 1,\,2,\,3)\,, \\ [0.1cm]
\sigma\uo_{\an,j;\bn,k}(p\ad,p\bd',z) \!\! & = & \!\!
\cK\uo_{\an,j;\bn,k}(p\ad,p\bd',z) \quad (\anum)\,,
\end{eqnarray*}
of components
$\cJ\uo(z)$ and $\cK\uo(z)$ of the function $\cW \uo(z)$, which are
considered at fixed $\bn,$ $p\bd'$ and $z$.

The properties of Eqs.~(\ref{weq1}) and~(\ref{weq2}) are described
by Theorems~7.1 and~7.2 of the {\sc L.\,D. Faddeev} paper~\cite{Faddeev63}.
We combine these theorems in the following statement.

\begin{theorem}\label{ThFaddeev}
For all $z$ belonging to the complex plane ${\C}$
cut along the ray $[\lambda_{\rm min},+\infty)$, the operator $\bA(z)$
as well as all its powers $\bA^n(z)$, $n=2,3,\,\ldots,$
are defined in $\cB_{\theta\mu'}$ for
$3/2 < \theta < \theta_0$, $0 <\mu < {1}/{8}$
on a dense subset consisting of the elements $w\in\cB_{\theta\mu'}$,
$\mu' > \mu$. If $n\geq 5$ then  $\bA^n(z)$
admits a continuation over all $\cB_{\theta\mu}$
as a compact operator. A set of values of $z$
where the homogeneous equation $w=\bA(z) w$
has a nontrivial solution coincides,
with the possible exception of
its limit points, with the discrete spectrum $\sigma_d(H)$
of the Hamiltonian  $H$. Therefore, the Fredholm alternative may be
applied to Eqs.~{\rm (\ref{weq1})}
and~{\rm (\ref{weq2})} and thereby
these equations are solvable uniquely
in $\cB_{\theta\mu}$ for any  $z\not\in\overline{\sigma_d(H)}$
including the points $z=E\pm i\hs{0.02cm}0$,
$E\in(\lambda_{\rm min},+\infty)\setminus\overline{\sigma_d(H)}$.
\end{theorem}

\vs{0.3cm}Using Theorem \ref{ThFaddeev} one can establish as well
the properties of the operator $M(z)$ itself.
The corresponding statement
was presented by {\sc L.\,D. Faddeev} in \cite{Faddeev63},
Theorem~5.1. In our notations it is the following:

\begin{theorem}\label{ThFaddeev1}
{Eq.~{\rm (\ref{MFE})} is uniquely solvable at
$z\not\in\overline{\sigma_d(H)}$. Its solution $M(z)$ admits
the representation
\be
\label{MQW}
   M(z) \EQ \Sum_{n=0}^{3}\,\cQ\un(z) \,+\,\cW (z)\,,
\ee
where the operator\,--\,valued function $\cW (z)$ is holomorphic
in the variable $z$ at $z\not\in\overline{\sigma(H)}$
and its components $\cW\abd(z)$
belong to the classes $\cD\abd(\theta,\mu)$,
 $3/2 < \theta < \theta_0$, $0 < \mu < {1}/{8}$, \,\,
uniformly with respect to $z$ varying in  arbitrary bounded set
of the complex plane
${\C}$ cut along the ray $[\lambda_{\rm min},+\infty)$
and with removed neighbourhoods of the points of $\sigma_d(H)$.}
\end{theorem}
\newsubsection{Scattering matrices}
\label{SSmatrReal}
Let us begin now by recalling the structure of the
three\,--\,body scattering operator $\bS$ (see e.\,g. the
book~\cite{MF}, \S\hs{0.02cm}6 of Chapter\,I and \S\hs{0.02cm}3 of
Chapter III).
We introduce for \vs{0.05cm}this purpose the operator\,--\,valued function
$\cT(z)$, $\cT(z):\cH_0\oplus\cH_1\rightarrow\cH_0\oplus\cH_1$,
for $z\in{\C}\setminus\overline{\sigma(H)}$,
\be
\label{T3body}
\cT(z) \;\equiv \;\left(\!\!
\begin{array} {cc}
\Omega M(z)\Omega^{\dagger}   &    \Omega M(z)\Y\Psi  \\ [0.1cm]
\Psi^{*}\Y M(z)\Omega^{\dagger} &   \Psi^{*}(\Y\bv +\Y M(z)\Y)\Psi
\end{array}\!\!\right).
\ee
In the following the notation $\bv=\diag\hs{0.02cm}\{ v_1, v_2, v_3 \}$
will be used.
Note that $\cT_{00}(z)=\Om M(z)\Omt\!\equiv\! T(z)$ and 
$\cT_{00}(z)\!:\!\cH_0\rightarrow \!\cH_0$.
The remaining components $\cT_{01}(z)\!:\!\cH_1\!\rightarrow\!\cH_0$,\,\,
$\cT_{10}(z):\cH_0\rightarrow\cH_1$ and
$\cT_{11}(z):\cH_1\rightarrow\cH_1$ are expressed by
the three\,--\,body transition operators~\cite{MF}
(see also~\cite{Schmid})
$U_0(z)=\Om M(z)\Y$, $U_0^\dagger=\Y M(z)\Omt$
and $U(z)=\Y\bv +\Y M(z) \Y$:
$\cT_{01}=U_0\Psi$, $\cT_{10}=\Psis U_0^\dagger$
and $\cT_{11}=\Psis U\Psi.$
The operator $\cT(z)$ is a matrix integral operator
with the kernels
$\cT_{00}(P,P',z)$,  $\cT_{\an,i;\, 0}(p\ad,P',z)$,
$\cT_{0;\, \bn,j}(P,p'\bd,z)$ and $\cT_{\an,i;\,
 \bn,j}(p\ad,p'\bd,z)$ ($\an=1,2,3,$  $i=1,2,\,\ldots,n\ad,$
 $\bn = 1,2,3$, $j = 1,2,\ldots, n_{\beta}$)
the properties of which are determined for complex $z$ including
the limit points  $z=E\pm i\hs{0.02cm}0$, $E > \lambda_{\rm min}$
by Theorem~\ref{ThFaddeev1}.

By $\what{\cT}(z)$, $\what{\cT}(z):$
$\what{\cH}_0\oplus\what{\cH}_1\rightarrow\what{\cH}_0\oplus\what{\cH}_1$
we denote the analytic continuation in ${\C}^\pm$
(see Theorems~\ref{ThLTL}, \ref{ThJ0TJ0t} and~\ref{ThJ0MYPsiJ1t})
of a matrix operator\,--\,valued function
 of the variable $z$ whose components
have the following kernels at $z=E\pm i\hs{0.02cm}0$
$$
\begin{array}{rcl}
\bigl(\what{\cT}(E\pm i\hs{0.02cm}0) \bigr)_{00}
\big(\what{P},\what{P}'\big)\!\! & = & \!\!
\cT_{00}\big(\pm\sqrt{E}\hat{P},\,
\pm\sqrt{E}\hat{P}',E\pm i\hs{0.02cm}0)\,, \quad \hs{1.6cm}E \GT 0\,;
\\ [0.2cm]
\bigl(\what{\cT}(E\pm i0) \bigr)_{0;\,\bn,j}\big(\what{P},\what{p}'\bd\big)
\!\!& = & \!\!
\cT_{0;\,\bn,j}\big(\pm\sqrt{E}\what{P},\,
\pm\sqrt{E-\lambda\bjd}\,\hat{p}'\bd,E\pm i\hs{0.02cm}0\big)\,,
\quad E \GT 0\,; \\ [0.2cm]
\bigl(\what{\cT}(E\pm i\hs{0.02cm}0) \bigr)_{\an,i;0}\big(\hat{p}\ad,
\hat{P}'\big) \!\!& = & \!\!
\cT_{\an,i;0}\big(\pm\sqrt{E-\lambda_{\an,i}}\,\hat{p}\ad,\,
\pm\sqrt{E}\what{P}',E\pm i\hs{0.02cm}0\big)\,, \quad
E \GT 0 \,; \\ [0.2cm]
\bigl(\what{\cT}(E\pm i\hs{0.02cm}0) \bigr)_{\an,i;\,\bn,j}\big(\hat{p}\ad,
\hat{p}'\bd\big)
\!\!& = & \!\!\cT_{\an,i;\,\bn,j}\big(\pm\sqrt{E-\lambda_{\an,i}}\,
\hat{p}\ad,\,
\pm\sqrt{E-\lambda\bjd}\,\hat{p}'\bd,E\pm i\hs{0.02cm}0\big)\,,   \\ [0.2cm]
&  & \phantom{E > \Max \{ \lambda_{\an,i},\, \lambda\bjd\}}
\hs{1.1cm}    E \GT \max \hs{0.04cm} \{ \lambda_{\an,i},\, \lambda\bjd\}\,.
\end{array}
$$
We assume by definition that for $z=E \pm i\,0$
the product $\bigl(\bJ\cT\bJt\bigr)(z)$
coincides with $\what{\cT}(z)$,
\be
\label{mcteq}
(\bJ\cT\bJt)(z)\;=\;\left(\!\!
\begin{array}{cc}
\bigl(\rJo \cT_{00}\rJot \bigr)(z)  &
\bigl(\rJo \cT_{01}\rJt_1 \bigr)(z) \\ [0.15cm]
\bigl(\rJ_1 \cT_{10}\rJot \bigr)(z)  &
\bigl(\rJ_1 \cT_{11}\rJt_1 \bigr)(z)
\end{array}\!\!\right)
\; \equiv \; \what{\cT}(z)\,.
\ee
The elements of the matrix $(\bJ\cT\bJt)(z)$ are expressed in terms of
amplitudes of different processes taking place in the three\,--\,body
system under consideration (see Sec.~\ref{SNumerMethod}).

The three\,--\,body scattering operator $\bS $ is unitary
in the space $\cH_0 \oplus \cH_1 $ and has as well as $\cT$,
a natural block structure.  The kernels of its components
$\bS_{00}$, $\bS_{0;\beta,j}$, $\bS_{\alpha,i;0}$,
$\bS_{\alpha,i;\beta,j}$ read, repectively,
\begin{eqnarray}
\label{S00}
\bS_{00}(P,P') \!\!& = &\!\!
\delta(P-P')\,-\, 2\pi i\, \delta\big(P^2-P'^2\big) \cT_{00}
\big(P,P', P'^2+i\hs{0.02cm}0\big)\,,  \\ [0.2cm]
\label{S0b}
\bS_{0;\beta,j}\big(P,p'\bd\big)\!\! & = & \!\!
-\,2 \pi i\, \delta\big(P^2-p'^{2}\bd-\lambda\bjd\big)
\cT_{0;\beta,j}\big(P,p'\bd,\lambda\bjd+p'^{2}\bd+i\hs{0.02cm}0\big),
\\ [0.2cm]
\label{Sa0}
\bS_{\alpha,i;0}(p\ad, P') \!\! & = & \!\!
- \,2 \pi i\, \delta\big(\lambda_{\alpha,i}+p^{2}\ad - P'^2\big)
\cT_{\alpha,i;0}\big(p\ad, P',P'^2 +i\hs{0.02cm}0\big)\,,   \\ [0.2cm]
\label{Sab}
\phantom{MMM}\bS_{\alpha,i;\beta,j}(p\ad,p'\bd) \!\!& = & \!\!
\delta\abd\delta_{ij}\delta\big(p\ad-p'\bd\big)
\,-\,2 \pi i\,\delta\big(\lambda_{\alpha,i} +p\ad^2 -\lambda_{\beta,j} -
p\bd^2\big) \nonumber \\ [-0.25cm]
& & \!\! \mbox{} \\ [-0.25cm]
& & \!\! \times \;\;\cT_{\alpha,i;\beta,j}\big(p\ad,p'\bd,\lambda_{\beta,j} +
p\bd^2 +
i\hs{0.02cm}0\big)\,.
\nonumber
\end{eqnarray}

The scattering matrices arise from  $\bS$ in the spectral
decomposition of $H$ as operators acting in the ``cross
section'' (at fixed energy) of the space ${\cal H}_0\oplus {\cal
H}_1$ in the von Neumann direct integral~\cite{MerkDiss}.  As a
matter of fact, the extraction of the scattering matrix from
$\bS$ corresponds to the replacements $|P|^2\rightarrow E$,
$\lambda_{\alpha,i} +p\ad^2 \rightarrow E$ ($\an=1,2,3,$
$i=1,2,\,\ldots,n\ad$) in the expressions (\ref{S00}) -- (\ref{Sab})
and then to the  factorization of the dependence of the kernels
of $\bS$ on the energies $E$ and $E'$,
\be
\label{SEE}
\bS(E,E') \EQ -\,\pi \hs{0.02cm}i\hs{0.02cm} \delta(E-E')\hs{0.02cm}
\tilde{\vartheta}(E)\hs{0.02cm}S'(E+i\hs{0.02cm}0)
\hs{0.02cm}\tilde{\vartheta}(E')
\ee
where $\tilde{\vartheta}(E)$ is a diagonal matrix\,--\,function
constructed from the Heaviside unit step functions
${\vartheta}(E)$
and ${\vartheta}(E-\lambda\bjd)$:
$\tilde{\vartheta}(E)=$
$\diag\hs{0.03cm}\{\vartheta(E),\vartheta(E-\lambda_{1,1}),\,\ldots,
\vartheta(E-\lambda_{1,n_1}),$ \linebreak
$\vartheta(E-\lambda_{2,1}),\,\ldots,
\vartheta(E-\lambda_{2,n_2}),\vartheta(E-\lambda_{3,1}),\,\ldots,
\vartheta(E-\lambda_{3,n_3}) \}. $
At $z\in{\C}$ we understand by $S'(z)$
the operator\,--\,valued function defined by
$
S'(z)=A^{-1}(z)\, \hbI \,+\, \what{\cT}(z).
$
Hereafter,  $A(z)=\diag\hs{0.03cm}\{A_0(z), A_1(z)\}$ with
$A_0(z)=-\pi i z^2$
and
$A_1(z)=\diag\hs{0.03cm}\big\{ A^{(1)},A^{(2)},A^{(3)} \big\}$ where
             $$
A\au (z) \EQ \diag\hs{0.03cm}\{ A_{\an,1}(z),\,\ldots,A_{\an,n\ad}(z) \}\quad
\mbox{with} \quad A\ajd(z) \EQ-\,\pi i \,\sqrt{z-\lambda\ajd}\,.
             $$

Continuing the factorization,
$
S'(z)=S(z)A^{-1}(z)=A^{-1}(z)S^{\dagger}(z),
$
corresponding to separate in (\ref{SEE}) the multiplier
$-\pi i A^{-1}(E+i\hs{0.02cm}0)$
as a derivative of a measure in the von~Neumann integral
above~\cite{MerkDiss} for  $\cH_0 \oplus \cH_1 $, one arrives at
the scattering matrices
\be
\label{SMx}
S(z) \EQ \hbI \,+\, \bigl(\bJ\cT\bJt \, A \bigr)(z) \quad\mbox{and} \quad
\St(z)\EQ \hbI \,+\, \bigl(A\,\bJ\cT\bJt \bigr)(z) \,.
\ee
In  contrast to \cite{MerkDiss}, it is more convenient for
us to use precisely this nonsymmetrical form of the scattering
matrices.  The matrices $S(z)$ and $\St(z)$ are considered as
operators in $\what{{\cal H}}_0\oplus \what{{\cal H}}_1 $.  At
$z=E+i\hs{0.02cm}0$, $E>0$ these operators are unitary. At
$z=E+i\hs{0.02cm}0$,
$E<0$ there are certain truncations of $S(z)$ and $\St(z)$
determined by a number of open channels which are unitary in
$\what{{\cal H}}_0\oplus \what{{\cal H}}_1 $; namely, the matrices
$\tilde{S}(E)=
\hat{\bf I}+ \tilde{\vartheta}(E)\bigl(S(E+i\hs{0.02cm}0)-\hat{\bf I}\bigr)
\tilde{\vartheta}(E)$
and
$\tilde{\St}(E)=
\hat{\bf I}+ \tilde{\vartheta}(E)\bigl(\St(E+i\hs{0.02cm}0)-\hat{\bf I}\bigr)
\tilde{\vartheta}(E)$.
It follows from Eq.~(\ref{SMx}) that the operator $\cT$ may be
considered as a kind of a ``multichannel T\,--\,matrix'' (cf.
\cite{MotTMF}) for the system of three particles.

It should be noted
that the matrix $\cT(z)$ may be replaced in Eq.~(\ref{SMx})
with the matrix $\cT^\dagger(z)$ obtained from
$\cT(z)$ by the substitution $\Y\bv\rightarrow\bv\Y$
(respectively, \linebreak
 $U \rightarrow U^{\dagger}=\bv\Y+ \Y M \Y$)
in the second component of the lower row of (\ref{T3body}).
To prove the equality
$\bigl(\bJ\cT^\dagger\bJt \bigr)(z)=\bigl(\bJ\cT\bJt \bigr)(z)$
it suffices to observe that for
$z=E\pm i\hs{0.02cm}0$, $E>\lambda\ajd$ ($\alpha = 1,2,3,\;j =
1,2,\,\ldots,n_{\alpha}$)
\be
\label{JPvPJ}
\big(\rJ_1 \Psis\Y\bv\Psi\rJt_1\big)(z) \;= \;
\big(\rJ_1 \Psis\bv\Y\Psi\rJt_1\big)(z)\,.
\ee
Indeed, according to Eqs.~(\ref{phij})
and~(\ref{psiphij}),
\begin{small}
\be \hs{1.0cm}
\label{PYvP}
(\Psis\Y\bv\Psi)_{\an,i;\bn,j}(p\ad,p'\bd) \;= \;
-\Frac{1-\delta\abd}{|s\abd|^{3}} \,\cdot\,
\Frac{    \overline{\phi}_{\an,i}\big(\tk\ad\bu(p\ad,p'\bd)\big)\,
\phi\bjd\big(\tk\bd\au(p'\bd,p\ad)\big)    }
{ \big[\tk\ad\bu(p\ad,p'\bd)\big]^2-\lambda_{\an,i}   }\,,
\ee
\be \hs{1.0cm}
\label{PvYP}
(\Psis\bv\Y\Psi)_{\an,i;\bn,j}(p\ad,p'\bd) \;=\;
-\Frac{1-\delta\bad}{|s\abd|^{3}} \,\cdot\,
\Frac{    \overline{\phi}_{\an,i}\big(\tk\ad\bu(p\ad,p'\bd)\big)\,
\phi\bjd\big(\tk\bd\au(p'\bd,p\ad)\big)    }
{ \big[\tk\bd\au(p'\bd,p\ad)\big]^2-\lambda\bjd   }
\ee
\end{small}
where
\be
\label{ktilde}
\tk_{\gamma}^{(\delta)}(q,q')) \EQ
\Frac{-c_{\gamma\delta} q+q'}{s_{\gamma\delta} }\,, \quad
\gamma,\,\delta \EQ 1,\,2,\,3\,, \quad  q,\,q'\in\Rt \,,
\ee
$\big($we shall suppose later that $q,q'\in \Ct\big)$.
One can see easily that the denominators
of the fractions (\ref{PYvP})
and (\ref{PvYP})
coincide on the energy shells
$|p\ad|=\sqrt{E-\lambda_{\an,i}},$
$|p'\bd|=$ \linebreak $ = \sqrt{E-\lambda\bjd},$
$E>\lambda_{\an,i},$
$E>\lambda\bjd$:
\begin{eqnarray}
\label{k2lambda}
 \big(\tk\ad\bu\big)^2\,-\,\lambda_{\an,i}\!\! & = & \!\!
 \big(\tk\bd\au\big)^2\,-\, \lambda\bjd  \nonumber \\ [0.1cm]
& = & \!\! \Frac{1}{|s\abd|^2}\big(E-\lambda_{\an,i}+E-\lambda\bjd \\ [0.1cm]
 & & \!\! \hs{1.3cm}- \;
2c\abd\sqrt{E-\lambda_{\an,i}}\sqrt{E-\lambda\bjd}
\big(\hp\ad,\hp'\bd\big)-s\abd^2 E\big)\,.
\nonumber
\end{eqnarray}
Meanwhile, the expression  (\ref{k2lambda})
cannot become zero at $E>\lambda_{\an,i},$ $E>\lambda\bjd$
(see Lemma~\ref{LEqQuadr}).
It follows now from Eqs.~(\ref{PYvP}), (\ref{PvYP}) and~(\ref{k2lambda})
that the equality~(\ref{JPvPJ}) is true.

Along with  $S(z)$ and $\St(z)$  we shall consider further also
the truncated scattering matrices
\be
\label{Slcut}
 S_l(z) \;\equiv \;\hat{\bf I} \,+\, \big(\tL\bJ \cT \bJt L A\big)(z)
 \quad \mbox{and}\quad
\St_l(z) \;\equiv \; \hat{\bf I} \,+\, \big(A L\bJ \cT \bJt\tL\big)(z)\,,
\ee
where the multi\,--\,index
\be
\label{lmulti}
l\;=\;(l_0,l_{1,1},\,\ldots,l_{1,n_1},l_{2,1},\,\ldots,l_{2,n_2},l_{3,1},
\,\ldots,l_{3,n_3})
\ee
has components $l_0=0$ or  $l_0=\pm 1$
and $l\ajd=0$ or $l\ajd=1$, $\alpha = 1,2,3, \; j=1,2,\,\ldots,n_{\alpha}$.
By $L$ and $\tL$ in (\ref{Slcut}) and in the following we mean
the diagonal number matrices
\be
\label{Lmatrix}
L \; = \;\diag\hs{0.03cm}\{ l_0,l_{1,1},\,\ldots,l_{1,n_1},l_{2,1},\,\ldots,
l_{2,n_2},l_{3,1},\,\ldots,l_{3,n_3}    \}
\ee
and
\vs{0.2cm}\be
\label{Ltmatrix}
\tL \; = \;\diag\hs{0.03cm}\{ |l_0|,l_{1,1},\,\ldots,l_{1,n_1},l_{2,1},
\,\ldots, %
l_{2,n_2},l_{3,1},\,\ldots,l_{3,n_3}    \}\,,
\ee
corresponding to the multi\,--\,index $l$.
The matrix $\tL$ is evidently a projection in $\what{\cH}_0
\oplus \what{\cH}_1$ on the subspace $\what{\cH}_1^{(l)}$  if
$l_0=0$ or on the subspace $\what{\cH}_0 \oplus
\what{\cH}_1^{(l)}$ if $l_0\neq 0$.  Here $\what{\cH}_1^{(l)}
=\bigoplus_{l\ajd\neq 0} \what{\cH}\aju$ in both cases.

As can be seen from formulas (\ref{SMx}) and (\ref{T3body}) the
scattering matrices $S(z)$ and $\St(z)$ include the kernels
$M\abd(P,P',z)$ taken on the energy shells: their arguments $P\in\Rs$
and $P'\in \Rs$ are connected with the energy $z=E+i\hs{0.03cm}0$ by
Eq.~(\ref{EnergyS3}) at $E > 0$ or Eqs.~(\ref{EnergyS2}) at $E >
\lambda\ajd$. We establish below [see formula~(\ref{Ml3fin})] that
the analytic continuation of the matrix  $M(z)$ into unphysical sheets
of the energy $z$ is expressed in terms of the analytic continuation of
the truncated scattering matrices $S_l(z)$ or $\St_l(z)$ and
half\,--\,on\,--\,shell Faddeev components $M\abd(z)$ taken in the physical
sheet.  More precisely, along with  $S_l(z)$, the final
formula~(\ref{Ml3fin}) includes the matrices $\bigl( L_0 \bJo M
\bigr)(z)$, $\bigl( L_1\rJ_1 \Psis\Y M\bigr)(z)$ and $\bigl( M
\bJot L_0  \bigr)(z)$, $\bigl( M\Y\Psi\rJt_1 L_1 \bigr)(z)$.
Here, $l$ is
the multi\,--\,index (\ref{lmulti}) and $L=\diag\hs{0.03cm}\{L_0,L_1\}$,
the respective matrix (\ref{Lmatrix}) with  $L_0=l_0$
and
$L_1=\diag\hs{0.03cm}\{l_{1,1},\,\ldots,l_{1,n_1},l_{2,1},\,\ldots,
l_{2,n_2},l_{3,1},\,\ldots,l_{3,n_3}\}$.

Further, we  formulate some statements
(Theorems~\ref{ThLTL} -- \ref{ThJ0MYPsiJ1t}) concerning the existence
of an analytic continuation of the above matrices and their domains
of holomorphy.  Proofs of these statements will be based on
analysis~\cite{Faddeev63} of the Faddeev equations~(\ref{MFE}).  For
this, one must pay a special attention to the study of
the domains of holomorphy in $z$ of the functions
\be
\label{Znam}
\left[ p\ad^2+{p\bd'}^2-2c\abd(p\ad,p'\bd)-s\abd^2 z \right]^{-1}
\ee
with one or both arguments $p\ad$ and $p'\bd$ situated on the
energy shells (\ref{EnergyS3}) or (\ref{EnergyS2}).  The
functions (\ref{Znam}) arise when iterating Eq.~(\ref{MFE})
because of the presence of the multiplier $\bR_0$ in the operator
$-\bt\bR_0\Y$. Also, the functions (\ref{Znam})
display the singularities (\ref{psiphij}) of the
eigenfunctions
$\psi\ajd,$ $\alpha = 1,2,3, \; j=1,2,\,\ldots, n_{\alpha}$.

In the case where the arguments
$p\ad$ and/or $p'\bd$ are taken on the shells
(\ref{EnergyS2}), \,\, $p\ad=$ $= \sqrt{z-\lambda\aid}\,\hat{p}\ad$ and
$p'\bd=\sqrt{z-\lambda\bjd}\,\hat{p}'\ad$,
the holomorphy domains of the functions (\ref{Znam})
with respect to the variable $z$ are described by the following
simple lemmas.
\begin{lemma}\label{LEqParab}
For any $\rho\geq 0$,\,\, $-1\leq \eta\leq 1$, the domain
\be
\label{DomainParab}
  \Real z \GT \Frac{\lambda}{c^2}\,+\,\Frac{c^2}{4s^2|\lambda|}\,(\Img z)^2
\ee
contains no roots $z$ of the equation
\be
\label{quadr0}
 z\,-\,\lambda \,+\, \rho \,+\, 2c\,\sqrt{z-\lambda}\,\sqrt{\rho}\:\eta
  \,-\,s^2 z \EQ 0\,,
\ee
with $\lambda < 0,$\,\, $0 < |c| < 1$ and $s^2=1-c^2.$
For any number $z\in{\C}$ outside of
the domain {\rm (\ref{DomainParab})}
one can always find  values of the parameters
$\rho\geq 0$  and  $\eta,$ $-1\leq \eta \leq 1,$
such that the left\,--\,hand side of Eq.~{\rm (\ref{quadr0})}
becomes equal to zero at the point $z$.
\end{lemma}
\begin{lemma}\label{LEqQuadr}
Consider the equation
\be
\label{quadr}
 z\,-\,\lambda_1 \,+ \,z \,-\, \lambda_2 \,+\,
 2c\,\sqrt{z-\lambda_1}\,\sqrt{z-\lambda_2}\:\eta
  \,-\,s^2 z \EQ 0
\ee
where $\eta\in[-1,\, 1]$,
 $\lambda_1 \leq \lambda_2 < 0,$ $0 < c <1 $ and $s^2 =1 -c^2$.
Then the following assertions hold:

\rm 1) \it If $|\lambda_2|>c^2 |\lambda_1|$,
then for all $\eta\in[-1,1]$
Eq.~{\rm (\ref{quadr})} has a unique root $z$ and this root is real.
Moreover, $z=z_{+}$ if  $\eta\geq 0$,
and $z=z_{-}$ if $\eta\leq 0$ with
\begin{small}
\be
\label{koren}
  z_{\pm}\;=\;\Frac{
      \big(1+c^2-2c^2\eta^2\big)(\lambda_1+\lambda_2) \pm
     2\sqrt{
           c^2\eta^2[\lambda_1 \lambda_2\, s^4-
        (\lambda_2-\lambda_1)^2 c^2 (1-\eta^2)]     }  }
   {(1+c^2)^2 - 4 c^2 \eta^2}\,.
\ee
\end{small}
If $\eta$ runs through the interval $[-1,1]$, the roots $z_\pm$ fill
the interval $ [z_{\rm lt}, z_{\rm rt}] $ whose ends are
\be
\label{zl}
 z_{\rm lt} \;=\;\Frac{1}{s^2} \big[-|\lambda_1|-|\lambda_2|-
 2c\sqrt{|\lambda_1| \cdot |\lambda_2|}\,\big]
\ee
and
\be
\label{zr}
 z_{\rm rt}\;=\;\Frac{1}{s^2} \,\big[-|\lambda_1|-|\lambda_2|+
 2c\sqrt{|\lambda_1| \cdot |\lambda_2|}\,\big]\,, \quad z_{\rm rt}\LT
 \lambda_1\,.
\ee

 \rm 2) \it If $|\lambda_2|=  c^2 |\lambda_1|$,
then Eq.~{\rm (\ref{quadr})} has two real roots:

\hs{0.5cm}\rm a) \it the root $z=\lambda_1$ existing for all $\eta\in [-1,1];$

\hs{0.5cm}\rm b) \it the root $ z=z_{-}$
given by~{\rm (\ref{koren})} which exists for  $-1\leq \eta\leq 0$ only.

For $z_{\rm lt} \EQ -\,|\lambda_1|\left(1+{2c^4}/{s^2}\right)$,
$-1\leq\eta\leq 1$ these roots together fill the interval
$[z_{\rm lt},\lambda_1]$.

\rm 3) \it If $|\lambda_2| <  c^2 |\lambda_1|$, then

\hs{0.5cm}\rm a) \it for $-1\leq\eta\leq\eta^{*},$
$
    \eta^{*}=({\sqrt{c^2-\rho}\,\sqrt{1-c^2\rho}})/({c(1-\rho)})
$
and $\rho={|\lambda_2|}/{|\lambda_1|},$
Eq.~{\rm (\ref{quadr})} has two real roots $z_{\pm}$
given by {\rm (\ref{koren})} which fill the interval
$[z_{\rm lt}, z_{\rm rt}]$
with ends {\rm (\ref{zl})} and {\rm (\ref{zr})},
$z_{\rm rt}<\lambda_1;$

\hs{0.5cm}\rm b) \it for $\eta^{*} < \eta \leq 0 $  Eq.~{\rm (\ref{quadr})}
has two complex roots $z_{\pm}$ described
again by Eq.~{\rm (\ref{koren})}.
If $\eta$ varies, these roots fill an ellipse centered at the point
$$
   z_c \EQ -\,|\lambda_1|\!\left[ 1 +
  \Frac{(c^2-\rho)^2}{s^2(1+c^2)(1+\rho)} \right].
$$
The half\,--\,axes of this ellipse are given by
$$
 a \EQ |\lambda_1| \cdot \Frac{(c^2-\rho)(1-c^2\rho)}
 {  (1+c^2)s^2(1+\rho)  }
$$
{\rm(}along the real axis{\rm)} and
$$
 b \EQ |\lambda_1| \cdot \Frac{(c^2-\rho)(1-c^2\rho)}
 {   (1+c^2)s^2(1-\rho)\sqrt{(1+c^2)^2 - 4c^2\eta^{*2}}  }
$$
{\rm(}along the imaginary axis{\rm)}. The right vertex of the ellipse
is located at the point
        $$
z^{({\rm e})}_{\rm rt} \EQ z_c+a \EQ -\,\Frac{|\lambda_1|+|\lambda_2|}{1+c^2}
$$
situated between  $\lambda_1$ and $\lambda_2$. Its left vertex is
              $$
z^{({\rm e})}_{\rm lt} \EQ z_c-a \LT z_{\rm rt}\,.
              $$
\end{lemma}\rm
\begin{lemma}\label{LEq23}
Let the parameters of the equation
\be
\label{QEq23}
z\,-\,\lambda \,+\, z\nu \,+\,
2c\,\sqrt{z}\,\sqrt{z-\lambda}\,\sqrt{\nu}\eta \,- \,s^2 z \EQ 0
\ee
satisfy the conditions $\nu\in[0,1]$, $\eta\in[-1,1]$, $\lambda<0$,
$c\in(0,1)$ and $ s^2 =1-c^2$. Then, if
$\nu$ and $\eta$ run through the above ranges, the roots $z$
of Eq.~{\rm(\ref{QEq23})} fill the ray
$\left(-\infty, {\lambda}/({1+c^4})\right]$ and the circle centered
at the point
$z_c={\lambda}/({1-c^4})$ the
radius of which is equal to ${c^2\lambda}/({1-c^4})$.
\end{lemma}

\vs{0.3cm}Also, we shall use
\begin{lemma}\label{LEqQuadrPrime}
Let the parameters  $\nu'$ and $\eta$ of the equation
$$
\rho \,+\, z\nu' \,+\, 2c\sqrt{z}\,\sqrt{\nu'}\sqrt{\rho}\eta \,-\,
s^2 z \EQ 0
$$
run through the intervals
$0\leq\nu'\leq 1$ and $-1\leq\eta\leq 1$, respectively,
and $c > 0,$  $ s^2 =1-c^2$,  $ z\in{\C}$ be fixed. Then
the roots $\rho$ of this equation
fill a set consisting of the line segment $[0,z]$ in
the complex plane ${\C}$ and a circle centered at the origin,
the radius of which is equal to $c^2 |z|.$
\end{lemma}

\newsubsection{Analytic continuation of  the matrices
        $M\Y\Psi\rJt_1$, $\rJ_1\Psis\Y M$ and $\cT_{11}$ }
\label{S41}
Let  $\Pibj$ be a domain in the complex plane
${\C}$ cut along the ray $[\lambda_{\rm min},+\infty)$
where the conditions
(\ref{DomainParab}) with $\lambda=\lambda\bjd$, $c=c\abd$
and the inequalities
\be
\label{Neqsb}
\Real z \GT \lambda\bjd\,-\,s\abd^2 b^2 \,+\,\Frac{1}{4 s\abd^2 b^2}\,
(\Img z)^2
\ee
are valid simultaneously for all $\an=1,2,3,$ $\an\neq\bn$.
In the case of the potentials (\ref{vpot}) one has to take
$b=+\infty$ in (\ref{Neqsb}).

By $\cR_{\an,i;\,\bn,j}$, $\an\neq\bn$
we denote a domain complementary
in ${\C}\setminus[\lambda_{\rm min},+\infty)$ to the set
filled by the roots of Eq.~(\ref{quadr}) in the case where
$\lambda_1={\rm min} \hs{0.03cm}\{\lambda\aid,\lambda\bjd\}$,
$\lambda_2=$ $= {\rm max} \hs{0.03cm}\{\lambda\aid,\lambda\bjd\}$,
$c=|c\abd|$
and  $\eta=\big(\hp\ad,\hp'\bd\big)$
runs through the interval $[-1,1]$.
\begin{theorem}\label{ThLTL}
The matrix integral operator $L'_1\hat{\cT}_{11}(z)L''_1$,
$z=E\pm i\hs{0.03cm}0$
acting in $\what{\cH}_1$ admits  analytic continuation in $z$
from the edges of the ray $E\in(\lambda,+\infty)$,\,\,
$$
\lambda \EQ  \Max{\mbox{   \scriptsize
                  $\begin{array}{c}
            l'_{\gn,k}\neq 0, \\ [0.1cm]
            l''_{\gn,k}\neq 0
                        \end{array}$
                   }
	     }
\!\!\!\! \lambda_{\gn,k}\,,
$$
on the domain
\be
\label{Pil1hol}
\Pi^{({\rm hol})}_{l'l''} \EQ
\left[ \Bigcap_{               \mbox{\scriptsize
                        $\begin{array}{c}
                      l'\aid\neq 0   \\ [0.1cm]
                      l''\bjd\neq 0
                       \end{array}$  }
                }
\cR_{\an,i;\,\, \bn,j}                              \right]\,
\Bigcap \,
\left[ \Bigcap_{               \mbox{\scriptsize
                        $\begin{array}{c}
                               l'_{\gn,k}\neq 0,   \\ [0.1cm]
                               l''_{\gn,k}\neq 0
                       \end{array}$  }
                }
\Pi_b^{(\gn,k)}  \right]
\Bigg \backslash \overline{\sigma(H)}
\ee
where
\begin{eqnarray*}
 l'_1 \!\! & = & \!\! \diag\left( l'_0,l'_{1,1},\,\ldots,l'_{1,n_1},
 l'_{2,1}, \,\ldots,l'_{2,n_2}, l'_{3,1},\,\ldots,l'_{3,n_3}\right),\\ [0.1cm]
l''_1 \!\! & = & \!\! \diag\left( l''_0,l''_{1,1},\,\ldots,
 l''_{1,n_1}, l''_{2,1},\,\ldots, l''_{2,n_2},
 l''_{3,1},\,\ldots, l''_{3,n_3}\right)
\end{eqnarray*}
with  $l'_0=l''_0=0.$
The nontrivial kernels
$
\big( L'_1\hat{\cT}_{11}(z)L''_1 \big)_{\an,i;\,\, \bn,j}\!
(\hp\ad,\hp\bd',z)$,
 $ l'\aid\neq 0$, \linebreak
   $ l''\bjd\neq 0 $
turn out to be functions holomorphic in
$z\in\Pi^{({\rm hol})}_{l'l''}$
and real\,--\,analytic with respect to  $\hp\ad,\hp\bd'\in S^2$.
\end{theorem}
\begin{remark}\label{NotePiSym}
{\rm  The domains
$\Pi^{({\rm hol})}_{l'l''}$ and $\Pi^{({\rm hol})}_{l''l'}$ coincide,
$\Pi^{({\rm hol})}_{l'l''}=\Pi^{({\rm hol})}_{l''l'}$. }
\end{remark}

If $l'=l''=l$, we use for $\Pi^{({\rm hol})}_{l'l''}$
the notation $\Pi^{({\rm hol})}_{l}$, and we have
\be
\label{Pilh1}
\Pi^{({\rm hol})}_{l} \EQ \Pi^{({\rm hol})}_{ll}\,.
\ee
\begin{theorem}\label{ThMYPsiJt}
Let $L_0=l_0=0.$ Then the matrices
 $\big(  M\Y\Psi\rJt_1 L_1\big)(z)$
and \linebreak
 $\bigl( L_1 \rJ_1 \Psis \Y M \bigr)(z)$,
 $z=E\pm i\hs{0.03cm}0$
admit  analytic continuation in $z$ from the edges of the ray
$E\in(\lambda,+\infty)$,
                 $$
\lambda \EQ {\max_{(\bn,j):
\,\lambda_{\beta,j \neq 0}} {\lambda_{\beta, j}}}
                 $$
to the domain
$\Pilh\setminus\overline{\sigma(H)}$
as bounded for
$ z\not\in[\lambda_{\rm min}, +\infty)$
operator\,--\,valued functions of the variable $z$,
                   $$
\bigl(  M\Y\Psi\rJt_1 L_1\bigr)(z) \;:\;
\what{\cH}_1\;\longrightarrow \; \cG_0 \quad   and \quad
\bigl( L_1 \rJ_1 \Psis \Y M \bigr)(z) \;: \;
\cG_0\;\longrightarrow \; \what{\cH}_1\,.
                   $$
\end{theorem}

\vs{0.3cm}First we prove Theorem~\ref{ThMYPsiJt},  then
Theorem~\ref{ThLTL}.

\vs{0.3cm}\Proof of Theorem~\ref{ThMYPsiJt}.
We give the proof for the case of the matrix $M\Y\Psi\rJt_1 L_1$.
It follows from Theorem~\ref{ThFaddeev} that the kernels
$$
\begin{array}{l}
\big( M\Y\Psi\rJt_1 \big)_{\an;\bn,j}(P,\hp'\bd,E\pm i\hs{0.03cm}0) \;
\equiv \; \Sum_{\gn\neq\bn}\int_{\Rt} dk'\bd\,  M_{\an\gn}(P,P',
E\pm i\hs{0.03cm}0)\,\psi\bjd\big(k'\bd\big)\,, \\
\;\;\, \;
\hspace*{6.7cm}P'\EQ \big(k'\bd,\pm\,\sqrt{E-\lambda\bjd}\,\hp'\bd\big)
\end{array}
$$
of the nontrivial elements
$\big( M\Y\Psi\rJt_1 \big)_{\an;\bn,j}(E\pm i\hs{0.03cm}0)$,
 $\an=1,2,3,$ $l\bjd\neq 0$
are defined at  \linebreak
$E>\lambda\bjd$.
As a whole, the matrix $\bigl( M\Y\Psi\rJt L_1\bigr)(z)$
at $z=E\pm i\hs{0.03cm}0,$
           $$
E \;>\; \lambda\;=\;\Max{(\bn,j):\, l\bjd\neq 0} \lambda\bjd
           $$
satisfies the Faddeev equation(s)
\be
\label{FeMYPsi}
\bigl( M\Y\Psi\rJt_1 L_1\bigr)(z)\EQ \bigl( \bt\Y\Psi\rJt_1 L_1\bigr)(z)
\,- \, \bigl( \bt\bRo\Y M\Y\Psi\rJt_1 L_1\bigr)(z)\,,
\ee
the absolute term $\bt\Y\Psi\rJt_1 L_1$ of which
at $l\bjd\neq 0$ has the kernels
\begin{eqnarray}
\nonumber
& & \!\! \bigl( \bt\Y\Psi\rJt_1 \bigr)_{\an;\bn,j}(P,\hp'\bd,z) \\ [0.1cm]
\label{tYPsiABj}
 &  = & \!\! (1-\delta\abd)\cdot\Frac{1}{|s_{\an\bn}|^3}\cdot
t\ad\big(k\ad,\tk\ad\bu\big(p\ad,\sqrt{z-\lambda\bjd}\, \hp'\bd\big),z-
p\ad^2\big) \\ [0.1cm]
\nonumber
 & & \!\! \times\;\; \psi\bjd\big(\tk\bd\au\big(\sqrt{z-\lambda\bjd}\,
 \hp'\bd\big), p\ad\big)\big)\,.
\end{eqnarray}
Evidently, these  admit
analytic continuation in $z$ in the usual sense
as smooth (even real\,--\,analytic)
functions of the arguments
$P\in{\R}^6$\, and $\hp'\bd\in S^2$
on the domains where the respective two\,--\,body eigenfunctions
$$
  \psi\bjd\Big(\tk\bd\au\Big)=- \,
  \Frac{\phi\bjd\Big(\tk\bd\au\Big)}{{\tk\bd\au}{^2}-\lambda\bjd}
$$
are regular. The condition of regularity of the functions
$\psi\bjd\big(\tk\bd\au\big)$ is equivalent to the requirements
\be
\label{Neqlam}
\Big[\tk\bd\au\big(\sqrt{z-\lambda\bjd}\, \hp'\bd, \, p\ad\big)\Big]^2  \;
\neq \; \lambda\bjd
\quad \mbox{for all} \quad p\ad\in\Rt\;\;\mbox{and}\;\; \hp'\bd\in S^2.
\ee
As follows from Eqs.~(\ref{ktilde}), the requirements contrary
to the conditions (\ref{Neqlam}), are equivalent to the
conditions (\ref{quadr0}) with  $\lambda=\lambda\bjd,$
$\mu=|p\ad|^2,$  $c=|c\bad|$,   $s=|s\bad|$ and
$\eta=(\hp\ad,\hp'\bd)$.  Thereby, on the basis of
Lemma~\ref{LEqParab}, we conclude that the inequalities
(\ref{Neqlam}) at $l\bjd\neq 0$ are satisfied only 
if the inequalities (\ref{DomainParab}) with
$\lambda=\lambda\bjd$, $c=c\abd$ are valid.

In the case of the potentials (\ref{vpotb}),
one has to require additionally  the variables \linebreak
$\tk\ad\bu$ and  $\tk\bd\au$
to belong to the analyticity strips $W_b$ and $W_{2b}$
of the kernels \linebreak
$t\ad\big(k\ad,\tk\ad\bu,z-p\ad^2\big)$
and form factors $\phi\bjd\big(\tk\bd\au\big)$, respectively,  i.\,e.,
\be \hs{1.0cm}
\label{Neqb}
\left|\Img \tk\ad\bu\big(p\ad,\sqrt{z-\lambda\bjd}\,\hp'\bd\big) \right| \LT
b\,, \qquad
\left|\Img \tk\bd\au\big(\sqrt{z-\lambda\bjd}\,\hp'\bd,p\ad\big) \right| \LT
2b\,.
\ee
At $p\ad\in\Rt$, $\hp'\bd\in S^2$,
the inequalities~(\ref{Neqb}) are certainly satisfied if
the conditions (\ref{Neqsb}) are satisfied.

Let $\bigl( M\Y\Psi\rJt_1 L_1\bigr)\bjd(z)$ be the $(\bn,j)$\,--\,th column
of the matrix  $\bigl( M\Y\Psi\rJt_1 L_1\bigr)(z)$, \linebreak
$ l\bjd\neq 0$.
The kernels
$\bigl( M\Y\Psi\rJt_1 L_1\bigr)(z)_{\an;\,\bn,j}(P,\hp'\bd,z)$,
$\an=1,2,3$
of its components at  fixed  $\hp\bd$ satisfy
a respective system of the Faddeev
equations following from (\ref{FeMYPsi}).
As we established, the absolute terms (\ref{tYPsiABj})
of this system admit analytic continuation in $z$
from the edges of the ray  $(\lambda\bjd,+\infty)$
over the whole domain $\Pibj$. The same may be stated also regarding
the iterations
$
   \cQ\un\bjd(z)=\bigl(-\bt(z)\bRo(z)\Y\bigr)^n
  \bigl(  \bt\Y\Psi\rJt_1 L_1\bigr)\bjd(z).
$
Embedding the Faddeev equations for the column
$\bigl(  M\Y\Psi\rJt_1 L_1\bigr)\bjd(z)$
into the Banach space $\cB_{\theta\mu'}$
in the same way as it was done for Eq.~(\ref{MFE}), one
finds that the kernels
$\bigl(  M\Y\Psi\rJt_1 L_1\bigr)_{\an; \,\bn,j}(P,\hp'\ad,z)$
may be continued on $\Pibj\setminus\overline{\sigma(H)}$
as analytic functions of the variable $z$
for all $P\in\bR^6$,\, $\hp'\bd\in S^2$.
It follows from the estimates of the rate of
decrease of these kernels that
they represent the operator
$\bigl(  M\Y\Psi\rJt_1 L_1\bigr)\bjd(z)$:
$\what{\cH}^{(\bn,j)}\rightarrow\cG_0,$
bounded for
$z\not\in[\lambda_{\rm min},+\infty)\cup\overline{\sigma_d(H)}$
and depending analytically on
$z\in\Pibj\setminus\overline{\sigma(H)}$.
Therefore, we have proved the statement of the theorem for the
matrix $\bigl(  M\Y\Psi\rJt_1 L_1\bigr)(z)$.  The statement for
$\bigl( L_1 \rJ_1 \Psis \Y M \bigr)(z)$ can be established quite
similarly.

The proof of Theorem~\ref{ThMYPsiJt} is completed.%
{\nopagebreak\mbox{}\hfill $\Box$\par\addvspace{0.25cm}

\Proof of Theorem~\ref{ThLTL}. Note first that when proving
Theorem~\ref{ThMYPsiJt} we found out in passing that the kernels
$\bigl(\cQ\un\Y\Psi\rJt_1  \bigr)_{\an;\,\bn,j}(P,\hp\bd,E\pm i\hs{0.03cm}0)$,
$E >\lambda\bjd,$
of the operators $\cQ\un\Y\Psi\rJt_1 $
admit already for $n=0$ an immediate analytic continuation
on the domain of $z\in\Pibj$
as real\,--\,analytic functions of the variables $P\in\Rs$
and $\hp'\bd\in S^2$.
The same may be stated as well for the kernels
$(\rJ_1\Psis\Y \cQ\un)_{\bn,j;\,\an}(\hp\bd,P',E\pm i\hs{0.03cm}0),$
$E >\lambda\bjd,$
of the operators $\rJ_1\Psis\Y \cQ\un$
when continuing them in $\Pibj$.

On the basis of this note one may consider the relation
\begin{eqnarray}
\nonumber
  L'_1 \what{\cT}_{11} L''_1 \!\! & = & \!\!
L'_1\rJ_1\Psis\Y\bv\Psi\rJt_1 L''_1   \\ [-0.2cm]
& & \mbox{ } \\ [-0.2cm]
\label{LTL1}
 & & \!\! +\;\; L'_1\rJ_1\Psis\Y \left(
\bt-\bt\bRo\Y\bt+\bt\bRo\Y M \Y\bRo\bt
\right)\Y \Psi\rJt_1 L''_1\,, \nonumber
\end{eqnarray}
following from Eq.~(\ref{MFEmn}) at $m=n=0$ and
\begin{eqnarray*}
L'_1 \!\! & = & \!\! \diag\left\{ l'_{1,1},\,\ldots,l'_{1,n_1},
 l'_{2,1},\,\ldots,l'_{2,n_2},
l'_{3,1},\,\ldots,l'_{3,n_3}\right\}, \\ [0.1cm]
L''_1 \!\! & = & \!\! \diag\left\{ l''_{1,1},\,\ldots,l''_{1,n_1},
 l''_{2,1},\,\ldots,l''_{2,n_2},
 l''_{3,1},\,\ldots,l''_{3,n_3}\right\},
\end{eqnarray*}
as a representation for analytic continuation
of the matrix $L'_1 \what{\cT}_{11} L''_1 :
\what{\cH_1}\rightarrow\what{\cH_1}$
in terms of the operator  $M(z)$
with the kernels $M\abd(P,P',z)$, the arguments
$P$ and $P'$ of which are real, $P,P'\in\Rs.$
Additionally, one knows already that the summand
$ L'_1\rJ_1\Psis\Y \cQ \Y \Psi\rJt_1 L''_1$
with $\cQ=\bt\bRo\Y M\Y \bRo\bt$
admits analytic continuation in $z$ on the domain
                $$
\Bigcap_{ \mbox{ \scriptsize
$\begin{array}{r}
             l'\bjd\neq 0, \\ [0.1cm]
            l''\bjd\neq 0
\end{array}$   }
 } \Pibj\setminus\overline{\sigma(H)}
                 $$
as a matrix integral operator in $\what{\cH_1}$, since
the operators
$ L'_1\rJ_1\Psis\Y \bt$ and $\bt\Y\Psi\rJt_1 L''_1$ may be continued
on this domain.
Now, one only has to find a domain of holomorphy for the rest of
the summands in (\ref{LTL1}).

The nontrivial kernels
\be
\label{LPYvPL}
\mbox{\phantom{MMM}}
\bigl( L'_1\rJ_1\Psis\Y \bv \Psi\rJt_1 L''_1\bigr)_{\an,i;\,\bn,j}
(\hp\ad,\hp\bd,z)\,, \quad
l'\aid\neq 0\,,\,\, l''\bjd\neq 0\,,\,\,\, \hp\ad\in S^2,\,\,  \hp'\bd\in S^2,
\ee
of the first summand in (\ref{LTL1}) look like (\ref{PYvP}) where one has
to take $p\ad=\sqrt{z-\lambda_{\an,i} }\, \hp\ad$,\,\,
$p'\bd=\sqrt{z-\lambda\bjd }\, \hp'\bd$.
Furthermore, the functions
$\overline{\phi}_{\an,i}\big(\tk\ad\bu\big)$
have to be replaced with their analytic continuations
$\tilde{\phi}_{\an,i}\big(\tk\ad\bu\big)$
in a domain of complex $\tk\ad\bu$.
It is clear that the kernels (\ref{LPYvPL}) are
holomorphic in that domain of the variable $z$
where their denominators  (\ref{k2lambda}) can be equal to zero for no
$\hp\ad,\hp'\bd\in S^2$.
This domain is described by Lemma~\ref{LEqQuadr}.
It follows from Eq.~(\ref{PYvP}) in accordance with this lemma that the
kernel (\ref{LPYvPL}) is a holomorphic function of
$z$ in the domain $\cR_{\an,i;\,\bn,j}$.

In the case of the potentials
(\ref{vpotb}) we require additionally
the arguments $\tk\ad\bu(p\ad,p'\bd)$
and $\tk\bd\au(p\bd',p\ad)$
of the formfactors $\tilde{\phi}\aid$ and $\phi\bjd$
to belong at $p\ad=\sqrt{z-\lambda_{\an,i} }\, \hp\ad$ and
$p'\bd=\sqrt{z-\lambda\bjd }\, \hp'\bd$
to the holomorphy strips $W_{2b}$.
Note meanwhile that if  $\lambda_1 < \lambda_2$
then $\Img\sqrt{z-\lambda_1} \leq \Img\sqrt{z-\lambda_2}.$
Thus, the conditions
$\big|\Img \tk\ad\bu\big| < 2b$,\,\, $\big|\Img \tk\bd\au\big| < 2b$
are satisfied at
\be
\label{Neqvb}
\Real z \GT \lambda \,-\, \Frac{4 |s\abd|^2}{(1+|c\abd|)^2}\, b^2 \,+ \,
\Frac{(1+|c\abd|)^2}{16 |s\abd|^2 \, b^2}\,(\Img z)^2
\ee
where $\lambda=\Max{} \hs{0.03cm} \{ \lambda\aid,\lambda\bjd \}.$
Since $1+|c\abd|<2,$ the inequalities
(\ref{Neqvb}) are obeyed automatically if
$$
z\in \!\! \Bigcap_{ \mbox{   \scriptsize
                  $\begin{array}{r}
            l'_{\gn,k}\neq 0, \\ [0.1cm]
            l''_{\gn,k}\neq 0
                        \end{array}$
                     }
             }
\!\! \Pi^{(\gn,k)}_b\,.
$$

The two remaining terms,
$ L'_1\rJ_1\Psis\Y\bt\Y\Psi\rJt_1 L''_1$
and $ L'_1\rJ_1\Psis\Y\bt\bRo\Y\bt\Y\Psi\rJt_1 L''_1$,
have respectively the kernels
\be \label{LtL}
\begin{array}{l}
\vs{0.2cm}\Sum_{\gn\neq\an,\bn} \Frac{1}{|s_{\an\gn} s_{\bn\gn}|^3}\,\,
\int_\Rt dq\, \tilde{\psi}\aid\big(\tk\ad^{(\gn)}(p\ad,q)\big)\,
{\psi}\bjd\big(\tk\bd^{(\gn)}(p'\bd,q)\big)     \\
\phantom{\Sum_{\gn\neq\an,\bn} \Frac{1}{|s_{\an\gn} s_{\bn\gn}|^3}\,\,
\int_\Rt}\times\;\; t\gd\bigl( \tk\ad\au(q,p\ad),\tk\gd\bu(q,p'\bd),z-q^2
\bigr)
\end{array}
\ee
and
\begin{small}
\be \label{LttL}
\vs{0.2cm}\begin{array}{c}
\Sum_{\mbox{\scriptsize
$\begin{array}{c}
\gn\neq\an,\,\, \delta\neq\bn\\
\gn\neq\delta
\end{array}   $
} }
\,\Frac{1}{|s_{\an\gn}|^3 |s_{\gn\delta}| |s_{\bn\delta}|^3}\,\cdot\,
\int_\Rt\,\, dq \int_\Rt\,\, dq'\,
\tilde{\psi}\aid\big(\tk\ad^{(\gn)}(p\ad,q)\big)\,
{\psi}\bjd\big(\tk\bd^{(\delta)}(p'\bd,q)\big)     \\
\phantom{
\Sum_{\mbox{\scriptsize$
\begin{array}{c}
\gn\neq\an,\,\, \delta\neq\bn\\
\gn\neq\delta
\end{array}$}  }}\times \;\;\Frac{
t\gd\bigl( \tk\ad\au(q,p\ad),\tk\gd^{(\delta)}(q,q'),z-q^2 \bigr)\cdot
t_\delta \bigl( \tk_{\delta}^{(\gn)}(q,q'),
\tk\gd^{(\delta)}(q',p\bd'),z-{q'}^2 \bigr)
}
{q^2 +{q'}^2-2c_{\gn\delta}(q,q')-s_{\gn\delta}^2\, z}\,.
\end{array}
\ee
\end{small}
Here one has to take $p\ad=\sqrt{z-\lambda\aid}\, \hp\ad$,\,\,
and $p'\bd=\sqrt{z-\lambda\bjd}\, \hp'\bd$, as well as to assume that
the indices $(\an,i)$ and $(\bn,j)$ are such
that $ l'\aid\neq 0$ and $ l''\bjd\neq 0$.
The kernels above turn out to be holomorphic functions
of $z$ on a set where the conditions
$$
\left[ \tk\ad^{(\gn)}(\sqrt{z-\lambda\aid}\, \hp\ad,q)\right]^2
-\lambda\aid \;\neq \;0\,, \quad
\left[ \tk\bd^{(\gn)}(\sqrt{z-\lambda\bjd}\, \hp'\bd,q')\right]^2
-\lambda\bjd \;\neq \;0
$$
are satisfied for any  $\gn\neq\an,$ $\delta\neq\bn,$
$q,q'\in\Rt$
and $\hp\ad,\hp\bd'\in S^2$.
With respect to the variables $\hp\ad,\hp\bd'$,
the kernels  (\ref{LtL})
and (\ref{LttL}) are real\,--\,analytic with these conditions.
Besides, in the case of the potentials (\ref{vpotb}),
the arguments of the functions $\tilde{\psi}\aid$
and $\psi\bjd$ have to belong to the strip $W_{2b}$,
and the respective arguments of the T\,--\,matrices $t\gd$ and
$t_\delta$ to the strip $W_b$.
It is easily to check that the mentioned conditions
are satisfied for all the kernels (\ref{LtL}) and (\ref{LttL}) if
$$
z\in \!\!  \Bigcap_{ \mbox{   \scriptsize
                  $\begin{array}{r}
            l'_{\nu,k}\neq 0, \\ [0.1cm]
            l''_{\nu,k}\neq 0
                        \end{array}
			$
                     }
             }
\!\! \Pi^{(\nu,k)}_b\,.
                       $$
This completes the proof of Theorem~\ref{ThLTL}.%
{\nopagebreak\mbox{}\hfill $\Box$\par\addvspace{0.25cm}}

\newsubsection{Analytic continuation of the matrices $\bJo M$ and
$M\bJot$,  $\bJo M\bJot$  and  $\hat{\cT}_{00}$,
     $\bJo M\Y\Psi\rJt_1$ and $\rJ_1\Psis\Y M\bJot$  }
\label{S42}
Continuation of the half\,--\,on\,--\,shell matrices
$(\bJo M)(z),$ $\big(M\bJot\big)(z),$ $z=E\pm i\hs{0.03cm}0,$\,\,\, $E>0,$
into a domain of complex $z$ is considered in the sense of
distributions over $\Osix$.  For example, in the case of
$M\bJot$ we consider continuation of the bilinear form
$$
{\bigl(F,\, \big(M\bJot\big)(E\pm i\hs{0.03cm}0) \bigr) \;\equiv} \;\;
\Sum_{\an,\bn}\int_\Rs dP \int_{S^5}d\hP'
\, F\ad(P)M\abd\big(P,\pm\sqrt{E}\,\hP',E\pm i\hs{0.03cm}0\big)\,
f\bd\big(\hP'\big)
$$
where $F=(F_1,F_2,F_3)$ with $F\ad\in\Osix$ and  $f=(f_1,f_2,f_3)$ with
$f\ad\in\what{\cH}_0$.

When constructing continuation of this form and the form for
$(\bJo M)(E\pm i\hs{0.03cm}0)$, we rely on Lemmas~\ref{LEq23}
and~\ref{LEqQuadrPrime}
describing domains of holomorphy of the function (\ref{Znam})
in the case where the argument $P'$ belongs to the three\,--\,body
energy shell (\ref{EnergyS3}) and therefore
$p'\bd=\sqrt{z}\,\nu'\hp'\bd$ with
$\nu'\in[0,\, 1]$. Using the statements of these lemmas we introduce the
following definition.

Let $\wtilde{\Pi}_b^{(0)\pm}$,
$\wtilde{\Pi}_b^{(0)\pm}\!\subset\!{\C}^\pm $, be domains
complementary in ${\C}^\pm$ to the totality of
circles having radii $r={c\abd^2|\lambda\ajd|}/\big({1-c\abd^4}\big)$ and
centered at the points
$z_c={\lambda\ajd}/\big({1-c\abd^4}\big)$
where $\an,\bn=1,2,3,$\, $\bn\neq\an$, and $j=1,2,\,\ldots,n\ad$.

In the case of the potentials (\ref{vpotb})
the domain $\wtilde{\Pi}_b^{(0)\pm}$ has
to satisfy for both signs ``$+$'' and ``$-$''
the following extra conditions
\be
\label{CondbRes0}
\Real z \GT -\,\Frac{|s\abd|^2 b^2}{(1+|c\abd|)^2}
\,+\,\Frac{(1+|c\abd|)^2}{4\,|s\abd|^2 b^2}\, (\Img z)^2
\ee
for all $\an,\bn=1,2,3,$\, $\bn\neq\an$.

\begin{theorem}\label{ThMJ0tJ0M}
The kernels of the matrices  $\big(M\bJot\big)(z)$
and $(\bJo M)(z)$, \,\, $z=E\pm i\hs{0.03cm}0$, $E>0$,
admit analytic continuation in $z$
on the domains
$\wtilde{\Pi}_b^{(0)+}$ and $\wtilde{\Pi}_b^{(0)-}$,
$\wtilde{\Pi}_b^{(0)\pm}\subset{\C}^\pm$, respectively.
The continuation of the kernels of the matrices
$\big( \cQ^{(n)}\bJot  \big)\!(z)$
and  $\big( \bJo\cQ^{(n)}\big)\!(z)$ $(n\leq 3)$
included in the representation {\rm (\ref{MQW})} for  $M(z)$
has to be understood in the sense of distributions over  $\Osix$.
At the same time the kernels
\begin{eqnarray}
\label{WJ0t}
& & \cF\abd\big(P,\sqrt{z}\,\hP',z\big)\,,  \quad \hs{1.0cm}
\cI_{\an,j;\,\bn}\big(p\ad,\sqrt{z}\,\hP',z\big)\,, \nonumber  \\ [0.1cm]
 & & \cJ_{\an;\,\bn,k}\big(P,\sqrt{z}\,\sqrt{\nu'}\,\hp'\bd,z\big)\,, \quad
\cK_{\an,j;\,\bn,k}\big(p\ad,\sqrt{z}\,\sqrt{\nu'}\,\hp'\bd,z\big) \\ [0.1cm]
& & \an,\bn \EQ 1,\,2,\,3, \quad j\EQ 1,\,2,\,\ldots,\,n\ad\,,\quad
k \EQ 1,\,2,\,\ldots,n\bd\,, \nonumber
\end{eqnarray}
of the matrices $\bigl(\cQ\un\bJot\bigr)(z)$ $(n\geq 4)$
and $\big(\cW\bJot\big)(z)$ as well as the kernels
\be
\label{J0W}
\begin{array}{ll}
\cF\abd\big(\sqrt{z}\,\hP,P',z\big)\,, &
\cI_{\an,j;\,\bn}\big(\sqrt{z}\,\sqrt{\nu}\,\hp\ad,P',z\big)\,,  \\ [0.1cm]
\cJ_{\an;\,\bn,k}\big(\sqrt{z}\,\hP,p'\bd,z\big)\,,
&  \cK_{\an,j;\,\bn,k}\big(\sqrt{z}\,\sqrt{\nu}\,\hp\ad,p'\bd,z\big)
\end{array}
\ee
of the matrices
$\bigl(\bJo\cQ\un\bigr)(z)$ $(n\geq 4)$ and  $(\bJo\cW)(z)$
can be continued on the domains $\wtilde{\Pi}_b^{(0)\pm}$
as usual holomorphic functions of the variable $z$.
Being H\"older functions of the variables
$\hP'\in S^5$ or $\sqrt{\nu'}\,\hp'\bd$, $0\leq\nu'\leq 1$,
$\hp'\bd\in S^2$
$\big[\hP\in S^5$ or $\sqrt{\nu}\,\hp\ad,$
 $0\leq\nu\leq 1$,
$\hp\ad\in S^2\big]$ with index $\mu'\in(0,\,\,1/8)$,
the kernels  {\rm (\ref{WJ0t})}\,\,
{\rm [}the kernels {\rm (\ref{J0W})]}
considered as functions of $P\in\Rs,$ $p\ad\in\Rt$
$\big(P'\in\Rs$, $p'\bd\in\Rt\big)$, can be embedded in their totality in
$\cB_{\theta\mu}$ with $\theta$ and $\mu$ being arbitrary
numbers such that $\theta\in (3/2,\theta_0)$
and $\mu\in(0, 1/8)$.
For $|\Img z|\geq \delta > 0$
one can take $\mu=1$.
\end{theorem}
\begin{proof}
Let us use the equations  (\ref{MFEmn}) for $m,n\geq 4$ keeping in
mind that continuation of the kernels of
$\big(\bJo\cQ^{(m)}\big)(z)$ and
$\big(\cQ^{(n)}\bJot\big)\!(z)$ for $m,n\leq 3$ is realized in
the sense of distributions. At the same time, for $m,n\geq 4$ one
can attach to the products  $\big(\bJo\cQ^{(m)}\big)\!(z)$ and
$\big(\cQ^{(n)}\bJot\big)(z)$ an operator sense,
$\big(\bJo\cQ^{(m)}\big)\!(z):\cG_0\rightarrow\what{\cG}_0$,
 $\big(\cQ^{(n)}\bJot\big)\!(z):\what{\cG}_0\rightarrow
{\cG}_0$.  Thus, as in the case of (\ref{LTL1}) one may use
Eqs.~(\ref{MFEmn}) as an instrument to obtain representations
for the half\,--\,on\,--\,shell kernels $(\bJo M)(z),$
$\big(M\bJot\big)(z)$ as well as for the on\,--\,shell kernels
$\what{\cT}_{00}(z),$ \vs{0.05cm}$\what{\cT}_{01}(z)$  and
 $\what{\cT}_{10}$ in terms of the Faddeev components
$M\abd(P,P',z)$ with real $P,P'\in\Rs$.

The exposition will be given for the case of
the matrices $\bigl(M\bJot \bigr)(z)$.

First, we find easily that
continuation on $\wtilde{\Pi}_b^{(0)\pm}$ of the form
$$
\big(F,\bigl(\cQ^{(0)}\bJot \bigr)(z)f\big) \EQ
\Sum\ad\big(F\ad,\bigl(\bt\ad\bJot \bigr)(z)f\ad\big)
$$
for the matrix $\bigl(\bt\ad\bJot \bigr)(z)$
is described by the equalities
\begin{eqnarray} \hs{1.0cm}
\label{FtJ0f}
& & \!\! \bigl(F\ad,\big(\bt\ad\rJot\big)(z)f\ad \bigr) \nonumber \\ [0.1cm]
& \equiv & \!\!\int_\Rt dk\ad \int_{S^2}d\hk'\ad      \int_{S^2}d\hp'\ad
\int_0^{\pi/2}   d\omega'\ad\sin^2\omega'\ad
\cos^2\omega'\ad     \\ [0.1cm]
& & \times\;\; t\ad\big(k\ad,\sqrt{z}\cos\omega'\ad\hk'\ad, z\cos^2\omega'\ad
\big)  F\ad\big(k\ad, \pm\sqrt{z}\sin\omega'\ad\hp'\ad\big)\cdot
f\ad\big(\omega'\ad,\hk'\ad,\hp'\ad\big)\,,           \nonumber
\end{eqnarray}
where $\omega'\ad,\hk'\ad,\hp'\ad$ are hyperspherical
coordinates~\cite{MF} of the point
$\hP'\in S^5$, $\omega'\ad\in\left[0,\pi/2\right]$,
$\hk'\ad,\,\hp'\ad \in S^2$.
Note  that
$$\hP' \EQ \lbrace \cos\omega'\ad\hk'\ad,\sin\omega'\ad\hp'\ad\rbrace
\quad \mbox{while} \quad
d\hP' \EQ \sin^2\omega'\ad\cos^2\omega'\ad \,d\omega'\ad \,d\hk'\ad \,
d\hp'\ad
$$
is a measure on $S^5$.  A holomorphy domain of the function
$\bigl(F\ad,\big(\bt\ad\rJot\big)(z)f \bigr)$ can be found from the
conditions that the poles of the T\,--\,matrix $t\ad(\,\cdot\,
,\,\cdot\, ,z\cos^2\omega'\ad)$ corresponding to the discrete
spectrum of the Hamiltonian $h\ad$ do not manifest themselves.
In other words, one has to require the equalities
\be
\label{poles23}
  z\cos^2\omega'\ad \EQ \lambda\ajd\,,\quad \anum\,,
\ee
to take place for no $\omega'\ad\in[0,\pi/2]$.  Evidently, the
last requirement is equivalent to making a cut along the ray
$(-\infty,\lambda_{\rm max}].$ The poles (\ref{poles23})
generate in (\ref{FtJ0f}) some integrals of the Cauchy type
analogous to (\ref{Phi}). Thus, a continuation of the
function (\ref{FtJ0f}) through the cut $(-\infty,\lambda_{\rm
max}]$ may be described using the representations (\ref{Phil})
while each point $\lambda\ajd,$
$j=1,2,\,\ldots,n\ad$ turns into a branch point of the Riemann surface
of the function (\ref{FtJ0f}).

In the case of the potentials (\ref{vpotb}),
additionally to the cut
$(-\infty,\lambda_{\rm max}]$, there appear additional
restrictions on the domain of holomorphy of this function
following from a requirement for the second argument
of  the T\,--\,matrix $t\ad$  to belong to the strip  $W_b$,
$|\Img\sqrt{z}\cos\omega'\ad\hk'\ad|<b$
for \vs{0.05cm}all  $\omega'\ad\in[0,\pi/2],$\, $\hk'\ad\in S^2$.
This means that $z$ has to be such that
$|\Img\sqrt{z}|<b$, i.\,e., $z\in\cP_b$ [see Eq.~(\ref{Pib})].

Note that for $z\neq\lambda\pm i\hs{0.03cm}0$, $\lambda \leq
\lambda_{\rm max}$ one can substitute any element
$f\in\what{\cH}_0$ in the form  (\ref{FtJ0f}).

We shall consider formulas (\ref{FtJ0f}) with $\an=1,2,3$ as
a definition of analytic continuation of the matrix
$\big(\bt\bJot\big)(z)$ on a domain of complex $z$.  Up to
now, we have established that immediate continuation of
$\big(\bt\bJot\big)(z)$ described by the formulas
(\ref{FtJ0f}) is possible on the domain
$\cP_b\setminus(-\infty,\lambda_{\rm max}]$.

Further, using Lemma~\ref{LEqQuadrPrime} we find that the
analytic continuation of the form
\linebreak
$\bigl(F,\bigl(\cQ^{(1)}\bJot \bigr)(E\pm i0)f\bigr)$
on a domain of complex $z\in{\C}^\pm$ is given by
\be
\label{Q1Distr12}
\big(F,\bigl(\cQ^{(1)}\bJot \bigr)(z)f\big) \EQ
\Sum_{\an,\bn,\, \an\neq\bn }
Q_{1,\an\bn}^\pm(z)\,+\,Q_{2,\an\bn}^\pm(z)
\ee
where
\be \hs{0.5cm}
\label{Q1cont}
\begin{array}{lcl}
\vs{0.2cm}
Q_{1,\an\bn}^\pm(z)\!\!\!\! & = & \!\!\!\!\pm\,\Frac{\sqrt{z}}{4}
\Frac{1}{|s\abd|}
\Int_{\Rt}\! d{k}\ad \Int_{S^2}\! d\hat{p}\ad
  \Int_{S^2} \!dk'\bd \Int_{S^2} \! d\hat{p}'\bd
  \Int_0^{1} \! d\nu\, \sqrt{\nu} \\
\vs{0.2cm} & & \!\! \times \;\;
\Int_0^{1} d\nu'\, \sqrt{\nu'}\,\sqrt{1-\nu'}    \,
   \Frac{ F\ad\big({k}\ad,\sqrt{z}\,\sqrt{\nu}\,\hat{p}\ad\big) \cdot
 f\bd\big(\sqrt{1-\nu'}\,\hat{k}'\bd,\sqrt{\nu'}\,\hat{p}'\bd) }
   {\nu+\nu'-2c\abd\,\sqrt{\nu}\,\sqrt{\nu'}\,\big(\hat{p}\ad,\hat{p}'\bd\big)
   -s\abd^2 \mp  i\hs{0.03cm}0}               \\
& & \!\!  \times \;\; t\ad\big({k}\ad,
   \tk\bu\ad\big(\sqrt{z}\,\sqrt{\nu}\hat{p}\ad,\sqrt{z}\,\sqrt{\nu'}\,
   \hat{p}'\bd\big),
   z(1-\nu)\big)       \\ [0.15cm]
 & & \!\!  \times \;\;
t\bd\big(\tk\au\bd(\sqrt{z}\,\sqrt{\nu'}\,\hat{p}'\bd,\sqrt{z}\,\sqrt{\nu}\,
\hat{p}\ad\big),
\sqrt{z}\,\sqrt{1-\nu'}\,\hat{k}'\bd, z(1-\nu')\big)
\end{array}
\ee
and
\be \hs{0.5cm}
\label{Q2cont}
\begin{array}{lcl}
 & & \!\! Q_{2,\an\bn}^\pm(z) \\ [0.2cm]
 \vs{0.2cm}& = & \pm\Frac{1}{4}
\cdot\Frac{1}{|s\abd|}
\Int_{\Rt}\! d{k}\ad \Int_{S^2}\! d\hat{p}\ad
  \Int_{S^2}\! dk'\bd \Int_{S^2} \!d\hat{p}'\bd
  \Int_{\Gamma_z^\pm} \! d\rho\, \sqrt{\rho} \\
\vs{0.2cm}& & \times\;\;  \Int_0^{1} \! d\nu'\, \sqrt{\nu'}\,\sqrt{1-\nu'}
                 %
  s\, \Frac{ F\ad({k}\ad,\pm\sqrt{\rho}\,\hat{p}\ad) \cdot
 f\bd\big(\sqrt{1-\nu'}\,\hat{k}'\bd,\sqrt{\nu'}\,\hat{p}'\bd\big) }
{\rho+z\nu'-2c\abd\sqrt{z}\,\sqrt{\rho}\,\sqrt{\nu'}\,\big(\hat{p}\ad,
\hat{p}'\bd\big)   -s\abd^2 z}    \\
& &   \times\;\; t\ad\big({k}\ad,
   \tk\bu\ad\big(\pm\sqrt{\rho}\,\hat{p}\ad,\sqrt{z}\,\sqrt{\nu'}\hat{p}'\,
   \bd\big),    z-\rho\big)         \\ [0.15cm]
& &   \times \;\;
t\bd\big(\tk\au\bd\big(\sqrt{z}\,\sqrt{\nu'}\,\hat{p}'\bd,\pm\sqrt{\rho}\,
\hat{p}\ad\big),
\sqrt{z}\,\sqrt{1-\nu'}\,\hat{k}'\bd, z(1-\nu')\big)\,.
\end{array}
\ee
Here, by $\Gamma_z^+$ ($\Gamma_z^-$) we understand the path of
integration beginning at the point $z$ and going clockwise
(counterclockwise) along the circumference $C_{|z|}$ having
radius $|z|$ and centered in the origin.  After the path crosses
the real axis, it goes further along this axis so that the rest
of $\Gamma_z^+$ ($\Gamma_z^-$) consists of the points
$\rho=\lambda + i\hs{0.03cm}0$\, ($\rho=\lambda - i\hs{0.03cm}0$),\,
$\lambda\in(|z|,+\infty).$ It should be noted
that as $\Gamma_z^\pm$ one can also
choose arbitrary equivalent paths having no intersections
with the line segment $[0,\,z]$ except at the point  $z$ and
with the circle of radius $c\abd^2|z|$ centered at the origin which
are spoken about in Lemma~\ref{LEqQuadrPrime}.

One can find easily that if $f\big(\hP\hs{0.02cm}\big)$ is a H\"{o}lder function
with the smoothness index $\mu > 0$ then the functions
$Q_{1,\an\bn}^+(z)$ and $Q_{2,\an\bn}^+(z)$ $\big[Q_{1,\an\bn}^-(z)$
and $Q_{2,\an\bn}^-(z)\big]$ are holomorphic in a domain of the
upper half\,--\,plane ${\C}^+$ [of the lower half\,--\,plane ${\C}^-$], a
boundary of which is determined  by the
numerators of the integrands in (\ref{Q1cont}) and
(\ref{Q2cont}).  In the case of the functions
$Q_{1,\an\bn}^\pm(z)$ this fact does not require  special
explanation.

Concerning the functions $Q_{2,\an\bn}^\pm(z)$ we find that,
according to Lemma~\ref{LEqQuadrPrime}, the denominators of
the respective expressions under the integration signs may become
zero for $\rho=z$ only.  The latter is possible for the points
corresponding to
$\nu'=c\abd^2$ and $\eta=(\hat{p}\ad,\hat{p}'\bd)=-1$.
Consequently, it suffices
to check holomorphy of the integral in a small
vicinity of this point; namely, the integral
$$
B(z)\EQ \Int_{c\abd^2- \varepsilon}^{c\abd^2+ \varepsilon} d\nu'
\Int_{-1}^{-1+\delta} d\eta
\Int_{z}^{z+z\zeta} d\rho \,
\Frac{\tilde{f}(\nu',\eta,\rho,z)}
{\rho+z\nu'
- 2 c\abd\,\sqrt{z} \,\sqrt{\rho}\,\sqrt{\nu'}\,\eta-s\abd^2 z}
$$
where $\tilde{f}$ stands for the numerator of
$Q_{2,\an\bn}^\pm(z)$ and  $\zeta$ is a certain complex number,
$ 0 < |\zeta| < 1,$ ${\rm arg}\hs{0.02cm}\zeta\sim -\Frac{\pi}{2}.$ At the
same time $\varepsilon,\delta$ are sufficiently small positive
numbers, $0<\varepsilon<\Min{}\hs{0.03cm}\big\{c\abd^2,s\abd^2\big\}$ and
$0<\delta<1$.

Let us rewrite $\tilde{f}$ in the form
$\tilde{f}(\nu',\eta,\rho,z)=$
$\tilde{f}\big(c^2,-1,z,z\big)+\delta\tilde{f}(\nu',\eta,\rho,z)$. Since
$\delta\tilde{f}(c^2,-1,z,z)=0,$
a contribution of the summand $\delta\tilde{f}$ to $B$
is a holomorphic function.
In the term generated by the summand $\tilde{f}(c^2,-1,z,z)$,
the latter may be transferred through the integration sign.
Making the substitution $\rho=z\nu'$ in the remaining integral
one gets
$$
\wtilde{B}(z) \EQ \Frac{1}{z}\,
\Int_{c\abd^2- \varepsilon}^{c\abd^2+ \varepsilon} d\nu
\Int_{-1}^{-1+\delta} d\eta
\Int_{1}^{1+\zeta} d\nu'\,
\Frac{1}{\nu +\nu'-
2 c\abd\,\sqrt{\nu}\,\sqrt{\nu'}\,\eta-s\abd^2 }
$$
where the integral converges being independent of $z$ at all.

Thus, we show the form
$\Big(F\ad,\big(\cQ^{(1)}\bJot\big)\abd(E\pm i\hs{0.03cm}0)f\bd\Big)$
admits, in correspondence with the sign ``$\pm$'',
analytic continuation in $z$ both in ${\C}^+$
and ${\C}^-$ on the domains of holomorphy of the above nominators.

Boundaries of these domains are found from those requirements
that the poles of the T\,--\,matrices $t\ad(\,\cdot\, ,\,\cdot\,
,z(1-\nu))$ and  $t\bd(\,\cdot\, ,\,\cdot\, ,z(1-\nu'))$ which
are present in the integral (\ref{Q1cont}) do not appear
in the above domains. Also, we require the same for the
poles of the T\,--\,matrices $t\ad(\,\cdot\, ,\,\cdot\, ,z-\rho)$
which appear in the integral (\ref{Q2cont}).  If
$z\not\in(-\infty,\lambda_{\rm max}]$ then the
conditions
$
z(1-\nu)=\lambda\ajd,\; j=1,2,\,\ldots,n\ad,\;
z(1-\nu')=$ $ =\lambda_{\bn,k},\; k=1,2,\,\ldots,n\bd
$
of appearance of the poles of the T\,--\,matrices $t\ad(\,\cdot\,
,\,\cdot\, ,z(1-\nu))$ and $t\bd(\,\cdot\, ,\,\cdot\,
,z(1-\nu'))$ are valid for no $\nu,\nu'\in[0,\,1]$.  The
appearance conditions $z-\rho=\lambda\ajd,$
$j=1,2,\,\ldots,n\ad,$ of the  poles of $t\ad(\,\cdot\, ,\,\cdot\,
,z-\rho)$ may be realized only if  the paths $\Gamma_z^\pm$
include  more than one fourth of the circumference
$C_{|z|}$.  However, their contribution to $Q_{2,\abd}^\pm(z)$
arising when the points $\rho=z-\lambda\ajd$ cross the contours
$\Gamma_z^\pm$ may always be taken into account using the
residue theorem.  We shall not present here respective formulas.
Note only that  taking residues at the points
$\rho=z-\lambda\ajd$ transforms the minor three\,--\,body pole
singularities of the integrand of $Q_{2,\an\bn}^\pm(z)$ into
those of type
$
\big(z-\lambda\ajd+z\nu'-2c\abd\sqrt{z}\,\sqrt{z-\lambda\ajd}
\,\sqrt{\nu'}\,\eta-s\abd^2 z\big)^{-1}.
$
The location of such singularities is described by
Lemma~\ref{LEq23}. As a result one finds that there are sets
$\wtilde{\Pi}_b^{(0)\pm}$ which are  holomorphy domains of the
form (\ref{Q1Distr12}). It should be noted only that the extra
conditions (\ref{CondbRes0}) arise as a result of the requirements
\begin{eqnarray*}
& \big|\Img \tk\ad\bu\big(\sqrt{z}\,\sqrt{\nu}\,\hp\ad,\sqrt{z}\,\nu'\,\hp'
\bd\big)\big| \LT b\,, &
\quad
\big|\Img \tk\bd\au\big(\sqrt{z}\,\sqrt{\nu'}\,\hp'\bd,\sqrt{z}\,\nu\,
\hp\ad\big) \big| \LT  b\,, \\ [0.2cm]
& \big|\Img \tk\ad\bu\big(\pm\sqrt{\rho}\,\hp\ad,\sqrt{z}\,\nu'\,\hp'\bd\big)
\big| \; \LT b\,, &
\quad
\big|\Img \tk\bd\au\big(\sqrt{z}\,\sqrt{\nu'}\,\hp'\bd,\pm\sqrt{\rho}\,\hp\ad
\big)\big| \LT b\,,
\end{eqnarray*}
where $\nu,\nu'\in[0,\,1]$, $\hp\ad,\hp'\bd\in S^2$,
$\rho\in\Gamma_z^\pm,$
and the condition
$$
\big|\Img \big(\sqrt{z}\,\sqrt{1-\nu'}\,\hk'\bd\big)\big| \LT b\,, \quad
\hk'\bd\in S^2\,,
$$
which have to be satisfied for the arguments of the T\,--\,matrices
$t\ad$ and $t\bd$ appearing in the expressions of
$Q_{1,\an\bn}^\pm$ and $Q_{2,\an\bn}^\pm$ under the integration
signs (since these arguments have to belong to the analyticity
strip $W_b$).

The matrix $\big( \cQ^{(2)}\bJot \big)(z)$, $z=E\pm i\hs{0.03cm}0,$\,
$E >0$, corresponding to the second iteration of the absolute
term of (\ref{MFE}), has the kernels
$$
\big( \cQ^{(2)}\bJot \big)\big(P,\hP',z\big) \;=
\Sum_{\gamma\neq\an,\, \gamma\neq\bn}
Q_{\an\gamma\bn,0}^{(2)}\big(P,\hP',z\big)
$$
where
\begin{eqnarray}
\nonumber
& & \!\! Q_{\an\gamma\bn,0}^{(2)}\big(P,\hP',z\big)\\ [0.1cm]
 & = & \!\!
\Frac{1}{|s_{\an\gamma}||s_{\gamma\bn}|}\int_\Rt \,dq\,
t\ad\bigl(k\ad,\tk\ad\gu(p\ad,q),z-p\ad^2\bigr)\nonumber \\  [-0.2cm]
 & &  \\ [-0.2cm]
\label{Q2agb0}
 & & \!\! \times\;\; \Frac{t\gd\bigl(\tk\gd\au(q,p\ad),
\tk\gd\bu(q,p\bd'),z-q^2\bigr)}
{\left(p\ad^2+q^2-2c_{\an\gamma}(p\ad,q)-s_{\an\gamma}^2 z\right)
\big({p'\bd}^2+q^2-2c_{\bn\gamma}(p'\bd,q)-s_{\bn\gamma}^2 z\big)} \nonumber
            \\ [0.1cm]
 & & \times\;\;
t\bd\big(\tk\bd\gu(p'\bd,q),k'\bd,z-{p'\bd}^2\big) \nonumber
\end{eqnarray}
with $\gn\neq\an$, $\gn\neq\bn$ and $P'=\pm\sqrt{E}\hP'.$
The existence of analytic continuation of these kernels onto domains
$\wtilde{\Pi}_b^{(0)\pm}$ may be proved, in the sense of
distributions over $\Osix$, following the same scheme as for the
kernels of the matrix $\big( \cQ^{(1)}\bJot \big)(z)$.
However, it follows from the results of~\cite{Faddeev63},
\cite{MF} that these kernels have ``better'' properties than
those of $\big( \cQ^{(1)}\bJot \big)(z)$.  We have in mind
now the fact that the components
$\cF\abd(P,P',z)$,
$\cI_{\an,j;\,\bn}(p\ad,P',z),$
$\cJ_{\an;\,\bn,k}\big(P,p'\bd,z\big)$
and
$\cK_{\an,j;\,\bn,k}\big(p\ad,p'\bd,z\big)$,
$P,P'\in\Rs$, $p\ad,p'\bd\in\Rt$,
of the iteration $\cQ^{(2)}(z)$ of the absolute term of
Eq.~(\ref{MFE}) turn out to be functions having weaker
singularities than the components of  $\cQ^{(1)}(z)$.  In
particular, the main singularities of the kernel
$\cF\abd(P,P',z)$ are described by
\be
\label{It2Sing}
\Frac{\pi^2 i}{D}\left\{
{\sf f}({\rm a})\ln\big(\sqrt{\xi}+\sqrt{\zeta}+D\big)-
{\sf f}({\rm b})\ln \big(\sqrt{\xi}+\sqrt{\zeta}-D\big)
\right\}
\ee
with ${\rm a}=c_{\an\gn} p\ad$, ${\rm b}=c_{\bn\gn} p'\bd,$
$\xi=s_{\an\gn}^2 (z-p\ad^2),$ $\zeta =s_{\bn\gn}^2
(z-{p'\bd}^2)$ and $D=\sqrt{({\rm a}-{\rm b})^2}$ for
$\gn\neq\an$, $\gn\neq\bn$.  The notation ${\sf f}(q)$ is used
here for the numerator of the expression under the integration sign
in (\ref{Q2agb0}) after the replacements
$t\ad\rightarrow\tilde{t}\ad$ and $t\bd\rightarrow\tilde{t}\bd$.
The expressions for the main singularities of the kernels
$\cI_{\an,j;\,\bn}(p\ad,P',z),$
$\cJ_{\an;\,\bn,k}(P,p'\bd,z)$, and
$\cK_{\an,j;\,\bn,k}(p\ad,p'\bd,z)$ may also be represented in the
form (\ref{It2Sing}) but one has to take ${\sf f}(q)$ as the
numerator of the expression in (\ref{Q2agb0})
replacing $t\ad\rightarrow\tilde{\phi}\ad\big(\tk\ad\gu(p\ad,q)\big)$
and/or $t\bd\rightarrow{\phi}\bd\big(\tk\bd\gu(p'\bd,q)\big)$.
It should be emphasized that the
singularities (\ref{It2Sing})
only appear for $p\ad^2\leq |z|$ simultaneously with ${p'\bd}^2\leq
|z|$.

Thereby, when continuing the form
$\big(F,\bigl(\cQ^{(2)}\bJot \bigr)(z)f\big)$
we get for it the representations which differ from
(\ref{Q1Distr12}) -- (\ref{Q2cont}) mainly in the
replacement of the distributions
$
   \big\{ z\big(\nu+\nu'-
   2c\abd\,\sqrt{\nu}\,\sqrt{\nu'}\,\big(\hp\ad,\hp'\bd\big) -s\abd^2 \mp
   i\hs{0.03cm}0    \big)\big\}^{-1}
$,
$0\leq \nu \leq 1,$  $0\leq \nu'\leq 1,$
with functions singular like
\be
\label{It20Sing}
\begin{array}{l}
\vs{0.2cm}\Frac{1}{z\, \big|c_{\an\gn}\nu\hp\ad-c_{\bn\gn}\nu'\hp'\bd\big|}
\\
\times \;\;\mbox{\large $\ln$}
\Frac{   \sqrt{s_{\an\gn}^2\big(1-\nu^2\big)}+\sqrt{s_{\bn\gn}^2
\big(1-{\nu'}^2\big)}+
\big|c_{\an\gn}\nu\hp\ad-c_{\bn\gn}\nu'\hp'\bd\big|   }
{   \sqrt{s_{\an\gn}^2(1-\nu^2)}+\sqrt{s_{\bn\gn}^2(1-{\nu'}^2)}-
\big|c_{\an\gn}\nu\hp\ad-c_{\bn\gn}\nu'\hp'\bd\big|   }\,.
\end{array}
\ee

The kernels $\cF\abd(P,P',z)$,
$\cI_{\an,j;\,\bn}(p\ad,P',z),$
$\cJ_{\an;\,\bn,k}\big(P,p'\bd,z\big)$,
and
$\cK_{\an,j;\,\bn,k}\big(p\ad,p'\bd,z\big)$\,
of the iteration
$\cQ^{(3)}(z)=\left(-\bt(z)\bRo(z)\Y\right)^3\bt(z)$ are still
singular. Though their singularities are weak we understand the
continuation of the kernels $\big(\cQ^{(3)}\bJot\big)(z)$ on
the domains $\wtilde{\Pi}_b^{(0)\pm}$ as before in the sense of
distributions over $\Osix$.  So, we realize this continuation
following the same scheme as for the continuation of
$\big(\cQ^{(1)}\bJot\big)(z)$ and
$\big(\cQ^{(2)}\bJot\big)(z)$.

As mentioned above, the components $\cF\abd(P,P',z),$
$\cI_{\an,j;\,\bn}(p\ad,P',z),$
$\cJ_{\an;\,\bn,k}(P,p'\bd,z)$,
and $\cK_{\an,j;\,\bn,k}(p\ad,p'\bd,z)$\,
of the subsequent iterations $\cQ\un(z)$
turn out to have no singularities.

Using the representations for the three\,--\,body singularities  of
$\cQ^{(2)}(z)$, one can show that the kernels $\cF\abd\big(P,\sqrt{z}\,\hP',
z\big),$
$\cI_{\an,j;\,\bn}\big(p\ad,\sqrt{z}\,\hP',z\big),$
$\cJ_{\an;\,\bn,k}\big(P,\sqrt{z}\,\sqrt{\nu'}\,\hp'\bd,z\big)$
and
$\cK_{\an,j;\,\bn,k}\big(p\ad,\sqrt{z}\,\sqrt{\nu'}\,\hp'\bd,z\big)$\, being
components of the matrix $\big( \cQ^{(4)}\bJot  \big)(z)$,
turn out at  $P\in\Rs$, $p\ad\in\Rt$, $\hat{P}'\in
S^5$, $\hat{p}'\bd\in S^2$  and  $0\leq \nu'\leq 1$
to be holomorphic functions of $z\in\wtilde{\Pi}_b^{(0)\pm}$.
The aggregate of these kernels may be embedded for fixed
$\what{P}'$, $\hat{p}'\bd$ and $ \nu'$ into the Banach
space $\cB_{\theta\mu}$, $\theta <\theta_0,$  $\mu<1/8.$
And this aggregate becomes a function continuous in $z$ with
respect to the norm in $\cB_{\theta\mu}$ right up to the edges of
the cut $z=E\pm i\hs{0.03cm}0$, $E > 0$.  All the subsequent iterations
$\big( \cQ^{(n)}\bJot \big)(z)$, $n\geq 5$, possess this
property, too.

Assertions similar to those obtained concerning the continuation
of the matrices $\big(\cQ^{(n)}\bJot\big)(z)$ also hold as
well for the matrices $\big(\bJo\cQ^{(n)}\big)(z)$.

Now, we use the equations (\ref{MFEmn}) for $m,n\geq 4$
and thereby complete the proof.%
\end{proof}

In the same way as Theorems~\ref{ThLTL} -- \ref{ThMJ0tJ0M},
the two following assertions may be proved.
\begin{theorem}\label{ThJ0TJ0t}
The matrix $\bigl(\bJo M\bJot \bigr)(z)$ {\rm(}the operator
$\bigl(\rJo T\rJot \bigr)(z)${\rm)} admits analytic continuation in
$z$ from the edges of the cut $z=E\pm i\hs{0.03cm}0$, $E> 0$, to
the domains $\wtilde{\Pi}_b^{(0)\pm}\in{\C}^\pm$ as a bounded operator in
$\what{\cG}_0$ $($in $\what{\cH}_0)$.  In addition, $\bigl(\bJo
M\bJot \bigr)(z)$, $z\in\wtilde{\Pi}_b^{(0)\pm}$ admits the
representation {\rm [cf.~(\ref{MQW})]}\,\,
$$
\bigl(\bJo M\bJot \bigr)(z) \EQ \Sum_{n=0}^3\, \bigl(\bJo \cQ\un\bJot \bigr)(z)
+\bigl(\bJo \cW \bJot \bigr)(z).
$$
The operators $\bigl(\bJo \cQ^{(0)}\bJot \bigr)(z)$ and
$\bigl(\bJo \cQ^{(1)}\bJot \bigr)(z)$ are bounded matrix
operators in $\what{\cG}_0$ with singular kernels.  Having weakly
singular kernels, the matrices $\bigl(\bJo \cQ\un\bJot
\bigr)(z)$, $n=2,3$, are compact operators in
$\what{\cG}_0$.  The kernels of the matrix $\bigl(\bJo
\cW \bJot \bigr)(z)$ are H\"older functions of their arguments
with the smoothness index $\mu\in\bigl(0,1/8 \bigr)$.
\end{theorem}
\begin{theorem}
\label{ThJ0MYPsiJ1t}
The operators
\begin{eqnarray*}
& \bigl(\bJo M\Y\Psi\rJt_1 \bigr)(z)\;: \;
\what{\cH}_1 \; \longrightarrow \; \what{\cG}_0\,, &
\bigl(\rJ_1\Psis\Y M\bJot \bigr)(z) \; : \;
\what{\cG}_0 \; \longrightarrow \; \what{\cH}_1\,, \\ [0.1cm]
& \what{\cT}_{01}(z)\;:\;\what{\cH}_1 \; \longrightarrow \; \what{\cH}_0\,, &
\what{\cT}_{10}(z) \;:\;\what{\cH}_0 \;\longrightarrow \; \what{\cH}_1
\end{eqnarray*}
admit analytic continuation from the edges of the cut
$z=E\pm i\hs{0.03cm}0,$ $E>0$,
onto the domains $\Pi_b^{(0)\pm}\subset{\C}^\pm$
including the points
$z\in\wtilde{\Pi}_b^{(0)\pm}\!\Bigcap\nolimits_{\bn,j}\Pi_b^{(\bn,j)}$
satisfying the additional conditions
$$
\Real z \GT  \Frac{|s_{\bn\gn}|^2}{(1+|c_{\bn\gn}|)^2}\,\lambda\bjd
\, + \, \Frac{(1+|c_{\bn\gn}|)^2}{4\,|s_{\bn\gn}|^2 |\lambda\bjd|}\,
(\Img z)^2
$$
for any $\bn,\gn=1,2,3,$ $\bn\neq\gn,$ and  $j=1,2,\,\ldots,n\bd.$
For all $z\in\Pi_b^{(0)\pm}$
including the boundary points $z=E\pm i\hs{0.03cm}0$, $E>0$,
these operators are compact.
\end{theorem}

\vs{0.3cm}As a comment to Theorem~\ref{ThJ0TJ0t} we present explicit formulas
for the kernels of the operators
$\big(\bJo \cQ^{(0)}\bJot \big)(z)$
and $\big(\bJo \cQ^{(1)}\bJot \big)(z)$.

The kernels of the first operator have the form
            $$
\big(\bJo \cQ^{(0)}\bJot \big)\abd\big(\hP,\hP',z\big) \EQ
\delta\abd\big(\rJo\bt\ad\rJot \big)\big(\hP,\hP',z\big)\,,
            \quad \an,\,\bn \EQ 1,\,2, \,3\,,
	    $$
where
\be
\label{JTJS}
\begin{array}{lcl}
\big(\rJo\bt\ad\rJot\big)\big(\hP,\hP',z\big) \!\! & = & \!\!
t\ad\big(\sqrt{z}\cos\omega\ad\hk\ad,\sqrt{z}\cos\omega{'}\ad\hk{'}\ad,
z\cos^2\omega\ad\big)\, \\ [0.15cm]
& & \!\!              \times\;\;
\delta\big(\sqrt{z}\sin\omega\ad\hp\ad - \sqrt{z}\sin\omega{'}\ad\hp'\ad
\big)\,.
\end{array}
\ee
Here,  $\omega\ad,\,\hk\ad,\,\hp\ad$ and $\omega'\ad,\,\hk'\ad,\,\hp'\ad$
are coordinates of the points
$\hP=\lbrace k\ad,p\ad\rbrace$ and \linebreak
$\hP'=\lbrace k'\ad,p'\ad\rbrace$
on the hypersphere $S^5$.
We mean here that
\be
\label{deltps}
\delta\big(\sqrt{z}\sin\omega\hs{0.02cm}\hp - \sqrt{z}\sin\omega{'}\hp'\big)
\EQ
{\rm Sign}\hs{0.02cm}\Img z\cdot
\Frac{\delta(\hp,\hp')\hs{0.02cm}\delta(\omega-\omega')}
{\big(\sqrt{z}\,\big)^3 \sin^2\omega\cos\omega}
\ee
where  $\delta(\hp,\hp')$ is the kernel of the identity operator in
$L_2\big(S^2\big)$.  The denominator \linebreak
$\big(\sqrt{z}\,\big)^3\sin^2\omega\cos\omega$ of the right\,--\,hand side of
Eq.~(\ref{deltps}) represents the analytic continuation of the
Jacobian corresponding to respective substitution of variables.

Therefore the operator $\bigl(\rJo\bt\ad\rJot\bigr)(z)$ acts at
$\Img z\neq 0$ on $f\in\what{\cH}_0$ as
\begin{eqnarray}
\nonumber
 \bigl(\bigl( \rJo\bt\ad\rJot \bigr)(z)f \bigr)\big(\hP\big)\!\! & =
& \!\! \Frac{{\rm Sign}\hs{0.02cm}\Img z}{\bigl(\sqrt{z} \bigr)^3}\cdot
\Int_{S^2} d\hk'\ad\,\,   \\ [0.1cm]
\label{J0t2J0t}
&& \!\! \times\;\;
t\ad\big(\sqrt{z}\cos\omega\ad\hk\ad,\,\sqrt{z}\cos\omega\ad\hk'\ad,\,
z\cos^2\omega\ad\big) \\ [0.1cm]
& & \! \! \times \;\; f\big(\cos\omega\ad\hk'\ad,\sin\omega\ad\hp\ad\big)\,.
\nonumber
\end{eqnarray}

The operators $\bigl(\bJo\cQ^{(1)}\bJot \bigr)(z)$,
$z\in\wtilde{\Pi}_b^{(0)\pm}$, have the kernels
$$
%
\bigl(\bJo\cQ^{(1)}\bJot \bigr)\abd\big(\hP,\hP',z\big)\!=\!
\Frac{1}{z}\cdot\Frac{1-\delta\abd}{|s\abd|}\cdot
 \Frac{ t\ad\bigl(k\ad ,k\ad\bu ,z(1-\nu) \bigr)
\,\,\, t\bd\bigl(k\bd\au ,k'\bd, z(1-\nu') \bigr)  }
{\nu+\nu'-2c\abd\,\sqrt{\nu}\,\sqrt{\nu'}\,(\hp\ad,\hp'\bd)-s\abd^2\mp
i\hs{0.03cm}0}\,,
%
$$
where $k\ad=\sqrt{z}\,\sqrt{1-\nu}\,\hk\ad,$
$k'\bd=\sqrt{z}\,\sqrt{1-\nu'}\,\hk'\bd,$
$
k\ad\bu=\tk\ad\bu\bigl(\sqrt{z}\,\sqrt{\nu}\,\hp\ad,\,
\sqrt{z}\,\sqrt{\nu'}\,\hp'\bd \bigr)
$ \linebreak
and
$
k\bd\au=\tk\bd\au\bigl(\sqrt{z}\,\sqrt{\nu'}\,\hp'\bd,\,
\sqrt{z}\,\sqrt{\nu}\,\hp\ad \bigr).
$
At the same time $\nu=\sin^2\omega\ad$ and $\nu'=\sin^2\omega'\bd$.

The main singularities of the kernels
$\bigl(\bJo\cQ^{(2)}\bJot\bigr)\abd\big(\hP,\hP',z\big)$ in $\hP$,
$\hP'$ are described by Eqs.~(\ref{It20Sing}).  The singularities of
the kernels $\bigl(\bJo\cQ^{(3)}\bJot\bigr)\abd\big(\hP,\hP',z\big)$ are
weaker.

Later, we shall use the notation
\be
\label{Pil0hol}
\Pi_{l^\pm}^{({\rm hol})} \; \equiv \; \Pi_b^{(0)\pm}\,\cap \,
\Pi_{l^{(1)}}^{({\rm hol})},
\ee
where $l^\pm=\big(l_0^\pm,l_{1,1},\,\ldots,l_{1,n_1},l_{2,1},\,\ldots,
l_{2,n_2},l_{3,1},\,\ldots,l_{3,n_3}\big)$
with  $l_0^\pm=\pm 1$,
 $l\ajd=1,$ $\alpha = 1,2,3, \;j = 1,2,\,\ldots ,n_{\alpha}$,
and  $l^{(1)}=$ $\big(0,\,\,l_{1,1},\,\ldots,l_{1,n_1},l_{2,1},\,\ldots,
l_{2,n_2},l_{3,1},\,\ldots,l_{3,n_3}\big)$ with the same $l\ajd$ as
in $l^\pm$.
Remember that the sets
$\Pi_{l^{(1)}}^{({\rm hol})}\equiv\Pi_{l^{(1)} l^{(1)}}^{({\rm hol})}$
were defined by Eqs.~(\ref{Pil1hol}).

As follows from Theorems~\ref{ThLTL}, \ref{ThJ0TJ0t}
and~\ref{ThJ0MYPsiJ1t}, the total three\,--\,body scattering matrix
$S(z)$, $z=E\pm i\hs{0.03cm}0$, $E>0$, admits the analytic continuation
as a holomorphic operator\,--\,valued function $S(z):
\what{\cH}_0\oplus\what{\cH}_1\rightarrow\what{\cH}_0\oplus\what{\cH}_1$
on the domain  $\Pi_{l^+}^{({\rm hol})}\subset{\C}^+.$ For any
$z\in\Pi_{l^+}^{({\rm hol})}$ the operator $S(z)$ is bounded.
In equal degree the same is true for $\St(z)$.

\newsection{Description of (part of) the three-body
             Riemann surface}
\label{SRiemannSurface}
By the three\,--\,body energy Riemann surface we mean
the Riemann surface of the kernel $R(P,P',z)$
of the resolvent $R(z)$ of the Hamiltonian $H$
consideration as a function of the
parameter $z$, the energy of the three\,--\,body system.

One has to expect that this surface, like that of the free Green function
$R_0(P,P',z)$, consists of an infinite number of sheets already because
the threshold $z=0$ is a logarithmic branching point.  Actually the
Riemann surface of $R(P,P',z)$ is much more complicated than that of
$R_0(P,P',z)$, since besides  $z=0$ it has a lot of additional
branching points. For example, the two\,--\,body thresholds $z=\lambda\ajd$,
$\alpha = {1,2,3 }, \;j = 1,2,\,\ldots,n_{\alpha}$ become square root
branching points of this surface.  Also,
the resonances of pair subsystems turn into such points. Extra
branching points are generated by boundaries of the supports of the
function (\ref{Znam}) singularities which were described in
Lemmas~\ref{LEqParab}, \ref{LEqQuadr} and \ref{LEq23}.

In the present paper we restrict ourselves to consider
only of a ``small'' part of the total three\,--\,body Riemann surface
for which we succeeded to find the explicit representations
expressing analytic continuation of the Green function
$R(P,P',z)$, the kernels of the matrix $M(z)$, as well as
the scattering matrix $S(z)$ in terms of the physical sheet
[see respective formulas (\ref{Ml3fin}), (\ref{Slfin}) and
(\ref{R3l})].  Namely, in the Riemann surface of $R(P,P',z)$ we
consider two neighboring ``three\,--\,body'' unphysical sheets
immediately joint with the physical one along the
three\,--\,body branch of the continuous spectrum
$[0,\,\,+\infty)$.  In addition, we examine all the
``two\,--\,body'' unphysical sheets, i.\,e., the sheets where the
parameter $z$ may be carried if going around the two\,--\,body
thresholds $z=\lambda\ajd$, $\alpha = 1,2,3,\; j= 1,2,\,\ldots,n_{\alpha}$, is
permitted but crossing
the ray $[0,\,\,+\infty)$ is forbidden.  Evidently, the part
of the three\,--\,body surface described includes all the sheets
neighboring the physical one.  The neighboring sheets are
of most interest in applications, since only resonances situated
on these ones are accessible for immediate experimental
observation.

We give a concrete description of the part under consideration
using the auxiliary vector\,--\,function
$\sff(z)=(\sff_0(z),\sff_1(z),\sff_2(z),\sff_3(z))$,
where
$\sff_0(z)=\ln z$
while
$$\sff\ad(z) \EQ \big((z-\lambda_{\an,1})^{1/2},
(z-\lambda_{\an,2})^{1/2},\,\ldots,(z-\lambda_{\an,n\ad})^{1/2}\big)\,, \quad
\an \EQ 1,\,2,\,3\,,$$
are again vector\,--\,functions.

The Riemann surface of $\sff(z)$ consists of an infinite number of
copies of the complex plane ${\C}'$ cut along the ray
$[\lambda_{\rm min},+\infty)$. These sheets are pasted together
in a suitable way along edges of the cut segments between
neighboring points in the set of thresholds $\lambda\ajd,$
$\alpha =1,2,3, \; j=1,2,\,\ldots, n_{\alpha}$ and $\lambda_0=0.$ The sheets
$\Pi_{l_0 l_1 l_2 l_3}$
are identified by indices of branches of the functions
$\sff_0(z)=\ln z$ and $\sff\ajd(z)=\big(z-\lambda\ajd\big)^{1/2}$ in
such a manner that  $l_0$ is integer and $l\ad$, $\an=1,2,3$ are
multi\,--\,indices, $l\ad=(l_{\an,1},l_{\an,2},\,\ldots,l_{\an,n\ad})$,
$l\ajd=0,1$.  For the main branch of the function $\sff\ajd(z)$,
$\alpha =1,2,3, \; j=1,2,\,\ldots, n_{\alpha}$, we take $l\ajd=0$, and
otherwise $l\ajd=1$. In case
there exist coinciding thresholds, i.\,e.,
$\lambda_{\an,i}=\lambda_{\bn,j}$ at $\an\ne \bn$ and/or $i\ne
j$ (this means that the discrete spectra of the pair Hamiltonians
coincide at least partly for two pair subsystems or at least one
of the pair subsystems has a multiple discrete spectrum), then
for each sheet  $\Pi_{l_0 l_1 l_2 l_3}$ the indices $l_{\an,i}$ and
$l_{\bn,j}$ coincide too,  $l_{\an,i}$=$l_{\bn,j}$.  As $l_0$
we choose the number of the function $\ln z$ branch, $\ln
z=\ln|z|+i \hs{0.02cm}\varphi_0+i \hs{0.02cm}2\pi l_0$ with $\varphi_0$, the
argument
of $z$, $z=|z|\,{\rm e}^{i\varphi_0},$ $\varphi_0\in[0,2\pi)$.
The sheets $\Pi_{l_0 l_1 l_2 l_3}$ are pasted together (along
edges of the cut) in such a way that if the parameter $z$ going
from the sheet $\Pi_{l_0 l_1 l_2 l_3}$ crosses the interval
between two neighboring thresholds $\lambda_{\an,i}$ and
$\lambda_{\bn,j}$, $\lambda_{\an,i} < \lambda_{\bn,j}$ (or
$\lambda_{\rm max}$ and $\lambda_0$) then it goes over the sheet
$\Pi_{l'_0 l'_1 l'_2 l'_3}$ where the indices $l_{\gamma,k}$
corresponding to $\lambda_{\gamma,k}\leq\lambda_{\alpha,i}$
$(\lambda_{\gamma,k}\leq\lambda_{\rm max})$ changed by 1.
If $l_{\gamma,k}=0$, then $l'_{\gamma,k}=1$; if
$l_{\gamma,k}=1$, then $l'_{\gamma,k}=0$.  The indices
$l_{\gamma,k}$ for $\lambda_{\gamma,k} > \lambda_{\alpha,i}$ and
$l_0$ stay unchanged:  $l'_{\gamma,k} = l_{\gamma,k}$,
$l'_0=l_0$.  In case the parameter $z$ crosses the cut on
the right from the three\,--\,body threshold $\lambda_0$ (at
$E>\lambda_0$), then all the indices $l_{\gamma,k}$ change as it was
described above.  In addition, the index $l_0$ changes by 1, too. If
at that, $z$ crosses the cut from below, then $l'_0=l_0+1$.
Otherwise $l'_0=l_0-1$.  Further, by $l$ we denote the
multi\,--\,index $l=(l_0,l_1,l_2,l_3)$.

Thus, we have described the Riemann surface of the auxiliary
vector\,--\,function $\sff(z)$.

As mentioned above we shall consider only a part of the
three\,--\,body Riemann surface which will be denoted by $\Re$.  We
include in $\Re$ all the sheets $\Pi_l$ of the Riemann surface
of the function $\sff(z)$ with $l_0=0$. Also, we include in
$\Re$ the upper half\,--\,plane $\Img z > 0$ of the sheet $\Pi_l$
with $l_0=+1$ and the lower half\,--\,plane $\Img{z}<0$ of the
sheet $\Pi_l$ with $l_0=-1$.  For these parts we keep the
previous notations  $\Pi_l$, $l_0=\pm 1$, assuming additionally
that the cuts are made on them along the rays belonging to the
set $Z_{\rm res}=\bigcup\nolimits_{\an=1}^3 Z_{\rm res}\au$.  Here,
$
Z_{\rm res}\au=
\big\{z:\, z=z_r \rho, \, 1\leq\rho <+\infty,\,
z_r\in\sigma\au_{\rm res} \big\}
$
is a totality of the rays beginning at the resonances
$z_r\in\sigma\au_{\rm res}$ of the pair subsystem
$\an$ and going to infinity
in the directions  $\hat{z}_r ={z_r}/{|z_{r}|}$.

The sheet $\Pi_l$ for which all the components of the
multi\,--\,index  $l$ are zero, $l_0=l\ajd=0$ $=0,$ $\alpha =1,2,3,
\; j=1,2,\,\ldots, n_{\alpha}$, is called the
physical sheet. The unphysical sheets $\Pi_l$ with $l_0=0$ are
called the two\,--\,body sheets since these ones may be reached by
only going around two\,--\,body thresholds and it is not necessary to
bypass the three\,--\,body threshold $\lambda_0$.  The sheets $\Pi_l$
at $l_0=\pm 1$ are called the three\,--\,body sheets.

\newsection{Analytic continuation of the Faddeev integral equations
                           into unphysical sheets }
\label{SKernels}
A goal of the present section consists in a continuation into the
unphysical sheets of the surface $\Re$ of the absolute terms and
kernels of the Faddeev equations (\ref{MFE}) and their iterations.
Continuation is realized in the sense of generalized functions
(distributions) over $\Osix$. Results of the continuation are
represented in terms related to the physical sheet only.

By $L\au$, $L\au=L\au(l),$ we denote the diagonal matrices
formed of components $l_{\an,1},$ $l_{\an,2},\,\ldots,$ $l_{\an,n\ad}$
of the multi\,--\,index $l$ of the sheet $\Pi_l\subset\Re$:
            $$
L\au \;= \;\diag\hs{0.03cm}\{ l_{\an,1}, l_{\an,2},\,\ldots,l_{\an,n\ad}\} \,.
            $$
Meanwhile
$L_1(l)=\diag\hs{0.03cm}\big\{ L^{(1)},L^{(2)},L^{(3)} \big\}$ and
$L(l)=\diag\hs{0.03cm}\{L_0,L_1\}$ with $L_0\equiv l_0$. Analogously,
                \begin{eqnarray*}
A\au(z) \!\! & = & \!\! \diag\hs{0.03cm}\{ A_{\an,1}(z),A_{\an,2},\,\ldots,
A_{\an,n\ad}(z)\}\,, \\ [0.1cm]
A_1(z) \!\! & = & \!\! \diag\hs{0.03cm}\big\{ A^{(1)}(z),A^{(2)}(z),
A^{(3)}(z)\big\}\,.
                \end{eqnarray*}
Thus $A(z)=\diag \hs{0.03cm} \{A_0(z),A_1(z)\}$.

By $\bs_{\an,l}(z)$ we understand an operator defined in
$\what{{\cal H}}_0$ by
\be
\label{spair6}
\bs_{\an,l}(z)\;=\;\hat{I}_0\,+\,\rJo(z)\bt\ad(z)\rJot(z)A_0(z)L_0\,,
\quad z\in\Pi_0\,.
\ee
It follows from Eq.~(\ref{spair6}) that $\bs_{\an,l}=\hat{I}_0$
at $l_0=0$.  If $l_0=\pm 1$, then, according to
Eqs.~(\ref{JTJS}) -- (\ref{J0t2J0t}), the operator
$\bs_{\an,l}(z)$ is defined for $z\in\cP_b\cap{\C}^\pm$
acting on  $f\in \what{\cH}_0$ by
\be
\label{spairh}
(\bs\adl(z)f)\big(\hP\hs{0.02cm}\big) \EQ
\int_{S^2} d\hk' s\ad\big(\hk\ad,\hk'\ad,z\cos^2\omega\big)
f\big(\cos\omega\ad\hk\ad',\sin\omega\ad\hp\ad\big)
\ee
where $\omega\ad,\,\hk\ad,\,\hp\ad$ stand for the coordinates~\cite{MF}
of the point $\hP$ on the hypersphere $S^5$,
$\omega\ad\in\left[0,\pi/2\right]$, $\hk\ad,\,\hp\ad \in S^2$ and
$\hP=\lbrace\cos\omega\ad\hk\ad,\,\,\sin\omega\ad\hp\ad\rbrace$.
By $s\ad$ we denote the scattering matrix (\ref{s2}) for the
pair subsystem $\alpha$.  Here we have taken into account the
fact that $l_0\cdot{\rm Sign}\hs{0.02cm}\Img z=1$ both for $l_0=1$ and
$l_0=-1$.  Recall that at $l_0=1$ the set $\Pi_l$
represents the upper half\,--\,plane and at $l_0=-1$, the lower one
(in accordance with our choice in Sec. 5 of
the part $\Re$ of a total Riemann surface in the problem of
three particles). Therefore, one can see that the operators
$\bs_{\an,l}$ are described by the same formula (\ref{spairh})
in both three\,--\,body sheets $\Pi_l$, $l_0=\pm 1$.  As a matter
of fact, $\bs\adl$ represents the pair scattering matrix $s\ad$
rewritten in the three\,--\,body momentum space.

It follows immediately from Eq.~(\ref{spairh}) that if
$z\in\cP_b\cap{\C}^\pm\setminus {Z}\au_{\rm res}$, then
the bounded inverse operator $\bs\adl^{-1}(z)$ exists and
$$
(\bs\adl^{-1}(z)f)\big(\hP\hs{0.02cm}\big) \EQ
\int_{S^2} d\hk' s\ad^{-1}\big(\hk\ad,\hk'\ad,z\cos^2\omega\ad\big)
f\big(\cos\omega\ad\hk'\ad,\sin\omega\ad\hp\ad\big)
$$
with $s\ad^{-1}\big(\hk,\hk',\zeta\big)$ the kernel of the inverse
scattering matrix $s^{-1}\ad(\zeta)$.

The operator $\bs\adl^{-1}(z)$ becomes unbounded
at the boundary points $z$ situated on the edges of the cuts
(the ``resonance'' rays) included in $Z\au_{\rm res}$.
\begin{theorem}\label{ThKernels}
The absolute terms $\bt\ad(P,P',z)$ and kernels $(\bt\ad
R_0)(P,P',z)$ of the Faddeev equations {\rm (\ref{MFE})} admit
analytic continuation in the sense of distributions over $\Osix$
both into the two\,--\,body
and three\,--\,body sheets  $\Pi_l$ of the Riemann
surface  $\Re$.  The continuation into the sheet $\Pi_l$,
$l=(l_0,l_{1,1},\,\ldots,l_{1,n_1},l_{2,1},\,\ldots,
l_{2,n_2},l_{3,1},\,\ldots,l_{3,n_3})$,  $l_0=$ $=0$, $l\bjd=0,
 1,$ or $l_0=\pm 1$,  $l\bjd=1$
{\rm(}in the both cases $\beta = 1,2,3,\; j=1,2,\,\ldots, n_{\beta}${\rm)}, is
written as
\begin{eqnarray}
\label{tlRbig}
& & \bt^{l}\ad(z)\;  \equiv  \; \reduction{\bt\ad(z)}{\Pi_l} \EQ
\bt\ad\,-\,L_0 A_0\bt\ad\rJot\bs\adl^{-1}\rJo\bt\ad\,-\,
\Phi\ad\rJ^{(\an)t}L\au A\au \rJ\au \Phi\ad^{*}\,, \\ [0.2cm]
\label{tr0l}
& & \reduction{[\bt\ad(z)R_0(z)]}{\Pi_l} \EQ
\bt\ad^l (z)R_0^l(z)\,,
\end{eqnarray}
where
$
R_0^l(z)\equiv \reduction{R_0(z)}{\Pi_l}=
R_0(z)+L_0 A_0(z)\rJot(z)\rJo(z)
$
is the continuation~{\rm (\ref{R06})} on $\Pi_l$ of the free Green
function $R_0(z)$. If $l_0=0$ {\rm(}and hence $\Pi_l$ is a two\,--\,body
unphysical sheet{\rm)}, then the continuation in the form
{\rm(\ref{tlRbig}), (\ref{tr0l})} is possible on the whole sheet
$\Pi_l$. For $l_0=\pm 1$, {\rm(}i.\,e., in the case where $\Pi_l$ is a
three\,--\,body sheet{\rm)}
the continuation in the form {\rm (\ref{tlRbig}), (\ref{tr0l})}
is possible on the domain $\cP_b\cap\Pi_l$.  All
the kernels on the right\,--\,hand side of Eqs.~{\rm (\ref{tlRbig})}
are taken in the physical sheet.
\end{theorem}
\begin{proof}
We prove the theorem for the case of the most complicated
continuation into the three\,--\,body unphysical sheets $\Pi_l$ with
$l_0=\pm 1$.  For the sake of definiteness we consider the case
$l_0=+1$.  For $l_0=-1$ the proof is quite analogous.

Let us consider at $z\in\Pi_0$, $\Img z<0$ the bilinear form
\be
\label{ftr0f}
(f,\bt\ad R_0(z)f') \EQ
\int_{\Rt}dk\int_{\Rt}dk'\int_{\Rt}dp\,
\Frac{t\ad\big(k,k',z-p^2\big)}{k'^2+p^2-z}\,
\tilde{f}(k,k',p)
\ee
with $\tilde{f}(k,k',p)=f(k,p)f'(k',p)$,
$f,f'\in\Osix$,
$k=k\ad,$ $k'=k\ad'$, $p=p\ad$.
Making the substitutions $|k'|\rightarrow\rho=|k'|^2$,
$|p|\rightarrow \lambda=z-|p|^2$ the integral
(\ref{ftr0f}) becomes
\begin{small}
\be
\label{ftr0f1}
\mbox{\phantom{MMM}}
\Frac{1}{4}\,\Int_\Rt \!dk\Int_{S^2}\!d\hk'\Int_{S^2}\!d\hp
\Int_{z-\infty}^{z}\!d\lambda\sqrt{z-\lambda}\Int_0^\infty
d\rho\sqrt{\rho}\,\Frac{t\ad\big(k,\sqrt{\rho}\hk',\lambda\big)}
{\rho-\lambda}\,
\tilde{f}\big(k,\sqrt{\rho}\,\hk',\sqrt{z-\lambda}\,\hp\big)\,.
\ee
\end{small}
The existence of an analytic continuation of the kernel $(\bt\ad
R_0)(z)$ into the sheet $\Pi_l$, $l_0=\pm 1$ follows from the
possibility of continuously deform the path of integration
in the variable $\rho$ to an arbitrary sector of the holomorphy
domain $\cP_b\cap\sigma\au_{\rm res}$ of the integrand in
the variable $\lambda$ in the way demonstrated in
Fig.~\ref{figContour-mu}.


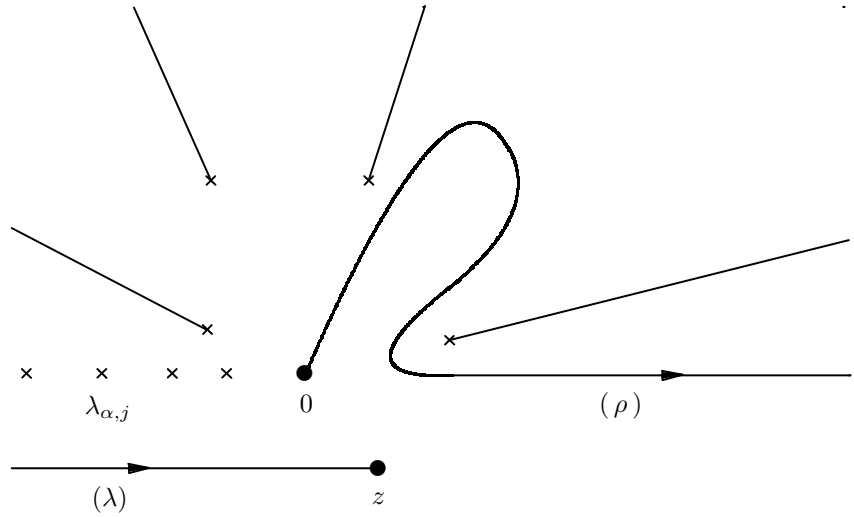
\begin{figure}
\centering
\unitlength=0.70mm
\special{em:linewidth .75pt}
\linethickness{.75pt}
\begin{picture}(159.33,95.54)
\emline{0.00}{7.67}{1}{69.67}{7.67}{2}
\emline{2.01}{24.85}{3}{3.66}{26.50}{4}
\emline{1.96}{26.50}{5}{3.72}{24.83}{6}
\emline{29.68}{24.85}{7}{31.33}{26.50}{8}
\emline{29.63}{26.50}{9}{31.39}{24.83}{10}
\emline{40.01}{24.85}{11}{41.66}{26.50}{12}
\emline{39.96}{26.50}{13}{41.72}{24.83}{14}
\emline{36.34}{33.18}{15}{37.99}{34.83}{16}
\emline{36.30}{34.83}{17}{38.06}{33.16}{18}
\emline{16.31}{24.87}{19}{17.96}{26.51}{20}
\emline{16.26}{26.51}{21}{18.02}{24.84}{22}
\bezier{372}(55.67,25.67)(83.33,89.33)(94.33,68.00)
\bezier{132}(94.33,68.33)(100.33,57.00)(85.00,43.67)
\bezier{372}(55.67,25.67)(83.33,89.33)(94.33,68.00)
\bezier{372}(55.67,25.67)(83.33,89.33)(94.33,68.00)
\bezier{132}(94.33,68.33)(100.33,57.00)(85.00,43.67)
\emline{83.67}{25.33}{23}{159.33}{25.33}{24}
\emline{159.33}{25.33}{25}{159.00}{25.33}{26}
\emline{82.34}{31.18}{27}{83.99}{32.83}{28}
\emline{82.30}{32.83}{29}{84.06}{31.16}{30}
\emline{67.01}{61.52}{31}{68.66}{63.17}{32}
\emline{66.96}{63.17}{33}{68.72}{61.50}{34}
\emline{37.01}{61.52}{35}{38.66}{63.17}{36}
\emline{36.96}{63.16}{37}{38.72}{61.50}{38}
\emline{0.00}{53.33}{39}{37.00}{34.00}{40}
\emline{23.26}{95.21}{41}{37.94}{62.19}{42}
\emline{37.94}{62.19}{43}{37.94}{62.28}{44}
\emline{159.00}{51.00}{45}{82.97}{31.95}{46}
\bezier{228}(85.04,43.54)(59.29,24.05)(83.47,25.33)
\put(69.55,7.76){\circle*{3.00}}
\put(55.67,25.67){\circle*{3.00}}
\put(56.00,20.00){\makebox(0,0)[cc]{0}}
\put(69.67,2.00){\makebox(0,0)[cc]{$z$}}
\put(115.67,19.66){\makebox(0,0)[cc]{(\,$\rho$\,)}}
\put(18.67,1.66){\makebox(0,0)[cc]{($\lambda$)}}
\put(18.33,19.00){\makebox(0,0)[cc]{$\lambda_{\alpha,j}$}}
\emline{123.65}{26.14}{47}{123.65}{24.38}{48}
\emline{123.65}{24.38}{49}{127.44}{25.40}{50}
\emline{127.44}{25.40}{51}{123.51}{26.14}{52}
\emline{123.65}{25.74}{53}{127.17}{25.40}{54}
\emline{127.17}{25.40}{55}{127.10}{25.40}{56}
\emline{127.10}{25.40}{57}{123.58}{24.86}{58}
\emline{123.58}{24.86}{59}{123.65}{25.06}{60}
\emline{123.65}{25.06}{61}{126.76}{25.33}{62}
\emline{126.76}{25.33}{63}{123.58}{25.53}{64}
\emline{22.32}{8.48}{65}{22.32}{6.72}{66}
\emline{22.32}{6.72}{67}{26.10}{7.73}{68}
\emline{26.10}{7.73}{69}{22.18}{8.48}{70}
\emline{22.32}{8.07}{71}{25.83}{7.73}{72}
\emline{25.83}{7.73}{73}{25.77}{7.73}{74}
\emline{25.77}{7.73}{75}{22.25}{7.19}{76}
\emline{22.25}{7.19}{77}{22.32}{7.39}{78}
\emline{22.32}{7.39}{79}{25.43}{7.66}{80}
\emline{25.43}{7.66}{81}{22.25}{7.87}{82}
\emline{158.00}{95.33}{83}{158.00}{95.33}{84}
\emline{158.00}{95.33}{85}{158.00}{95.33}{86}
\emline{78.18}{95.36}{87}{78.18}{95.36}{88}
\emline{67.82}{62.40}{89}{78.52}{95.54}{90}
\end{picture}
\caption{
Deformation of the integration path in the variable $\rho$.  The
integration paths in $\rho$ and $\lambda$ are denoted by letters
in brackets. The cross ``$\times$'' denotes the eigenvalues
$\lambda_{\alpha,j}$ of $h_\alpha$ on the negative half-axis of
the physical sheet and the pair resonances belonging to the set
$\sigma_{\rm res}^{(\alpha)}$ of the sheet $\Pi_l$,\,\,
$l_0=+1$. Also, the cuts on $\Pi_l$,\,\, $l_0=+1$ beginning at
the points of $\sigma_{\rm res}^{(\alpha)}$ are shown in the
figure.
}
\label{figContour-mu}
\end{figure}

Besides, this is connected with the possibility when taking $z$
from $\Pi_0$ to $\Pi_l$, $l_0=+1$ to make a necessary
deformation of the integration path in $\lambda$ in such a way
that this path is separated from the integration contour in
$\rho$.

To obtain the representation (\ref{tr0l}) at a concrete point
$z=z_0$ we choose special final locations of the integration
paths in the variables $\lambda$ and $\rho$ after their consistent
deformation (see Fig.~\ref{figContour-fin}).


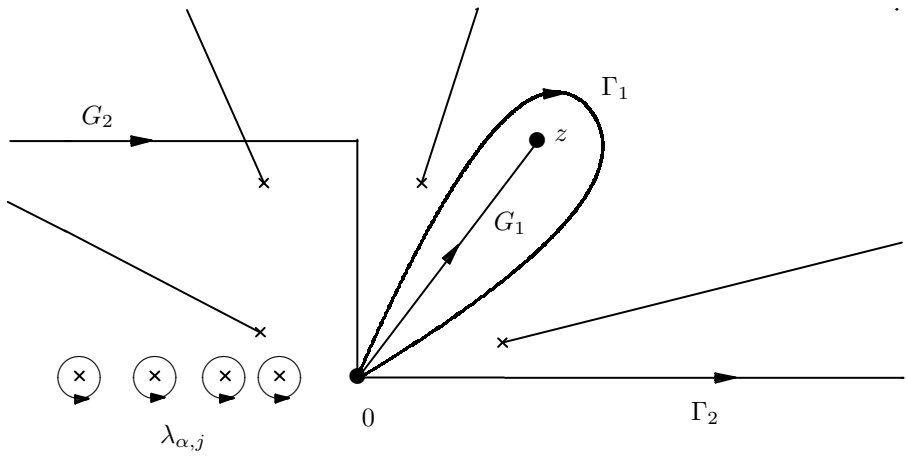
\begin{figure}
\centering
\unitlength=0.70mm
\special{em:linewidth .75pt}
\linethickness{.75pt}
\begin{picture}(159.33,95.54)
\emline{2.01}{24.85}{1}{3.66}{26.50}{2}
\emline{1.96}{26.50}{3}{3.72}{24.83}{4}
\emline{29.68}{24.85}{5}{31.33}{26.50}{6}
\emline{29.63}{26.50}{7}{31.39}{24.83}{8}
\emline{40.01}{24.85}{9}{41.66}{26.50}{10}
\emline{39.96}{26.50}{11}{41.72}{24.83}{12}
\emline{36.34}{33.18}{13}{37.99}{34.83}{14}
\emline{36.30}{34.83}{15}{38.06}{33.16}{16}
\emline{16.31}{24.87}{17}{17.96}{26.51}{18}
\emline{16.26}{26.51}{19}{18.02}{24.84}{20}
\emline{83.67}{25.33}{21}{159.33}{25.33}{22}
\emline{159.33}{25.33}{23}{159.00}{25.33}{24}
\emline{82.34}{31.18}{25}{83.99}{32.83}{26}
\emline{82.30}{32.83}{27}{84.06}{31.16}{28}
\emline{67.01}{61.52}{29}{68.66}{63.17}{30}
\emline{66.96}{63.17}{31}{68.72}{61.50}{32}
\emline{37.01}{61.52}{33}{38.66}{63.17}{34}
\emline{36.96}{63.16}{35}{38.72}{61.50}{36}
\emline{0.00}{53.33}{37}{37.00}{34.00}{38}
\emline{23.26}{95.21}{39}{37.94}{62.19}{40}
\emline{37.94}{62.19}{41}{37.94}{62.28}{42}
\emline{159.00}{51.00}{43}{82.97}{31.95}{44}
\put(55.67,25.67){\circle*{3.00}}
\put(57.67,17.67){\makebox(0,0)[cc]{0}}
\put(121.67,18.66){\makebox(0,0)[cc]{$\Gamma_2$}}
\put(22.60,13.90){\makebox(0,0)[cc]{$\lambda_{\alpha,j}$}}
\emline{123.65}{26.14}{45}{123.65}{24.38}{46}
\emline{123.65}{24.38}{47}{127.44}{25.40}{48}
\emline{127.44}{25.40}{49}{123.51}{26.14}{50}
\emline{123.65}{25.74}{51}{127.17}{25.40}{52}
\emline{127.17}{25.40}{53}{127.10}{25.40}{54}
\emline{127.10}{25.40}{55}{123.58}{24.86}{56}
\emline{123.58}{24.86}{57}{123.65}{25.06}{58}
\emline{123.65}{25.06}{59}{126.76}{25.33}{60}
\emline{126.76}{25.33}{61}{123.58}{25.53}{62}
\emline{158.00}{95.33}{63}{158.00}{95.33}{64}
\emline{158.00}{95.33}{65}{158.00}{95.33}{66}
\emline{78.18}{95.36}{67}{78.18}{95.36}{68}
\emline{67.82}{62.40}{69}{78.52}{95.54}{70}
\bezier{348}(55.43,25.51)(82.83,88.75)(97.70,78.10)
\bezier{376}(55.59,25.19)(115.70,61.20)(98.01,77.95)
\emline{55.28}{25.35}{71}{83.45}{25.35}{72}
\emline{55.75}{25.35}{73}{90.03}{70.43}{74}
\put(89.79,70.60){\circle*{3.00}}
\emline{55.59}{70.59}{75}{55.59}{25.35}{76}
\emline{55.67}{70.33}{77}{0.00}{70.33}{78}
\emline{91.07}{80.19}{79}{91.07}{78.43}{80}
\emline{91.07}{78.43}{81}{94.86}{79.45}{82}
\emline{94.86}{79.45}{83}{90.93}{80.19}{84}
\emline{91.07}{79.79}{85}{94.59}{79.45}{86}
\emline{94.59}{79.45}{87}{94.52}{79.45}{88}
\emline{94.52}{79.45}{89}{91.00}{78.91}{90}
\emline{91.00}{78.91}{91}{91.07}{79.11}{92}
\emline{91.07}{79.11}{93}{94.18}{79.38}{94}
\emline{94.18}{79.38}{95}{91.00}{79.58}{96}
\emline{91.10}{78.65}{97}{94.77}{79.46}{98}
\emline{94.77}{79.46}{99}{91.04}{80.01}{100}
\emline{12.73}{71.19}{101}{12.73}{69.43}{102}
\emline{12.73}{69.43}{103}{16.52}{70.45}{104}
\emline{16.52}{70.45}{105}{12.59}{71.19}{106}
\emline{12.73}{70.79}{107}{16.25}{70.45}{108}
\emline{16.25}{70.45}{109}{16.18}{70.45}{110}
\emline{16.18}{70.45}{111}{12.66}{69.91}{112}
\emline{12.66}{69.91}{113}{12.73}{70.11}{114}
\emline{12.73}{70.11}{115}{15.84}{70.38}{116}
\emline{15.84}{70.38}{117}{12.66}{70.58}{118}
\emline{12.77}{69.65}{119}{16.44}{70.46}{120}
\emline{16.44}{70.46}{121}{12.71}{71.01}{122}
\put(40.78,25.42){\circle{8.75}}
\put(30.29,25.42){\circle{8.75}}
\put(17.17,25.42){\circle{8.75}}
\put(2.74,25.42){\circle{8.75}}
\emline{2.01}{22.02}{123}{2.01}{20.42}{124}
\emline{2.01}{20.42}{125}{4.51}{21.36}{126}
\emline{4.51}{21.36}{127}{2.01}{21.98}{128}
\emline{1.97}{21.71}{129}{4.51}{21.32}{130}
\emline{4.51}{21.32}{131}{2.01}{20.69}{132}
\emline{2.01}{20.97}{133}{4.51}{21.28}{134}
\emline{1.97}{21.40}{135}{4.54}{21.32}{136}
\emline{16.44}{22.31}{137}{16.44}{20.71}{138}
\emline{16.44}{20.71}{139}{18.93}{21.65}{140}
\emline{18.93}{21.65}{141}{16.44}{22.27}{142}
\emline{16.40}{22.00}{143}{18.93}{21.61}{144}
\emline{18.93}{21.61}{145}{16.44}{20.98}{146}
\emline{16.44}{21.26}{147}{18.93}{21.57}{148}
\emline{16.40}{21.69}{149}{18.97}{21.61}{150}
\emline{29.55}{22.17}{151}{29.55}{20.57}{152}
\emline{29.55}{20.57}{153}{32.05}{21.50}{154}
\emline{32.05}{21.50}{155}{29.55}{22.13}{156}
\emline{29.51}{21.85}{157}{32.05}{21.46}{158}
\emline{32.05}{21.46}{159}{29.55}{20.84}{160}
\emline{29.55}{21.11}{161}{32.05}{21.42}{162}
\emline{29.51}{21.54}{163}{32.09}{21.46}{164}
\emline{40.05}{22.17}{165}{40.05}{20.57}{166}
\emline{40.05}{20.57}{167}{42.55}{21.50}{168}
\emline{42.55}{21.50}{169}{40.05}{22.13}{170}
\emline{40.01}{21.85}{171}{42.55}{21.46}{172}
\emline{42.55}{21.46}{173}{40.05}{20.84}{174}
\emline{40.05}{21.11}{175}{42.55}{21.42}{176}
\emline{40.01}{21.54}{177}{42.58}{21.46}{178}
\emline{71.84}{48.08}{179}{73.16}{47.05}{180}
\emline{73.16}{47.05}{181}{74.85}{50.59}{182}
\emline{74.85}{50.59}{183}{71.80}{48.04}{184}
\emline{72.05}{47.92}{185}{74.85}{50.55}{186}
\emline{74.85}{50.55}{187}{72.91}{47.18}{188}
\emline{72.71}{47.34}{189}{74.68}{50.43}{190}
\emline{72.25}{47.75}{191}{74.68}{50.35}{192}
\emline{-10.33}{70.33}{193}{1.33}{70.33}{194}
\emline{0.41}{53.05}{195}{-10.78}{58.73}{196}
\put(6.33,75.00){\makebox(0,0)[cc]{$G_2$}}
\put(84.67,53.33){\makebox(0,0)[cb]{$G_1$}}
\put(94.33,71.33){\makebox(0,0)[cc]{$z$}}
\put(104.67,80.33){\makebox(0,0)[cc]{$\Gamma_1$}}
\end{picture}
\caption{
Final location of the integration paths in the variables
$\rho$\,\, $(\Gamma_1\cup\Gamma_2)$ and $\lambda$ \,\,
$(G_1\cup G_2)$.  The path  $\Gamma_1$ represents a loop
going clockwise around the path $G_1$, the line segment
$[0,\,\, z]$;\, $\Gamma_2=[0,\,+\infty)$; \, $G_2=(z-\infty,\,
i\,{\rm Im}\, z]\cup[i\,{\rm Im}\, z,\, 0)$.
}
\label{figContour-fin}
\end{figure}

The singularity of the inner integral (in the variable $\rho$)
remains integrable after such a deformation
due to the presence of the factor  $\sqrt{\rho}$.
As a whole, the integral (\ref{ftr0f1}) becomes
\begin{small}
\be
\label{formpf}
\begin{array}{cl}
\vs{0.2cm}& \!\! \ds \Frac{1}{4}\Int_\Rt dk\Int_{S^2} d\hk'\Int_{S^2}
d\hp
\\
\vs{0.2cm}& \!\! \times \;\;\left\{\,
\Int_{G_1} d\lambda\sqrt{z-\lambda} \Int_{\Gamma_1\cup\Gamma_2}
d\rho\sqrt{\rho}\,
\Frac{t'\ad\big(k,\sqrt{\rho}\,\hk',\lambda\big)}{\rho-\lambda}\,
\tilde{f}\big(k,\sqrt{\rho}\,\hk',\sqrt{z-\lambda}\,\hp\big)  \right.\\
\vs{0.2cm} & \!\! + \;\;\Int_{G_2} d\lambda\sqrt{z-\lambda}
\Int_{\Gamma_1\cup  \Gamma_2}
d\rho\sqrt{\rho}\,
\Frac{t\ad\big(k,\sqrt{\rho}\,\hk',\lambda\big)}{\rho-\lambda}\,
\tilde{f}\big(k,\sqrt{\rho}\,\hk',\sqrt{z-\lambda}\,\hp\big)  \\
& \!\!  \left. + \;\;\Sum_{j=1}^{n\ad}\,2\pi i \,\sqrt{z-\lambda\ajd}
\Int_{0}^{+\infty}d\rho\sqrt{\rho}\,
\Frac{\phi\ajd(k)\,\overline{\phi}\ajd(k')}{\rho-\lambda\ajd}\,
\tilde{f}\big(k, \sqrt{\rho} \,\hk',\sqrt{z-\lambda\ajd}\,\hp\big)
\!\right\}
\end{array}
\ee
\end{small}
where $t'\ad$ denotes the pair T\,--\,matrix $t\ad(z)$ continued
into the unphysical sheet (with respect to $t\ad(\lambda)$ the path
$G_1$, $\lambda\in G_1$ just belongs to this sheet). The last term
arises as a result of taking residues at the points
$\lambda\ajd\in\sigma_d(h\ad)$.

Evidently, a domain of the variable $z\in\Pi_l$, $l_0=+1$,
where one can continue the function (\ref{ftr0f}) analytically
in the form  (\ref{formpf}) is determined by the conditions
$\Gamma_1\subset\cP_b$ and $\Gamma_1\cap Z\au_{\rm
res}=\emptyset$. These conditions can only be satified for
$z\in\cP_b$.

Note that the values of the inner integrals along $\Gamma_1$ are
determined for $\lambda\in G_1$ by the residues at the points
$\rho=\lambda$. At the same time
$\int_{G_2}d\lambda\,\ldots\int_{\Gamma_1}\,\ldots=0$, since at $\lambda\in
G_2$ the functions under the integration sign are holomorphic in
$\rho\in{\rm Int\,}\Gamma_1$.  Therefore,
\begin{eqnarray}
& & \!\!\reduction{(f,\bt\ad R_0(z)f')}{z\in\Pi_l,\, l_0=+1}
\nonumber          \\ [0.1cm]
& = &
\Frac{1}{4}\int_\Rt dk\int_{S^2} d\hk'
\int_{S^2} d\hp\, \nonumber  \\ [0.1cm]
 & & \!\! \times \;\;\left\{\,
\Int_{G_1} d\lambda\sqrt{z-\lambda}\, (-2\pi i) \sqrt{\lambda}\,
t'\ad\big(k,\sqrt{\lambda}\,\hk',\lambda\big)
\tilde{f}\big(k,\sqrt{\lambda}\,\hk',\sqrt{z-\lambda}\,\hp\big) \right.
\nonumber \\ [0.1cm]
 & & \!\!+ \;\;\Int_{G_1} d\lambda\sqrt{z-\lambda} \Int_{\Gamma_2}
d\rho\sqrt{\rho}  \\ [0.1cm]
 & & \!\!\nonumber
\times \;\; \Frac{t\ad\big(k,\sqrt{\rho}\,\hk',\lambda\big)+\pi i
\sqrt{\lambda}\,\tau\ad\big(k,\sqrt{\rho}\,\hk',\lambda\big)}{\rho-\lambda}
\,\tilde{f}\big(k,\sqrt{\rho}\,\hk',\sqrt{z-\lambda}\,\hp\big)
\nonumber\\ [0.1cm]
 & & \!\! + \;\; \Int_{G_2} d\lambda\sqrt{z-\lambda} \Int_{\Gamma_2}
d\rho\sqrt{\rho}\,
\Frac{t\ad\big(k,\sqrt{\rho}\,\hk',\lambda\big)}{\rho-\lambda}
\,\tilde{f}\big(k,\sqrt{\rho}\,\hk',\sqrt{z-\lambda}\,\hp\big)
\label{tR0fin}  \nonumber \\ [0.1cm]
 & & \!\!\left. + \;\; \Sum_{j=1}^{n\ad}\,2\pi i \sqrt{z-\lambda\ajd}
\Int_{0}^{+\infty}d\rho\sqrt{\rho}\,
\Frac{\phi\ajd(k)\,\overline{\phi}\ajd(k')}{\rho-\lambda\ajd}\,
\tilde{f}\big(k,\sqrt{\rho}\,\hk',\sqrt{z-\lambda\ajd}\,\hp\big)
\right\}. \nonumber
\end{eqnarray}
In the second summand of Eq.~(\ref{tR0fin}) we have used the
representation (\ref{3tnf}).

Joining the summands involving $t\ad$ on the physical sheet in
a separate integral $\int_{G_1\cup G_2}\,\ldots$  and then using
the holomorphy property of the integrand in $\lambda$  we
straighten the path $G_1\cup G_2$ turning it into the ray
$(z-\infty,z]$. As a result we get the bilinear form
corresponding to the product $\big(\bt\ad R_0\big)(z)$ taken
in the physical sheet.

The last term of the expression (\ref{tR0fin}) corresponds to
the kernel of the product \linebreak
$-\,\Phi\ad\rJ^{(\an)\dagger}L\au A\au\rJ\au \Phi^{*}R_0$.

Returning the reminding summands involving $t'\ad$
and $\tau\ad$ to the initial variables $k'$, $p'$  and utilizing then
the definition (\ref{spair6}), we find these summands
correspond to the expression
$$
L_0 A_0\!\!\left[\bt\ad -L_0 A_0 \rJot
\bs\adl^{-1}\rJo\bt\ad\right]\!\rJot\rJo
\,-\, L_0 A_0 \bt\ad\rJot\bs\adl^{-1}\rJo\bt\ad R_0\,.
$$
Gathering the results obtained we reveal that the analytical
continuation of $\bt\ad R_0$ into the sheet $\Pi_l,$ $l_0=+1$,
reads
\begin{eqnarray} \hs{1.0cm}
\reduction{[\bt\ad R_0(z)]}{\Pi_l}
 \!\! & = & \!\! \big(\bt\ad-L_0 A_0\bt\ad\rJot\bs\adl^{-1}\rJo\bt\ad-
\Phi\ad\rJ^{(\an)\dagger}L\au A\au \rJ\au \Phi^{*}\ad\big)
\nonumber\\ [-0.2cm]
& & \\ [-0.2cm]
\label{polufinal}
 & & \!\! \times \;\;\big(R_0+L_0 A_0\rJot\rJo\big) +
    L_0 A_0 \Phi\ad\rJ^{(\an)\dagger}
L\au A\au \rJ\au \Phi^{*}\ad\rJot\rJo\,. \nonumber
\end{eqnarray}
To be convinced of the factorization (\ref{tr0l}) it suffices
to observe that the last summand of (\ref{polufinal})
is equal to zero.  Indeed, for $\Img z\neq 0$ or $\Img z=0$ and
$z>\max_{\,j}{ \lambda\ajd}$ the following equalities hold
\be
\label{JFJ}
\big(\rJ\au\Phi\ad^{*}\rJot\big)(z)\EQ 0\,,\quad
\big(\rJo\Phi\ad\rJ\aut\big)(z) \EQ 0\,.
\ee
To prove, say, the first of them one can consider the $j$\,--\,th
component of the matrix\,--\,column $\rJ\au\Phi^{*}\rJot$,
$\big(\rJ\au\Phi^{*}\rJot\big)_j(z)=\rJ\ajd(z)
\langle\,\cdot\,,\phi\ajd\rangle\rJot(z).$ This component acts
on  $f\in L_2\big(S^5\big)$ as follows
\begin{eqnarray*}
& & \big(\big(\rJ\au\Phi^{*}\rJot\big)_{j}f\big)\big(\hk\ad,z\big) \\ [0.1cm]
&  = & \!\!\Int \!dk\ad''\,\overline{\phi}\ad\big(k\ad''\big)
\int \!d\hP' \delta\big(k\ad''-\sqrt{z}\cos\omega\ad'\hk\ad'\big)
\delta\big(\sqrt{z-\lambda\ajd}\,\hp\ad-\sqrt{z}\sin\omega\ad'\hp\ad'\big)
f\big(\hP'\hs{0.02cm}\big)
\end{eqnarray*}
where we use again the hyperspherical coordinates, $\hP'\sim
\big(\omega\ad', \hk\ad',\hp\ad'\big).$ It is clear that only the points
$\hP'\!\in \!S^5$ with
$\sqrt{z}\sin\omega\ad'\hp\ad'=\sqrt{z-\lambda\ajd}\,\hp\ad$ may
give a nontrivial contribution to the integral.  The last
equality means that the condition
$z-\lambda\ajd=$ $= z\sin^2\omega\ad'$ has to be satisfied for some
$\omega\ad'\in[0,\, {\pi}/{2}]$. This condition is
equivalent to the requirement $z\cos^2\omega\ad=\lambda\ajd$
which may be obeyed for real $z\leq \lambda\ajd$ only.
However, by the definition of the surface $\Re$ such $z$ do not
belong to the sheets $\Pi_l$, $l_0=\pm 1$ (see
Sec.$\:$ 5).  Consequently,
$\big(\big(\rJ\au\Phi^{*}\rJot\big)_{j}f\big)\big(\hk\ad,z\big)=0$ for any
$j$ and the
first equality (\ref{JFJ}) holds. The second one
may be proved analogously.

As already mentioned above, it follows from Eq.~(\ref{JFJ}) that
the last summand \linebreak
of (\ref{polufinal})  disappears and hence,
Eq.~(\ref{tr0l}) is true.  This completes the proof of the \linebreak
theorem.
\end{proof}

Using Eqs.~(\ref{tlRbig}) and (\ref{tr0l}) one can present the
Faddeev equations (\ref{MFE}) continued into the sheet $\Pi_l$ in
the matrix form
\be
\label{MFEl}
M^l(z) \EQ \bt^l(z)\,-\,\bt^l(z)\hs{0.02cm}\bRo^l(z)\hs{0.02cm}\Y \hs{0.02cm}
M^l(z)
\ee
where
\vs{0.2cm}\begin{eqnarray}
\label{tltot}
\bt^l(z) \!\! & = & \!\!\bt\,-\,
L_0A_0 \bt\bJot\bs_l^{-1}\bJo\bt\,-\,\Phi\rJt_1 L_1 A_1\rJ_1\Phi^{*},\\ [0.2cm]
\label{R0ltot}
\bRo^l(z) \!\! & = & \!\!\bRo(z)\,+\,L_0 A_0(z)\bJot(z)\bJo(z)\,.
\end{eqnarray}
Here,
 $\bs_l(z)=\diag\{\bs_{1,l}(z),\bs_{2,l}(z),\bs_{3,l}(z)\}$.  By
$M^l(z)$ we understand a supposed analytic continuation
of the matrix $M(z)$ into the sheet $\Pi_l$.
\begin{lemma}\label{LPossibility}
For each two\,--\,body unphysical sheet $\Pi_l$ of the surface $\Re$
there exists a path from the physical sheet $\Pi_0$ to the
domain $\Pi_l^{\rm (hol)}$ of $\Pi_l$
which only passes trough two\,--\,body
unphysical sheets $\Pi_{l'}$ and, moving on this path, the
parameter $z$ always stays in the respective domains
$\Pi_{l'}^{\rm (hol)}\subset\Pi_{l'}$.
\end{lemma}

\begin{proof}
Let us make a use of the principle
of mathematical induction.  To
this end, we rearrange the branching points
$\lambda\ajd,$ $\alpha =  1,2, 3,  \;j= 1,2,\,\ldots,n_{\alpha}$, in
nondecreasing order redenoting them
as $\lambda_1,\lambda_2,\,\ldots,\lambda_m$, $m\leq\sum_{\an} n\ad,$
$\lambda_1 <\lambda_2 <$ $< \cdots< \lambda_m$, and putting
$\lambda_{m+1}=0.$ Let the multi\,--\,index $l=(l_1,l_2,\,\ldots,l_m)$
correspond temporarily namely to this enumeration.  As
previously, $l_j=0$ if the sheet $\Pi_l$ is related to the main
branch of the function $(z-\lambda_j)^{1/2}$ otherwise $l_j=1$.  The
index $l_0$ is omitted in these temporary notations.

It is clear that the transition of $z$ from the physical sheet
$\Pi_0$ through the segment $(\lambda_1,\lambda_2)$ in the
neighboring unphysical sheet $\Pi_{l^{(1)}}$ $\big($into
$\Pi_{l^{(1)}}^{\rm hol}\big)$,
$l^{(1)}=\Big(l_1^{(1)},l_2^{(1)},\,\ldots,l_m^{(1)}\Big)$ with $l_1^{(1)}=1$
and $l_j^{(1)}=0$ for $j\neq 1$ is possible by definition  of
the domain $\Pi_{l^{(1)}}^{\rm hol}$ (see  Sec.$\:$4).
According to Lemmas~\ref{LEqParab} and~\ref{LEqQuadr}, if $z$
belongs to $\Pi_{l^{(1)}}^{\rm (hol)}$, it may be led to the
real axis in the interval $\big(\lambda^{(1)},+\infty\big)$ with certain
$\lambda^{(1)} <\lambda_1.$ Remaining in
$\Pi_{l^{(1)}}^{\rm(hol)}$, the point $z$ may even go around the
threshold $\lambda_1$ crossing the real axis in the segment
$\big(\lambda^{(1)}, \lambda_1\big). $ Thus, the parameter $z$ may be
led from the sheet $\Pi_{l^{(1)}}$ into each neighboring
unphysical sheet and, in particular, into the sheet $\Pi_l$ identified
by  $l_1=0,$ $l_2=1$, $l_j=0,$ $j\geq 3.$ Transition of $z$ from
$\Pi_0$ through the segment $(\lambda_2,\lambda_3)$ into the sheet
$\Pi_l$ with $l_1=l_2=1,$    $l_j=0,$   $j\geq 3,$ is always
possible.

We suppose further that the parameter $z$ may be carried in this
manner from $\Pi_0$ into all the two\,--\,body unphysical sheets
$\Pi_{l^{(k)}}$ determined by the conditions $l_j^{(k)}=0,$
$j>k. $ It is assumed also that during this motion $z$ always
remains in the domains $\Pi_{l^{(k)}}^{\rm (hol)}$ of these
sheets and does not visit other sheets.  It follows from
Lemmas~\ref{LEqParab} and~\ref{LEqQuadr} that if $z$ stays in
the domain $\Pi_{l^{(k)}}^{\rm (hol)}$ of the sheet
described, then it can be led to the real axis in the segment
$\big(\lambda^{(k)},+\infty\big)$ with a certain
$\lambda^{(k)}<\lambda_k.$ Hence, the parameter $z$ from each of
the sheets $\Pi_{l^{(k)}}$ may be carried through the interval
$(\lambda_k,\lambda_{k+1})$ into the neighboring unphysical sheet
$\Pi_{l^{(k+1)}}$ with $l_j^{(k+1)}=1- l_j^{(k)}$, $j\leq k,$
$l_{k+1}^{(k+1)}=1$  and  $l_j^{(k+1)}=0,$  $j>k+1$.  This just means
that  $z$ may be carried from $\Pi_0$ into all the
two\,--\,body unphysical sheets $\Pi_{l^{(k+1)}}$ with
$l_j^{(k+1)}=0,$   $j >k+1$. The whole time the parameter $z$
remains in the holomorphy domains $\Pi_{l^{(k+1)}}^{\rm (hol)}$
and does not visit the sheets $\Pi_{l^{(s)}}$ with $s>k+1$. By
the principle of mathematical induction we conclude that the
parameter $z$ may be carried into all the two\,--\,body
unphysical sheets.

The proof is completed.%
\end{proof}

Using results of Sec. 4, as well as
Lemma~\ref{LPossibility}, we can prove the following important
statement.

\begin{theorem}\label{ThIter}
The iterations $\cQ\un(z)=\bigl((-\bt\bRo\Y)^n\bt\bigr)(z)$ for
$n\geq 1$  of the absolute terms in the Faddeev
equations~{\rm(\ref{MFE})} admit analytic continuation on the
domain $\Pilh$ of each unphysical sheet $\Pi_l\subset\Re$ in
the sense of distributions over $\Osix$. This continuation is
described by the equalities
$
%
  \reduction{\cQ\un(z)}{\Pi_l}=\bigl((-\bt^l\bRo^l\Y)^n \bt^l \bigr)(z).
$
\end{theorem}
\begin{remark}\label{NIterF}
{\rm The products  $L_1\rJ_1\Psis\Y\cQ^{(m)}$,
$\cQ^{(m)}\Y\Psi\rJt_1 L_1$,   $\tilde{L}_0\bJo\cQ^{(m)}$,
$\cQ^{(m)}\bJot\tilde{L}_0$,
$L_1\rJ_1\Psis\Y \cQ^{(m)}\Y\Psi\rJt_1 L_1$,
$\tilde{L}_0\bJo\cQ^{(m)}\bJot\tilde{L}_0,$
$L_1\rJ_1\Psis\Y\cQ^{(m)}\bJot\tilde{L}_0$ and
$\tilde{L}_0\bJo\cQ^{(m)}\Y\Psi\rJt_1 L_1$,
$0\leq m < n$,
arising after substitution of the relations {\rm (\ref{tltot})}
and {\rm (\ref{R0ltot})}  into
$ \reduction{\cQ\un(z)}{\Pi_l}$, have to be understood
in the sense of the definitions from
Sec. 4.}
\end{remark}

\begin{proof}
Theorem~\ref{ThIter} will be proved in the case of the analytic
continuation of the iteration $\cQ^{(1)}(z)$. It will be clear
from this proof that the iterations $\cQ\un(z)$ with $n\geq 2$
could be considered in the same way as $\cQ^{(1)}(z)$. We shall
not expound here on the correspoding computations for $n\geq 2$, since
they are too cumbersome.

So, let us consider the bilinear forms
\begin{eqnarray}
Q\abd(z)\!\! & = & \!\!(f,\bt\ad(z) R_0(z)\bt\bd(z)f')\nonumber\\ [0.1cm]
&  = & \!\! \Frac{1}{|s\abd|}\,\ds \int_\Rt dk\ad \int_\Rt dp\ad
\int_\Rt dk'\bd \int_\Rt dp'\bd \,   f(P)f'(P')    \nonumber \\ [0.1cm]
& & \!\! \times \;\; \Frac{t\ad\big(k\ad,\tk\bu\ad\big(p\ad,p\bd'\big),z-
p\ad^2\big) t\bd\big(\tk\bd\au\big(p\ad,p'\bd\big),k'\bd,z-{p'}\bd^{2}\big) }
{p\ad^{2}+{p'}\bd^2- 2 c\abd\big(p\ad,p'\bd\big)-s\abd^2 z}
         \nonumber
\end{eqnarray}
corresponding to the components $\cQ\abd^{(1)}(z)=-\bt\ad(z)
R_0(z)\bt\bd(z)$,\,\, $\bn\neq\an$, of the iteration
$\cQ^{(1)}(z)$, $\Img z\neq 0$.  It is assumed that
$f,f'\in \Osix$.

Using the spherical coordinates  $p\ad\rightarrow\rho=|p\ad|^2,$
$\hat{p}\ad,$\, $p'\bd\rightarrow\rho'=|p'\bd|^2,$
$\hat{p}'\bd$, in the integrals in the variables $p\ad$ and
$p'\bd$ we get
\begin{eqnarray} \hs{1.0cm}
& & \!\!  Q\abd(z) \nonumber \\ [0.1cm]
& = & \!\! \Frac{1}{4}\cdot\Frac{1}{|s\abd|}
\int_\Rt dk\ad \int_\Rt dk'\bd
\int_\Stwo d\hat{p}\ad \int_\Stwo d\hat{p}'\bd
\int_0^{\infty} d\rho\sqrt{\rho} \int_0^{\infty} d\rho'\sqrt{\rho'}
                               \nonumber \\ [-0.2cm]
& & \!\!\times \;\;f\big(k\ad,\sqrt{\rho}\,\hat{p}\ad\big)\,
f'\big(k'\bd,\sqrt{\rho'}\,\hat{p}'\bd\big)
                   \label{Qab}     \\
& & \!\! \times \;\;
\Frac{t\ad\big(k\ad,\tk\ad\bu\big(\sqrt{\rho}\,\hat{p}\ad,\sqrt{\rho'}
\,\hat{p}\bd'\big),z-\rho\big)
t\bd\big(\tk\bd\au\big(\sqrt{\rho}\,\hat{p}\ad,
\sqrt{\rho'}\,\hat{p}\bd'\big),k'\bd,z-\rho'\big) }
{\rho+\rho- 2 c\abd\,\sqrt{\rho}\,\sqrt{\rho'}
\big(\hat{p}\ad,\hat{p}'\bd\big)-s\abd^2 z}\,.
\phantom{\Int_0^0}                                       \nonumber
\end{eqnarray}
Let us begin with continuation of the functions $Q\abd(z)$
across the cut $(\lambda_{\rm min},+\infty)$ to the left from
the three\,--\,body threshold $\lambda_0=0.$ We realize this
continuation in the same way as the continuation of the kernels of
$\bt\ad(z)$ and $(\bt\ad R_0)(z)$ in the proof of
Theorem~\ref{ThKernels}.  Additionally, we use the fact that for
$z<0$ the denominator of the expression under the integration
sign in (\ref{Qab}) cannot become equal to zero since
         $$
\rho\,+\,\rho'\,-\, 2 c\abd\,\sqrt{\rho}\,\sqrt{\rho'}\big(\hat{p}\ad,
\hat{p}'\bd\big) \GE (1-|c|_{\alpha \beta})(\rho + \rho') \quad
\fr \quad \rho,\,\rho' \GT 0\,, \;\;\hat{p}\ad, \,\hat{p}'\bd\in\Stwo\,.
          $$
Thus, when continuing across the segment $(\lambda_{\rm
min},0)$, only residues at the poles $\lambda\ajd,\,\lambda\bkd$
give a nontrivial contribution to the Cauchy type integrals in the
variables $\rho,\,\rho'$ generated in (\ref{Qab}) by the two\,--\,body
singularities (the poles of $g\ajd(z)$ [see Eq.~(\ref{tcd})]).

Let us continue the function (\ref{Qab}) across the segment
$(\lambda^{(1)},\lambda^{(2)})$ where
$\lambda^{(1)},\lambda^{(2)}, $
$\lambda^{(1)} < \lambda^{(2)},$
are some neighboring points of the set
$\sigma_d^{(2)}=\bigcup_{\an=1}^{3} \sigma_d(h\ad)$,
$\big(\lambda^{(1)}, \lambda^{(2)}\big) \cap$ $\cap \sigma_d^{(2)}
=\emptyset.$ Then the totality of the kernels  $\bt\ad(z) R_0(z)
\bt\bd(z)$, $\an,\,\bn=1,2,3,$ $\bn\neq\an$, may be
continued in $z$ into the two\,--\,body unphysical sheet
$\Pi_l\subset\Re$ neighboring the physical one and such
that its indices $l_0=0$, $l\gjd =1$ if
$\lambda\gjd>\lambda^{(1)}$, and $l\gjd=0$ if
$l\gjd\geq\lambda^{(2)}$, $\gn=1,2,3,$ $j=1,2,\,\ldots,n\gd$.  One
can easily check that the continuation of the functions $Q\abd(z)$
into these sheets $\big($in a vicinity of the segment
$\big(\lambda^{(1)},\lambda^{(2)}\big)\big)$ corresponds exactly to the
iteration $\bt^l(z)\bRo(z)\Y\bt^l(z)$ of the absolute term of
the continued Faddeev equations (\ref{MFEl}).  Recall that in
the two\,--\,body sheets $\bt^l(z)=\bt(z)-\Phi \rJt_1 L_1 A_1 \rJ_1
\Phi^{*}$.  Hence the analyticity domain of the kernels
$\bt^l(z) \bRo(z)\Y\bt^l(z)$ in these sheets is determined by the
set of those points where the functions
\begin{eqnarray*}
    F_{\an,j;\,\bn,k}(z,\eta) \!\! & = & \!\!
    \big(z-\lambda\ajd+z-\lambda\bkd-
 2 c\abd \sqrt{z-\lambda\ajd}\,\sqrt{z-\lambda\bjd}\,\eta-s\abd^2
 z\big)^{-1}\,, \\ [0.1cm]
    F_{\an,j}(z,\rho',\eta) \!\! & = & \!\!\big(z-\lambda\ajd+\rho'-
    2 c\abd \sqrt{z-\lambda\ajd}\,\sqrt{\rho'}\:\eta-s\abd^2 z\big)^{-1}
\,,  \\ [0.1cm]
    F_{\bn,k}(z,\rho,\eta)\!\! & = & \!\!\big(\rho+z-\lambda\bkd-
    2 c\abd \sqrt{\rho}\,\sqrt{z-\lambda\bjd}\:\eta-s\abd^2 z\big)^{-1}\,,
\end{eqnarray*}
are holomorphic in $z$ for all $\rho,\,\rho'>0$, $\eta\in [-1,1]$.
The latter arise in (\ref{Qab}) due to the presence of the factors
(\ref{Znam}) as a result of taking residues at the poles
$\rho=z-\lambda\ajd$ and/or $\rho'=z-\lambda\bkd,$
$\lambda\ajd, \,\lambda\bkd \leq\lambda^{(1)},$
$\an,\,\bn=1,2,3,$ $\an\neq\bn.$ The domains where the
singularities of the above functions are situated have been
described in Lemmas~\ref{LEqParab} and~\ref{LEqQuadr}.  It
follows from these lemmas that  the product $\bt^l(z)
\bRo(z)\Y\bt^l(z)$ describes the analytic continuation of the
iteration $\cQ^{(1)}(z)$  on the domain $\Pi_l^{\rm (hol)}$ of
each neighboring (with respect to $\Pi_0$) two\,--\,body sheet
$\Pi_l$.  Note that the singularities of the functions $
F_{\an,j;\,\bn,k}(z,\eta),$ $ F_{\an,j}(z,\rho',\eta)$ and $
F_{\bn,k}(z,\rho,\eta)$ are, as a matter of fact, three\,--\,body ones
though being situated in the two\,--\,body unphysical sheets.
Indeed, making the substitution
$\eta=(\hat{p}\ad,\hat{p}\bd')$
in (\ref{Qab})  one finds that the integral in
the variable $\eta$ turns out to be a Cauchy type integral.  This
means that, e.\,g., the points $z_{\rm rl}, z_{\rm rt}$ (see
Lemma~\ref{LEqQuadr}) are extra logarithmic branching points. After
crossing the cuts on $\Pi_l$ along the segments $[z_{\rm rl},
z_{\rm rt}]$ as well as the root ellipses from
Lemma~\ref{LEqQuadr}, the representation of the analytic
continuation of $\cQ^{(1)}(z)$ in the form of the product
$\bt^l\bRo^l\bt^l$ becomes invalid. In the present paper we
restrict ourselves to considering only those domains of the
unphysical sheets where the correctness of such representations is
not violated.

Let us now use Lemma~\ref{LPossibility} and carry out a
continuation of the form $Q\abd(z)$ into the rest of the two\,--\,body
unphysical sheets.  Boundaries of the holomorphy domains
$\Pi_l^{\rm (hol)}$ of $Q\abd(z)$ of these sheets are determined
again only by the indices $\an,j$ and $\bn,k$ of the functions
$F_{\an,j;\bn,k},$   $F_{\an,j}$ and $  F_{\bn,k}$ included in
the kernels of the operators
\begin{eqnarray*}
L_1\rJ_1\Phi^{*}\bRo\Y\Phi\rJt_1 L_1 \!\! & \equiv & \!\!
L_1 \rJ_1\Psis\Y\bv\Psi\rJt_1 L_1\,, \\ [0.1cm]
 L_1 \rJ_1 \Phis\bRo\Y\bt \!\! & \equiv & \!\!  L_1 \rJ_1 \Psis\Y\bt\,,
 \\ [0.1cm]
\bt\bRo\Y\Phi\rJt_1 L_1 \!\! & \equiv & \!\! \bt\Y\Psi\rJt_1 L_1
\end{eqnarray*}
arising in
the product $\bt^l\bRo^l\Y\bt^l$.
\par
Thus we can state that the kernels of the iteration
$\cQ^{(1)}(z)$ admit an immediate analytic continuation as
holomorphic generalized functions over $\Osix$ on the domains
$\Pi_l^{\rm (hol)}$ of all the two\,--\,body unphysical sheets where
$$\reduction{\cQ^{(1)}(z)}{\Pi_l} \EQ \bt^l(z)\bRo(z)\Y\bt^l(z)\,.$$

Let us consider now the continuation of the iteration
$\cQ^{(1)}(z)$ into the three\,--\,body unphysical sheets
$\Pi_l$ with $l_0=\pm 1$.  It is clear that the two\,--\,body
singularities will give the same contribution to the continued
kernels as before when continuing this iteration into the
two\,--\,body sheets. Therefore we assume here for the sake of
simplicity that these singularities are absent or, in other
words, that the two\,--\,body subsystems have no discrete
spectrum.

Let us deal, say, with the continuation of $\cQ^{(1)}(z)$ into
the sheet  $\Pi_l$ with $l_0=+1$. This means that we study a
crossing of the ray  $(0,+\infty)$ from below going upward. We
begin with taking in (\ref{Qab}) the limit $z\rightarrow
E-i\hs{0.03cm}0$, $E > 0,$ and rewriting the limit values of the
T\,--\,matrix $t\ad$ ($t\bd$) at the points
$E-\rho-i\hs{0.03cm}0$ $(E-\rho'-i\hs{0.03cm}0)$ of the segment
$(0,E)$ in terms of the continued kernel
$t\ad\big(k\ad,k'\ad,z\big)$
$\big(t\bd\big(k\bd,k'\bd,z\big)\big)$ on the unphysical sheet
using the representation (\ref{3tnf}). For $t\ad$ we have
\begin{eqnarray*}
& & \!\!t\ad\big(k\ad,\tk\bu\ad\big(\sqrt{\rho}\,\hat{p}\ad,
\sqrt{\rho'}\,\hat{p}'\bd\big),E-\rho-i\hs{0.03cm}0\big)\\ [0.1cm]
 \!\! &= & \!\!
t\ad\big(k\ad,\tk\bu\ad\big(\sqrt{\rho}\,\hat{p}\ad,
\sqrt{\rho'}\,\hat{p}'\bd\big),E-\rho+i\hs{0.03cm}0\big) \\ [0.1cm]
& & \!\! + \;\;
\pi i \,\sqrt{\rho+i\hs{0.03cm}0}\,
\tau\ad\big(k\ad,\tk\bu\ad\big(\sqrt{\rho}\,\hat{p}\ad,
\sqrt{\rho'}\,\hat{p}'\bd\big),E-\rho+i\hs{0.03cm}0\big)
\end{eqnarray*}
and analogously for $t\bd$.  At  $\rho>E$ $(\rho'> E)$ the limit
values of the T\,--\,matrix \linebreak
$t\ad(\,\ldots,E-\rho-i\hs{0.03cm}0)$
$\big(t\bd\big(\,\ldots,E-\rho-i\hs{0.03cm}0\big)\big)$ from below coincide
with the limit
values  $t\ad(\,\ldots,E-\rho+i\hs{0.03cm}0)$ $\big(t\bd\big(\,\ldots,
E-\rho+i\hs{0.03cm}0\big)\big)$ from
above, in view of analyticity of $t\ad(z)$ ($t\bd(z)$) in
$z$ at $z\not\in \R^{+}$.  In the same way we rewrite as well
the denominator of the expression under the integration sign in
(\ref{Qab}),
\begin{eqnarray*}
\Frac{1}{{\rm F}\big(\rho,\rho',
\hat{p}\ad,\hat{p}'\bd\big)-s\abd^2+i\hs{0.03cm}0}\!\! & = & \!\!
\Frac{1}{{\rm F}\big(\rho,\rho',
\hat{p}\ad,\hat{p}'\bd\big)  -  s\abd^2-i\hs{0.03cm}0} \\ [0.1cm]
& & \!\! - \;\;
2\pi i\delta\big({\rm F}\big(\rho,\rho',
\hat{p}\ad,\hat{p}'\bd\big)-s\abd^2 E\big)
\end{eqnarray*}
where ${\rm F}=\rho+\rho'-2c\abd\sqrt{\rho}\,\sqrt{\rho'}
\big(\hat{p}\ad,\hat{p}'\bd\big).$

Note right away that the kernel $\delta({\rm F}(\rho,\rho',\hat{p}\ad,\hat{p}'
\bd)-s\abd^2 E)$
corresponds to the distribution $\big(\rJot\rJo\big)(P,P',E\pm i\hs{0.03cm}0).$
Also, it is easy to find that all the terms of (\ref{Qab})
including the $\delta$\,--\,fun\-c\-tion $\delta\big({\rm
F}\big(\rho,\rho',\hat{p}\ad,\hat{p}'\bd\big)-s\abd^2 E\big)$ generate as a
sum a bilinear form corresponding to the kernel $l_0
A_0(E)\big(\bt\ad^l\rJot\rJo\bt\bd^l\big)(P,P',E+i\hs{0.03cm}0)$ admitting
analytic
continuation on $\Pi_l\cap\cP_b,$\,\, $l_0=1$ in the
sense of distributions over $\Osix$.

Then we consider the terms of (\ref{Qab}) including the factor
$1/\big({\rm F}\big(\rho,\rho',\hat{p}\ad,\hat{p}'\bd\big)-s\abd^2 E-
- i \hs{0.03cm}0\big).$
The simplest of the summands includes the fraction
$$
\Frac{1}{|s\abd|}\cdot
\Frac{t\ad(\,\ldots,E-\rho+i\hs{0.03cm}0)t\bd(\,\ldots,E-\rho'+i\hs{0.03cm}0)}
{\rho+\rho'-2c\abd\sqrt{\rho}\,\sqrt{\rho'}
\big(\hat{p}\ad,\hat{p}'\bd\big)-s\abd^2 E -i\hs{0.03cm}0}\,.
$$
In all this, the integration in $\rho$ as well as in
$\rho'$ is carried out along the interval $(0, +\infty)$.
Evidently, this summand represents a boundary value
at $z=E+i\hs{0.03cm}0$ of the bilinear form for the product
$\bt\ad(z) R_0(z)\bt\bd(z).$

We consider contributions of the summands which include
the products $\tau\ad(\,\ldots\,)t\bd(\,\ldots\,)$ and
$t\ad(\,\ldots\,)\tau\bd(\,\ldots\,)$ 
for the case of the first of such summands.

Let us rewrite the respective bilinear form,
\begin{equation}   \hs{0.5cm}
\begin{array}{lcl}
\vs{0.2cm}      Q\abd^{(t\tau)}(E) \!\! & =
 & \!\!\ds  \Frac{1}{4}\cdot \Frac{1}{|s\abd|}\int_\Rt dk\ad
\int_\Rt dk'\bd
\int_{S^2} d\hat{p}\ad \int_{S^2} d\hat{p}'\bd
\int_0^{E} d\rho\,\sqrt{\rho} \int_0^{\infty} d\rho'\sqrt{\rho'}  \\
\vs{0.2cm}& & \!\!  \times \;\;  \pi i\,\sqrt{\rho}
  f\big(k\ad,\sqrt{\rho}\,\hat{p}\ad\big)\,
  f'\big(k'\bd,\sqrt{\rho'}\hat{p}'\bd\big)
                 \label{Qabtt} \\
\vs{0.2cm}& & \!\! \times  \;\;  \Frac{\tau\ad\big(k\ad,
\tk\bu\ad\big(\sqrt{\rho}\,
\hat{p}\ad,\sqrt{\rho'}\,\hat{p}\bd'\big),E-\rho+i0\big)       }
{\rho+\rho- 2 c\abd\,\sqrt{\rho}\,
\sqrt{\rho'}\big(\hat{p}\ad,\hat{p}'\bd\big)-s\abd^2 E-i0}\,.\\
& & \!\! \times \;\; \Frac{t\bd\big(\tk\au\bd(\sqrt{\rho}\,\hat{p}\ad,
\sqrt{\rho'}\,\hat{p}\bd'),k'\bd,E-\rho'+i0\big)}
{\rho+\rho- 2 c\abd\,\sqrt{\rho}\,
\sqrt{\rho'}\big(\hat{p}\ad,\hat{p}'\bd\big)-s\abd^2 E-i0}
\end{array}
\end{equation}
We use Lemma~\ref{LEqQuadrPrime} to prove the existence of an
analytic continuation of the function $Q\abd^{(t\tau)}(E)$ in
the domain $\Img z>0$. To apply this lemma we
divide the interval of integration in the variable
$\rho'$  in (\ref{Qabtt}) into two intervals
$[0,E]$  and  $(E, +\infty)$. Then we
get in (\ref{Qabtt}) two terms including
$ \ldots\int_0^E d\rho \int_0^E d\rho'\,\ldots$
and $ \ldots\int_0^E d\rho \int_E^{+\infty} d\rho'\,\ldots\: .$
In the first term we make two changes of variables,
$\rho\rightarrow\nu,$ $\rho=\nu E,$
and $\rho'\rightarrow\nu',$ $\rho'=\nu' E,$
and in the second term, only the first change, $\rho=\nu E.$
As a result we find that
$$ Q\abd^{(t\tau)}(E) \EQ
 Q\abd^{(1)}(E+i\hs{0.03cm}0)\,+\, Q\abd^{(2)}(E+i\hs{0.03cm}0)$$
where by  $ Q\abd^{(1)}(E+i\hs{0.03cm}0)$ and $Q\abd^{(2)}(E+i\hs{0.03cm}0)$
we understand the boundary values (at $ z= E+i\hs{0.03cm}0,$ $ E > 0 $)
of the functions
$$
\begin{array}{lcl}
\vs{0.2cm}Q\abd^{(1)}(z)\!\! &  = & \!\! \ds
\Frac{1}{4}\cdot\Frac{(\sqrt{z})^5}{|s\abd|}
\int_\Rt dk\ad \int_\Rt dk'\bd
\int_{S^2} d\hat{p}\ad \int_{S^2} d\hat{p}'\bd
\int_0^{1} d\nu\sqrt{\nu} \int_0^{1} d\nu'\sqrt{\nu'}  \\
\vs{0.2cm}& & \!\! \times\;\; f\big(k\ad,\sqrt{z}\,\sqrt{\nu}\,\hat{p}\ad\big)
\cdot f'\big(k'\bd,\sqrt{z}\sqrt{\nu'}\hat{p}'\bd\big) \\
& & \!\! \times \;\;\Frac{       \pi i\,\sqrt{\nu}
\:\tau\ad(\,\ldots,z(1-\nu)) \cdot t\bd(\,\ldots, z(1-\nu'))       }
{\nu+\nu- 2 c\abd\,\sqrt{\nu}\,
\sqrt{\nu'}\big(\hat{p}\ad,\hat{p}'\bd\big)-s\abd^2 -i\hs{0.03cm}0}\,,
\end{array}
$$
and
\be
\label{Qab2}
\begin{array}{lcl}
Q\abd^{(2)}(z)\!\! & = & \!\! \ds
\Frac{1}{4}\cdot \Frac{z^2}{|s\abd|}\int_\Rt dk\ad \int_\Rt dk'\bd
\int_{S^2} d\hat{p}\ad \int_{S^2} d\hat{p}'\bd
\int_0^{1} d\nu\,\sqrt{\nu} \int_{\Gamma_z} d\rho'\,\sqrt{\rho'}  \\[0.3cm]
& & \!\!    \times \;\;f\big(k\ad,\sqrt{z}\,\sqrt{\nu}\,
\hat{p}\ad\big)\cdot
f'\big(k'\bd,\sqrt{\rho'}\,\hat{p}'\bd\big)  \phantom{\Int_0^0}   \\ [0.15cm]
& & \!\! \times \;\;    \Frac{       \pi i\,\sqrt{\nu}
\, \tau\ad(\,\ldots,z(1-\nu)) \cdot t\bd(\,\ldots,z-\rho')       }
{\nu z+\rho'-
2 c\abd\,\sqrt{\nu}\,\sqrt{z}\,\sqrt{\rho'}\big(\hat{p}\ad,\hat{p}'\bd\big)-
s\abd^2 z}\,,
\end{array}
\ee
respectively. For $\Gamma_z$ we take the same path
as in (\ref{Q2cont}).

By the same reasoning as in the consideration of the forms
(\ref{Q1cont}) and (\ref{Q2cont}), we conclude that the
functions $Q^{(1)}\abd(z)$ and $Q^{(2)}\abd(z)$ admit
continuations in the domain $\Img z\! >\! 0$.

Therefore, we have proved that the function  $Q\abd^{(t\tau)}(E)$
admits an analytic continuation on the domain $\Pilh$ of the sheet
$\Pi_l,$ $l_0=+1.$
Analogously, an analytic continuation on $\Pilh,$ $l_0=+1$
exists as well for the bilinear
form $Q\abd^{(\tau t)}(E)$ corresponding
to the contribution in (\ref{Qab})
from the product $t\ad(\,\ldots\,)\tau\bd(\,\ldots\,)$.

To this end, let us consider the contribution
of the product $\tau\ad(\,\ldots)\tau\bd(\,\ldots)$.
The respective bilinear form $Q\abd^{(\tau\tau)}(E)$ reads
\begin{eqnarray}
  Q\abd^{(\tau\tau)}(E) \!\! & = & \!\!
 \Frac{1}{4}\cdot\Frac{-\pi^2}{|s\abd|} \int_\Rt dk\ad \int_\Rt dk'\bd
\int_{S^2} d\hat{p}\ad \int_{S^2} d\hat{p}'\bd
\int_0^{E} d\rho {\rho}
\int_0^{E} d\rho' {\rho'}     \nonumber \\  [0.1cm]
& & \!\! \times \;\; f\big(k\ad,\sqrt{\rho}\,\hat{p}\ad\big)
 f'\big(k'\bd,\sqrt{\rho'}\,\hat{p}'\bd\big)
                         \nonumber \\ [-0.2cm]
			 & & \mbox{ } \\ [-0.2cm]
& & \!\!  \times \;\; \ds {{
\tau_{\alpha} \big( k_{\alpha},{{\tilde k}_{\alpha}^{(\beta)}} \big(
\sqrt{\rho}\,{\hat{p}}_{alpha}, \sqrt{\rho'}{\hat{p}}'_{\beta}\big),
 E - \rho + i \hs{0.03cm} 0\big)} \over
{\rho+\rho'- 2 c\abd\,\sqrt{\rho}
\,\sqrt{\rho'}\big(\hat{p}\ad,\hat{p}'\bd \big)-s_{\alpha \beta}^2 E-
i\hs{0.03cm} 0}}
\nonumber \\ [0.1cm]
& & \!\! \times \;\; \Frac{\tau\bd\big(\tk\au\bd\big(\sqrt{\rho}\,\hat{p}\ad,
\sqrt{\rho'}\,\hat{p}\bd'\big),k'\bd,E-\rho'+i\hs{0.03cm} 0 ) }
{\rho+\rho'- 2 c\abd\,\sqrt{\rho}
\,\sqrt{\rho'}\big(\hat{p}\ad,\hat{p}'\bd\big)-s\abd^2 E-i\hs{0.03cm} 0}\,.
            \nonumber
\end{eqnarray}
Making the change of variables
$\rho\rightarrow\nu$, $\rho'\rightarrow\nu'$,
$\rho=E\nu,$ $\rho'=E\nu',$
we get the integral
\begin{eqnarray*}
 Q\abd^{(\tau\tau)}(z) \!\! &= & \!\!
 \Frac{-\pi^2}{|s\abd|}\cdot
\Frac{z^3}{4} \int_\Rt dk\ad \int_\Rt dk'\bd
\int_{S^2} d\hat{p}\ad \int_{S^2} d\hat{p}'\bd
\int_0^{1} d\nu {\nu} \int_0^{1} d\nu' {\nu'} \nonumber \\ [0.1cm]
& & \!\! \times \;\; f\big(k\ad,\sqrt{z}\,\sqrt{\nu}\,\hat{p}\ad\big)
\, f'\big(k'\bd,\sqrt{z}\,\sqrt{\nu'}\hat{p}'\bd\big)  \nonumber\\ [0.1cm]
& & \!\! \times \;\;    \Frac{
\tau\ad(\,\ldots,z(1-\nu))\tau\bd(\,\ldots,z(1-\nu')) }
{\nu+\nu'- 2 c\abd\sqrt{\nu}\,
\sqrt{\nu'}\big(\hat{p}\ad,\hat{p}'\bd\big)-s\abd^2 -i0}\,,
\nonumber
\end{eqnarray*}
where the denominator of the expression under the integral
sign includes no dependence on the parameter $z$.  In view of
holomorphy in $z$ of the numerator of this expression, the
integral $Q\abd^{(\tau\tau)}(z)$ admits an immediate analytic
continuation on $\Pilh,$\,\, $l=+1.$

Summarizing the above, we can assert that the kernels of
the iteration $\cQ^{(1)}(z)$ admit analytic continuation $\big($in the
sense of distributions over $\Osix \big)$ on the domain $\Pilh$ of
the three\,--\,body unphysical sheet $\Pi_l,$   $l_0=+1$.  A similar
assertion holds as well for the three\,--\,body sheet $\Pi_l$
with  $l_0=-1$.  Also, we can state that the result of
continuation may be represented as
\begin{eqnarray}
\quad \quad
\label{tRt0}
\reduction{\cQ^{(1)}}{\Pi_l} \!\! & = & \!\! \bt^l\bRo^l\Y\bt^l
\nonumber \\ [-0.2cm]
& & \\ [-0.2cm]
&  = & \!\! \big(\bt - L_0 A_0\bt\bJot\bs_l^{-1}\bJo\bt\big)
\big(\bRo + L_0 A_0 \bJot \bJo\big)\Y \big(\bt - L_0 A_0\bt\bJot\bs_l^{-1}
\bJo\bt\big)\,.
\nonumber
\end{eqnarray}

When studying a continuation of the iteration $\cQ^{(1)}(z)$
into the three\,--\,body unphysical sheets $\Pi_l$,  $l_0=\pm 1$
in the general case where the pair subsystems may have eigenstates
one arrives again at the formula
$
%
      \reduction{\cQ^{(1)}}{\Pi_l}=\bt^l\bRo^l\Y\bt^l.
$
However, in contrast to (\ref{tRt0}) one must now use for $\bt^l(z)$
the total expressions (\ref{tltot}).

The proof is completed.%
\end{proof}

\begin{remark}\label{NDomainF}
{\rm Theorem~\ref{ThIter} means that one can
pose the continued Faddeev equations  (\ref{MFEl}) only in the domains
$\Pilh\subset\Pi_l$.
}
\end{remark}

\newsection{Representations for the analytic
            continuation of the matrix $M(z)$
             in unphysical sheets}
\label{STRepres}
In the present section we use the continued Faddeev
equations~(\ref{MFEl}) to obtain  representations for the
matrix $M^l(z)$ in the domains $\Pilh$ of the unphysical sheets
$\Pi_l\subset\Re$.  The representations will be given in terms
of the matrix $M(z)$ components themselves taken in the physical
sheet or, more precisely, in terms of the half\,--\,on\,--\,shell matrix
$M(z)$ as well as the inverse operators of the truncated
scattering matrices $S_l(z)$ and $\St_l(z)$.  As a matter of
fact, the construction of the representations for  $M^l(z)$
consists in explicitly ``solving'' the continued Faddeev
 equations~(\ref{MFEl}) in the same way as in~\cite{MotTMF},
\cite{MotYaF} where representations of the type~(\ref{3tnf}) had
been found for analytic continuation of the T\,--\,matrix in the
multichannel scattering problem with binary channels.  We
consider derivation of the representations for $M^l(z)$ as a
constructive proof of existence $\big($in the sense of
distributions over $\times_{\an=1}^3 \Osix\big)$ of analytic
continuation of the matrix  $M(z)$ into the unphysical sheets
$\Pi_l$ of the surface $\Re$.

So, let us consider the Faddeev equations~(\ref{MFEl}) in the
sheet $\Pi_l$ with  $l_0=0$ or $l_0=\pm 1$ and $l\bjd=0$ or
$l\bjd=1,$  $\beta = 1,2,3, \; j =1,2,\,\ldots, n_{\beta}$.  Using the
expressions~(\ref{tltot})
for $\bt^l(z)$ and~(\ref{R0ltot}) for $\bRo^l(z)$, we transfer
all the summands including $M^l(z)$ but not $\bJo$ and $\rJ_1$,
to the left side of Eqs.~(\ref{MFEl}).  Making then a
simple transformation based on the identity
$
\bs_l^{-1}(z)=\hbIo -\bs_l^{-1}(z)\bJo(z)\bt(z)\bJot(z)A_0(z)L_0
$
we rewrite~(\ref{MFEl}) in the form
\be
\label{MFEltr}
\mbox{\phantom{MMM}}
   (\bI + \bt \bRo\Y)M^l \EQ
 \bt\!\left[ \bI-\Aol\bJot\bs_l^{-1}\bJo\bt - \Aol\bJo\Xol \right]\,-\,
\Phi\rJt_1 \Ail \Big(\rJ_1\Phi^{*}+\Xil\Big)
\ee
where
$
            \Aol(z)=L_0 A_0(z),
$
$
 \Ail(z)=L_1 A_1(z).
$
Besides, we denote
\be
\label{DefX0}
\begin{array}{lcl}
 \Xol \!\! & = & \!\!|L_0| \bs_l^{-1}\bJo(\bI-\bt\bRo)\Y M^l\,,\\ [0.15cm]
   \Xil \!\! & = & \!\!-\,L_1\!
\!\left[ \rJ_1\Phis\bRo+\Aol\rJ_1\Phis\bJot\bJo
   \right]\!\Y M^l\,.
\end{array}
\ee

It should be noted that
\vs{0.2cm}\be
\label{JFR0}
 \rJ_1\Phis\bRo \EQ -\,\rJ_1\Psis\,.
\ee
Indeed,
\begin{eqnarray*}
 \big( \rJ\ajd \langle\;\cdot\; ,\phi\ajd\rangle  R_0   \big)
\big(\hp\ad,P'\big) \!\! & = & \!\!
 \int_\Rt dk''\ad \int_\Rt dp''_\alpha      
\Frac{\delta\big(\sqrt{z-\lambda\ajd}-|p''\ad|\big)}{|p''\ad|^2}
\\ [0.1cm]
& & \!\! \times \;\; \delta\big(\hp\ad,\hp''\ad\big)\,\overline{\phi\ajd}
\big(k''\ad\big)
  \ldots  \Frac{\delta\big(k''\ad-k'\ad\big)\delta
\big(p''\ad-p'\ad\big)}%
{ {k'}^2\ad+ {p'}^2\ad -z}        \\ [0.1cm]
& = & \!\!\Frac{ \overline{\phi\ajd}(k'\ad)
\delta\big(\sqrt{z-\lambda\ajd}-|p'\ad|\big)   }
{ |p'\ad|^2 \, ({k'\ad}^2+z-\lambda\ajd-z) }  \cdot
\delta\big(\hp\ad,\hp'\ad\big) \\ [0.1cm]
& = & \!\! \Frac{ \overline{\phi\ajd}(k'\ad)  }{ {k'\ad}^2-\lambda\ajd}
\,\rJ\ajd\big(z,\hp\ad,p'\ad\big)\,.
\end{eqnarray*}
Then it follows from Eq.~(\ref{psiphij}) that
$
\rJ\ajd \langle\;\cdot\; ,\phi\ajd\rangle  R_0 =
-\rJ\ajd \langle\;\cdot\; ,\psi\ajd\rangle,
$
and thereby the equality~(\ref{JFR0}) is really true.

Along with~(\ref{JFR0}) the equalities
\be
\label{JFJtot}
\Big( \rJ_1\Phis\bJot \Big)\! (z) \EQ 0\,, \quad
\Big( \bJo\Phi\rJt_1   \Big)\! (z) \EQ 0\,,
\ee
hold in accordance with~(\ref{JFJ}) for all
$z\in{\C}\setminus(-\infty,\lambda_{\rm max}]$.

Note that the condition  $z\not\in(-\infty,\lambda_{\rm max})$
necessary for Eq.~(\ref{JFJtot}) to be valid, does not apply to the
two\,--\,body unphysical sheets $\Pi_l$, $l_0=0,$ since in these sheets
$\Aol(z)=0$ and consequently, the terms including the products
$\bJot\bJo$ are absent in~(\ref{MFEltr}).  Meanwhile, the
points $z\in(-\infty,\lambda_{\rm max}]$ were excluded from the
three\,--\,body sheets $\Pi_l,$\,\, $l_0=\pm 1$, by definition.

Using Eq.~(\ref{JFR0}) and the first of Eqs.~(\ref{JFJtot})
one can rewrite $\Xil$ in the form
\be
\label{DefX1}
\Xil \EQ L_1 \rJ_1\Psis\Y M^l\,.
\ee

Notice further that the operator $\bI+\bt\bRo\Y$ admits an
explicit inversion in terms of  $M(z)$,
\vs{0.2cm}\be
\label{inv}
(\bI+\bt\bRo\Y)^{-1} \EQ \bI\,-\,M\Y\bRo\,,
\ee
for all $z\in\Pi_0$ which do not belong to the discrete spectrum
$\sigma_d(H)$ of the Hamiltonian $H$, and
\vs{0.2cm}\be
\label{invt}
(\bI-M\Y\bRo)\bt \EQ M\,.
\ee
The equality~(\ref{inv}) is a simple consequence of the
Faddeev equations~(\ref{MFE}) and the identity $\bRo\Y=\Y\bRo$.
The relation~(\ref{invt}) represents the alternative variant
(\ref{MFEAl}) of these equations.
Now, we can rewrite Eqs.~(\ref{MFEltr})
in equivalent form
\be
\label{Ml3in}
\begin{array}{lcl}
M^l \!\! & = & \!\!M\!\Big(\bI-\Aol\bJot\bs_l^{-1}\bJo\bt-\Aol\bJot\Xol
\Big)   \\ [0.15cm]
& & \!\!  - \;\;(\bI-M\Y\bRo)\Phi\rJt_1\Ail\!\Big(\rJ_1\Phis+\Xil\Big).
\end{array}
\ee
Eq.~(\ref{Ml3in}) means that the matrix $M^l(z)$
is expressed in terms of the quantities $\Xol(z)$ and $\Xil(z)$.
The main goal of this section consists in the
representation of these quantities in terms of the matrix
$M(z)$ considered in the physical sheet.

To obtain for $\Xol$ and $\Xil$ a closed system of equations
we use the definitions~(\ref{DefX0}) and~(\ref{DefX1}) and
apply the operators $\bs_l^{-1}\bJo(\bI-\bt\bRo)\Y$
and $\rJ_1\Psis$ to  both parts of Eq.~(\ref{Ml3in}).
At the moment we use also the identities
\be
\label{Help1}
[ \bI -\bt\bRo]\Y M \EQ M_0 \,-\,\bt\,, \quad
[ \bI -\bt\bRo]\Y [\bI-M\Y\bRo] \EQ [\bI-M_0\bRo]\Y
\ee
where $M_0=\Omt\Om M=(\bI+\Y)M.$
The relations~(\ref{Help1})
are another easily verified consequence
of the Faddeev equations~(\ref{MFE}).
Along with Eqs.~(\ref{Help1}) we use here also the second of
the equalities~(\ref{JFJtot}).
As a result we come to the desired system
of equations for  $\Xol$ and $\Xil$:
\begin{eqnarray}
\Xol \!\!&=&\!\!|L_0| \,\bs_l^{-1}\bJo\!\left[ (M_0-\bt)\Big(\bI-
\Aol\bJot\bs_l^{-1}\bJo\bt -
\Aol\bJot\Xol\Big)\right]  \nonumber \\ [-0.2cm]
& & \\ [-0.2cm]
& & \!\! - \, |L_0| \,\bs_l^{-1}\bJo M_0\Y\Psi\rJt_1\Ail
\Big(\rJ_1\Phis+\Xil\Big), \nonumber
\label{sys1}         \\ [0.2cm]
\Xil \!\! & =& \!\!L_1 \rJ_1\Psis\Y
M\Big(\bI-\Aol\bJot\bs_l^{-1}\bJot\bt-\Aol\bJot\Xol\Big)
\nonumber          \\ [-0.2cm]
& & \mbox{ } \\ [-0.2cm]
&& \!\! - \;\; L_1 \rJ_1\Psis\Y[\Phi+M\Y\Psi]\rJt_1\Ail\Big(\rJ_1\Phis+\Xil
\Big).
\label{sys2} \nonumber
\end{eqnarray}
It is convenient to rewrite this system in matrix form
             $$
\Btl\Xl \EQ \Dtl\,,\quad \Xl \EQ \Big(\Xol,\Xil\Big)^\dagger
             $$
with $\Btl=\big\{ \Btl_{ij}  \big\},\; i,j=0,1,$
the matrix consisting
of the operators appearing in the unknowns $\Xol$
and $\Xil$. By
$\Dtl,$  $\Dtl=\Big(\Dtl_0, \Dtl_1\Big)^\dagger, $
we understand a column constructed of the absolute terms
of Eqs.~(\ref{sys1}) and~(\ref{sys2}). Since
$\bs_l=\hbIo+\Aol\bJo \bt\bJot$ we find
\begin{eqnarray*}
\Btl_{00} \!\!& = &\!\! \hbIo \,+\,\Aol \bs_l^{-1}\bJo(M_0-\bt)\bJot
\\[0.1cm]
&=& \!\! \bs_l^{-1}\! \Big(\hbIo + \Aol\bJo\bt\bJot +
\Aol\bJo M_0\bJot -\Aol\bJo\bt\bJot\Big)                     \\ [0.1cm]
& = & \!\! \bs_l^{-1} \!\Big(\hbIo +\Aol\bJo M_0\bJot\Big).
\end{eqnarray*}
At the same time
\begin{eqnarray*}
& & \Btl_{01}  \EQ  |L_0| \bs_l^{-1}\bJo M_0\Y\Psi\rJt_1\Ail\,, \\ [0.1cm]
& & \Btl_{10}  \EQ  L_1 \rJ_1\Psis\Y M\bJot \Aol\,, \\ [0.1cm]
& & \Btl_{11} \EQ \hIi+L_1 \rJ_1\Psis U\Psi\rJt_1\Ail
\end{eqnarray*}
because
$
\Y(\Phi+M\Y\Psi)=\Y\bv\Psi+\Y M\Y\Psi=(\Y \bv +\Y M\Y)\Psi =U\Psi
$
(see Sec.~\ref{SMSphys}).

The absolute terms are
\begin{eqnarray*}
\Dtl_0 \!\!& = & \!\!|L_0|\,\bs_l^{-1}\Big[\bJo(M_0-\bt)\Big(\bI-\Aol\bJot
\bs_l^{-1}\bJo\bt\Big)
\,-\,|L_0|\bJo M_0\Y\Psi\rJt_1\Ail\rJ_1\Phis\Big],  \\ [0.1cm]
\Dtl_1 \!\!& = &\!\! L_1 \rJ_1\Psis\Y M\Big(\bI-\Aol\bJot\bs_l^{-1}\bJo\bt\Big)
 \,- \, L_1 \rJ_1\Psis U\Psi\rJ_1\Ail\rJ_1\Phis.
\end{eqnarray*}

The operator $\bs_l(z)$,  $l_0=\pm 1$ is invertible for
all $z\in\cP_b$.  If $z\not\in Z_{\rm res}$, then $\bs_l^{-1}(z)$
is a bounded operator in $\what{\cG}_0$.  Therefore, applying
the operator $\bs_l$ to both parts of the first equation
$\Btl_{00}\Xol+\Btl_{01}\Xil=\Dtl_0$
of the system $\Btl\Xl=\Dtl$, and not changing the second equation,
we come to the equivalent system
\be
\label{BXD}
\Bl\Xl \EQ \Dl
\ee
where
\be
\label{Bl}
\Bl \EQ \left(\!\!\begin{array}{cc}
\hbIo+|L_0|\bJo M_0\bJot\Aol  & |L_0|\bJo M_0\Y \Psi\rJt_1\Ail   \\ [0.2cm]
L_1 \rJ_1\Psis\Y M\bJot\Aol & \hIi +L_1 \rJ_1\Psis U\Psi\rJt_1\Ail
\end{array}\!\!\right),
\ee
$\Bl(z):$
$\what{\cG}_0\oplus\what{\cH}_1\rightarrow \what{\cG}_0\oplus\what{\cH}_1$.
The absolute term $\Dl$ has components
$\Dl_0=\bs_l \Dtl_0$ and $\Dl_1=\Dtl_1$.
\begin{lemma}\label{LBlInv}
The inverse operator $\left(\Bl(z)\right)^{-1}$
exists for all $z\in\Pi_l^{(\rm hol)}$ where
the inverse operator $S_l^{-1}(z)$
of the truncated three\,--\,body scattering matrix
$S_l(z)$ given by the first of the equalities in~{\rm (\ref{Slcut})}
exists with
$L=\diag\hs{0.03cm}\{L_0,L_1\}$, $\tL=\diag\hs{0.03cm}\{|L_0|,L_1\},$
and where the inverse operators $[S_l(z)]^{-1}_{00}$
and $[S_l(z)]^{-1}_{11}$ of
$[S_l(z)]_{00}=\hat{I}_0+\rJo T\rJot A_0 L_0$
and $[S_l(z)]_{11}=\hat{I}_1+L_1 \rJ_1\Psis U\Psi\rJt_1 A_1 L_1,$
respectively, exist. The components
$\big[ \big(\Bl(z) \big)^{-1}\big]_{ij},$ $ i,j=0,1,$
of the operator $\left(\Bl(z)\right)^{-1}$
admit the representations
\begin{eqnarray}
\label{Y00}
\Big[ \big(\Bl(z) \big)^{-1}\Big]_{00}\!\! & = & \!\!
\hbIo\,-\,\Omt\big[S_l^{-1}\big]_{00}  \\ [0.1cm]
\nonumber
 & & \!\! \times \;\;\left\{|L_0|\rJo T_0 -
 [S_l]_{01}[S_l]_{11}^{-1}L_1\rJ_1\Psis\Y M\right\}
\bJot\Aol, \\ [0.1cm]
\label{Y01}
\Big[ \big(\Bl(z) \big)^{-1}\Big]_{01}\!\! & = & \!\!\Omt\big[S_l^{-1}
\big]_{01}\,,\\ [0.1cm]
\label{Y10}
\Big[ \big(\Bl(z) \big)^{-1}\Big]_{10} \!\!& = &\!\! -\,\big[S_l^{-1}\big]_{11}
L_1 \rJ_1\Psis\Y M\bJot\Aol  \\ [0.1cm]
\nonumber  &  & \!\! \times
\left\{\hbIo -\Omt[S_l]^{-1}_{00}|L_0|\rJo T_0\bJot\Aol \right\},  \\ [0.1cm]
\label{Y11}
\Big[ \big(\Bl(z) \big)^{-1}\Big]_{11} \!\!& = & \!\!\big[S_l^{-1}\big]_{00}
\end{eqnarray}
where $T_0\equiv \Om M$.
\end{lemma}

\vs{0.3cm}Note that since $|L_0|$ and $\Aol$ are numbers which
become zero for $l_0=0$ simultaneously, the factors $|L_0|$ in
(\ref{Y00}) and (\ref{Y10}) may be omitted.

\begin{proof}
Let us find at the beginning, the components
$
\Big[ \big(\Bl(z) \big)^{-1}\Big]_{00}
$
and \linebreak
$\Big[ \big(\Bl(z) \big)^{-1}\Big]_{10},$
which will be denoted temporarily (for the sake of brevity)
by $Y_{00}$ and $Y_{10}$.  Using Eq.~(\ref{Bl}) we
write the system of equations for these components as
\begin{eqnarray}
\label{EQY00}
\big[\Bl\big]_{00}\, Y_{00}\,+\,\big[\Bl\big]_{01}\, Y_{10} \!\!&=&\!\!
\hbIo\\ [0.1cm]
\label{EQY10}
\big[\Bl\big]_{10}\, Y_{00}\,+\,\big[\Bl\big]_{11}\, Y_{10}\!\!&=&\!\! 0\,.
\end{eqnarray}
Eliminating the unknown $Y_{10}$  from the first equation~(\ref{EQY00})
with the help of~(\ref{EQY10}) we come to the following
equation including the element $Y_{00}$ only,
\be
\label{EY00}
\left\{ \hbIo+\Omt\!\!\left[|L_0|\rJo T_0\bJot\Aol -
[S_l]_{01}[S_l]^{-1}_{11} L_1 \rJ_1\Psis\Y M\bJot\Aol\right]   \right\}
Y_{00} \EQ \hbIo\,.
\ee

The operator-matrix on the left\,--\,hand side of Eq.~(\ref{EY00})
complementary to $\hbIo$ has three identical rows.
Thus one can apply
to Eq.~(\ref{EY00}) the inversion formula
\be
\label{C123}
\left[ \hbIo+\Omt(C_1,\, C_2,\, C_3)\right]^{-1}\EQ \hbIo -
\Omt\!\left[ \hIo+C_1+C_2+C_3\right]^{-1}\!(C_1,\, C_2,\, C_3)\,,
\ee
which is true for a wide class of operators $C_1,\, C_2$ and $C_3$.
The single essential requirement on
$C_1,\, C_2$ and $C_3$ evidently, is the existence of
$\big(\hat{I}_0+C_1+C_2+C_3\big)^{-1}$.

In the case concerned
$$
C\bd(z) \; \equiv \; \left\{         |L_0|  \rJo T_{0\bn}\rJot
 -[S_l]_{01}[S_l]_{11}^{-1} L_1 \rJ_1\Psis\Y[M]\bd\rJot  \right\}\Aol
$$
where $[M]\bd$ is the
$\beta$\,--\,th column of the matrix  $M$,
$[M]\bd=\big(M_{1\bn},M_{2\bn},M_{3\bn}\big)^\dagger$.
Thus
\begin{eqnarray*}
\hIo\,+\,C_1\,+\,C_2\,+\,C_3 \!\! &= & \!\! \hIo\,+\,\rJo T\rJot\Aol \,-\,
[S_l]_{01}[S_l]_{11}^{-1}
\rJ_1\Psis U_0^{\dagger}\rJot\Aol \\ [0.1cm]
& \equiv & \!\!
[S_l]_{00}\,-\,[S_l]_{01}[S_l]_{11}^{-1}[S_l]_{10}.
\end{eqnarray*}

Note that the components $\left[S_l^{-1}\right]_{ij}$,
$i,j=0,1,$ of $S_l^{-1}$  have the representations
\begin{eqnarray}
\label{Si00}
\Slis_{00}\!\! &=&\!\!\left( \Sls_{00}-
\Sls_{01}\Sls_{11}^{-1}\Sls_{10}\right)^{-1}     \\ [0.2cm]
\label{Si11}
\Slis_{11} \!\!& = &\!\! \left( \Sls_{11}-
\Sls_{10}\Sls_{00}^{-1}\Sls_{01}\right)^{-1}     \\ [0.2cm]
\label{Si10}
\Slis_{10}\!\! & = & \!\!-\, \Sls_{11}^{-1}\Sls_{10}\Slis_{00}   \\ [0.2cm]
\label{Si01}
\Slis_{01}\!\! & = & \!\!- \,\Sls_{00}^{-1}\Sls_{01}\Slis_{11}
\end{eqnarray}
in terms of the components $[S_l]_{ij}$.
It follows from~(\ref{Si00}) that
$\hIo+C_1+C_2+C_3=\left(\Slis_{00}\right)^{-1}$. Therefore,
in the conditions of the Lemma, the operator
$\hat{I}_0+C_1+C_2+C_3$ is invertible. Now,
an application of Eq.~(\ref{C123}) in~(\ref{EY00})
leads us immediately to the representation~(\ref{Y00}) for
$\Big[ \left( \Bl\right)^{-1} \Big]_{00}$.

When calculating $Y_{10}=\Big[ \left( \Bl\right)^{-1} \Big]_{10}$
we eliminate from the second equation~(\ref{EQY10})
the quantity $Y_{00}$ using Eq.~(\ref{EQY00}).
In all of this, we need to calculate the inverse operator of
$\hbIo+\bJo M_0\bJot\Aol$.
Here we apply again the relation~(\ref{C123}) and obtain
\begin{eqnarray}
\left( \hbIo+|L_0|\bJo M_0\bJot\Aol \right)^{-1} \!\!& = &\!\!
\left( \hbIo+\Omt|L_0|\rJo T_0\bJot\Aol \right)^{-1} \nonumber \\ [-0.2cm]
& & \mbox{ } \\ [-0.2cm]
 & = & \!\!\hbIo-\Omt\Sls_{00}^{-1}|L_0|\rJo T_0\bJot\Aol\,.
\nonumber \label{JM0J}
\end{eqnarray}

With the help of~(\ref{Slcut}) we can write the resulting equation for
$Y_{10}$ as
\begin{eqnarray}
\nonumber
& & \!\!\left\{ \Sls_{11}-
\Sls_{10}\Sls_{00}^{-1}\Sls_{01}\right\}Y_{10} \nonumber \\ [-0.2cm]
\label{EY10}
& & \mbox{ } \\ [-0.2cm]
&  = & \!\! -\, \rJ_1\Psis\Y M\bJot\Aol
\!\left[ \hbIo+\bJo M_0\bJot\Aol \right]^{-1}.
\nonumber
\end{eqnarray}
According to Eq.~(\ref{Si11}) the expression in
braces on the left\,--\,hand side of Eq.~(\ref{EY10}) coincides with
$\Slis_{11}^{-1}$. Then, from~(\ref{EY10}) we get
immediately~(\ref{Y01}).

The system of the equations
\begin{eqnarray}
\label{EQY01}
\big[B^{(l)}\big]_{00}\, Y_{01} \,+\, \big[B^{(l)}\big]_{01}\, Y_{11} \!\! & =
& \!\! 0
\\ [0.2cm]
\label{EQY11}
\big[B^{(l)}\big]_{10}\, Y_{01} \,+\, \big[B^{(l)}\big]_{11}\, Y_{11}
\!\! & = & \!\!  \hIi
\end{eqnarray}
for the components  $Y_{01}=\big[(B^{(l)})^{-1}\big]_{01}$
and $Y_{11}=\big[(B^{(l)})^{-1}\big]_{11}$
is solved analogously. The search for $Y_{11}$
is a simple problem, since application of
the inversion formula~(\ref{JM0J})
to Eq.~(\ref{EQY01})  immediately gives
$
Y_{01}=\Omt\Sls_{00}^{-1} \Sls_{01} Y_{11}.
$
Substituting this $Y_{01}$ in~(\ref{EQY11}) we find
$$
\left\{ \Sls_{11}-\Sls_{10}\Sls_{00}^{-1}\Sls_{01}\right\}Y_{11} \EQ \hIi\,.
$$
As in Eq.~(\ref{EY10}) the operator on the left\,--\,hand side
is just $\Slis_{11}^{-1}$. Inverting it, we come to
Eq.~(\ref{Y11}).

When calculating the unknown $Y_{01}$, we begin by
expressing  the unknown $Y_{11}$ in terms of it.
Using Eq.~(\ref{EQY11}) we find
\be
\label{XY11}
Y_{11} \EQ
\Sls_{11}^{-1}\!\left(\hIi-L_1\rJ_1\Psi\Y M\bJo\Aol Y_{01}\right).
\ee
Substituting~(\ref{XY11}) into Eq.~(\ref{EQY01}) we obtain an
equation with an operator in the position of $Y_{01}$, which may be
inverted with the help of Eq.~(\ref{C123}). Then we use the
chain of equalities
\begin{eqnarray*}
|L_0|\,\bJo M_0\Y\Psi\rJt_1\Ail \!\! & = & \!\! |L_0|\,
\bJo\Omt\Om M\Y\Psi\rJt_1\Ail
\\ [0.1cm]
& = & \!\! \Omt|L_0|\rJo\Om M\Y\Psi\rJt_1\Ail \\ [0.1cm]
& = & \!\! \Omt\!\Sls_{01},
\end{eqnarray*}
simplifying the absolute term as well as the summand
on the left\,--\,hand side of the equation for $Y_{01}$
appearing there due to~(\ref{XY11})
from the element $\left[\Bl\right]_{01}$.
Completing the transformations we find
$$
Y_{01} \EQ -\,\Omt\!\left\{
\Sls_{00}-\Sls_{01}\Sls_{11}^{-1}\Sls_{10}\right\}^{-1} \!
\Sls_{01}\!\Sls_{11}^{-1}.
$$
In view of~(\ref{Si01}), the expression appearing
after $\Omega^{\dagger}$ on the right\,--\,hand side
of the last equation coincides
exactly with that for $\Slis_{01}$. Therefore,
we obtain finally Eq.~(\ref{Y01}). Thus, all the components
of the inverse operator $\left(\Bl\right)^{-1}$
have already been calculated.

It follows from the representations~(\ref{Y00}) -- (\ref{Y11})
that  $\left(\Bl(z)\right)^{-1}$ exists for those
$z\in\Pi_l^{(\rm hol)}$ where the inverse operators of
$S_l(z)$,\, $\left[S_l(z)\right]_{00}$ and
$\left[S_l(z)\right]_{11}$ exist.  This completes the proof of
the Lemma.%
\end{proof}

Let us return to Eq.~(\ref{BXD}) and invert
the operator $\Bl(z)$ using the relations~(\ref{Y00}) -- (\ref{Y11}).
In this way we find the unknowns $\Xol$ and $\Xil$ which express
$M^l(z)$ [see Eq.~(\ref{Ml3in})].

When carrying out a concrete calculation of $\Xol=\Blis_{00}
\!\Dl_0+\Blis_{01}\!\Dl_1$ we use the relation
$|L_0| \Big[\big(B^{(l)}\big)\Big]_{00}\bJo M_0 =\Omt |L_0| \!\Slis_{00}\rJo
T_0$ that can be
checked with the help of~(\ref{Slcut}) and~(\ref{T3body}).  Along
with the identity
\be
\label{sjtl3}
\bJo\bt\!\Big(\hbIo-\Aol\bJot\bs_l^{-1}\bJo\bt\Big) \EQ \bs_l^{-1}\bJo\bt\,,
\ee
this relation simplifies essentially the transform of the
product $\Blis_{00}\!\Dl_0 $.  In addition, when calculating $
\Xol$ we use the equalities~(\ref{JFJtot}).  As a result we find
\be
\begin{array}{lcl}
\label{X0}
  \mbox{\phantom{MMM}}
\Xol \!\! & = & \!\! \Omt\!\left\{ |L_0|\Slis_{00}\rJo T_0 +\Slis_{01}
L_1 \left(\rJ_1\Psis\Y M+\rJ_1\Phis \right) \right\}
\\ [0.15cm]
& & \!\! - \;\; |L_0|\,\bs_l^{-1}\bJo\bt\,.
\end{array}
\ee

Now, to find
$
%
 \Xil=\Blis_{10}\Dl_0+\Blis_{11}\Dl_1
$
we observe additionally that the equality
$
\left\{\hbIo-\Omt\Slis^{-1}_{00}\rJo T_0\bJot\Aol \right\}
\bJo M_0 =\Omt\Slis_{00}^{-1}\rJo T_0
$
which simplifies the product $\Blis_{10} \Dl_0$ is valid.
The final expression for $\Xil$ reads as follows
\begin{small}
\be
\label{X1}
  \mbox{\phantom{MMM}}
\Xil \EQ L_1 \!\left\{ \Slis_{10} |L_0|\,\rJo\, T_0+
\Slis_{11} L_1\rJ_1\Psis\Y M-
\left( \hIi-\Slis_{11}\right) \!L_1\rJ_1\Psis\right\}.
\ee
\end{small}

To obtain now a representation for $M^l(z)$, one has only to
substitute in Eq.~(\ref{Ml3in}) the expressions~(\ref{X0})
for $\Xol$ and~(\ref{X1}) for $\Xil$.  Carrying out a series of
simple but rather cumbersome transformations of
Eq.~(\ref{Ml3in}) we arrive as a result at a statement analogous
to Theorem~\ref{Tht2body} concerning analytical continuation of
the two\,--\,body T\,--\,matrix.  The statement is the following.
\begin{theorem}\label{ThMlRepr}
The matrix $M(z)$ admits, in the sense of distributions over
$\Osix$, an analytic continuation in $z$ on the domains $\Pilh$
of the unphysical sheets $\Pi_l$ of the surface $\Re$.
The continuation is described by
\be \hs{0.5cm}
\label{Ml3fin}
M^l \EQ M-\left(M\Omt\rJot,   \Phi\rJt_1+M\Y\Psi\rJt_1\right)\!
LA\, S_l^{-1}  \tL  \!\left( \!\!
\begin{array}{c}
\rJo\Om M \\ [0.2cm]
\rJ_1\Psis\Y M+ \rJ_1\Phis
\end{array}       \!\!  \right)
\ee
where $S_l(z)$ stands for the truncated scattering matrix
{\rm(\ref{Slcut})},
\begin{eqnarray*}
& &  L \EQ \diag\hs{0.03cm}\{l_0,l_{1,1},\,\ldots,
l_{1,n_1}, l_{2,1},\,\ldots,
l_{2,n_2},l_{3,1},\,\ldots,l_{3,n_3}\}\,\\ [0.1cm]
& & \tL \EQ \diag \hs{0.03cm}\{|l_0|,l_{1,1},\,\ldots,
l_{1,n_1}, l_{2,1},\,\ldots,
l_{2,n_2},l_{3,1},\,\ldots\,,l_{3,n_3}\}\,.
\end{eqnarray*}
The kernels of all the operators on the right\,--\,hand side of
Eq.~{\rm(\ref{Ml3fin})}
are taken in the physical sheet.
\end{theorem}

\vs{0.3cm}Note that  $LA\,S_l^{-1}(z)\tL=\tL\,\big[S^\dagger_l(z)\big]^{-1}
AL$.
This means that the relations~(\ref{Ml3fin}) may also be rewritten
in terms of the scattering matrices
$S^\dagger_l(z)$.

\newsection{Analytic continuation of the scattering matrices}
\label{SSmxl}
Let $l=\left\{l_0,\, l_{1,1},\,\ldots,l_{1,n_1},l_{2,1},\,\ldots,l_{2,n_2},
 l_{3,1},\,\ldots,l_{3,n_3}\right\}$
with certain $l_0$, $l_0=0$ or $l_0=\pm 1$, and
 $l\ajd$, $l\ajd=0$  or $l\ajd=+1$, $\alpha= 1,2,3, \;j=1,2,\,\ldots,
n_{\alpha}$.
The truncated scattering matrices  $S_l(z)\!:\!
\what{\cH}_0\oplus\what{\cH}_1 \!\rightarrow \! \what{\cH}_0\oplus\what{\cH}_1$
and $\St_l(z)\!:\!
\what{\cH}_0\oplus\what{\cH}_1 \!\rightarrow \! \what{\cH}_0\oplus
\what{\cH}_1$,
given by formulas (\ref{Slcut}),
are operator\,--\,valued functions of
the variable $z$ which are holomorphic
in the domain $\Pilh$ of the physical sheet $\Pi_0$.
For $l_0=1$ and $l\ajd=1$, $\alpha = 1,2,3, \;j=1,2,\,\ldots,n_{\alpha}$,
these matrices coincide with the respective
total three\,--\,body scattering matrices:
$S_l(z)=S(z),$  $\St_l(z)=\St(z).$

We describe now the analytic continuation of
the truncated scattering matrices$^{4)}$%
\footnote{%
$^{4)}$Note that the analytic properties of the truncated
scattering matrices or, more exactly, the ($2\rightarrow 2$)
scattering amplitudes in the $N$\,--\,body system with $N\geq 3$ were
investigated in the paper~\cite{Derezinski} in the case of the type
(\ref{vpotb}) pair interactions. A proof is given
in~\cite{Derezinski} for existence of analytic continuation of these
amplitudes through the cut in vicinities of the branches of the
continuous spectrum below the first threshold of the system to
breakup into three clusters.
}
$S_{l'}(z)$ and $\St_{l'}(z)$
with a certain multi\,--\,index $l'$
in the unphysical sheets $\Pi_{l}\in\Re$.

We shall use here the representations (\ref{Ml3fin}) for
$\reduction{M(z)}{\Pi_{l}}$. As mentioned above, our goal is to
find explicit representations for $\reduction{S_l(z)}{\Pi_{l'}}$
and  $\reduction{\St_l(z)}{\Pi_{l'}}$ again in terms of the
physical sheet.

First, we notice that the function $A_0(z)$ is univalent. It looks like
$A_0(z)=-\pi i z^2$ on all the sheets $\Pi_l$.  At the same time
after continuing from $\Pi_0$ on $\Pi_l$ the function
$A\bjd(z)=-\pi i\sqrt{z-\lambda\bjd}$ keeps its form only if
$l\bjd=0$.  If $l\bjd=1$, this function turns into
$A'\bjd(z)=-A\bjd(z)$.  Analogous inversion takes (or does not
take) place with arguments $\hP$, $\hP'$, $\hp\ad$ and $\hp'\bd$
of the kernels of the operators $\rJo \Om M\Omt\rJot$, $\rJo \Om
M\Y\Psi\rJt_1$, $\rJ_1\Psis\Y M\Omt\rJot$ and
$\rJ_1\Psis(\Y\bv+\Y M\Y)\Psi\rJt_1$, too.  Recall that on the
physical sheet $\Pi_0$, the action of $\rJo(z)$
$\big(\rJot(z)\big)$ transforms $P\in\Rs$ in $\sqrt{z}\hP$
$\big(P'\in\Rs$ in $\sqrt{z}\hP'\big)$.  At the same time, $p\ad\in\Rt$
 $\big(p'\bd\in\Rt\big)$ turns under $\rJ_{\an,i}(z)$
$\big(\rJt_{\bn,j}(z)\big)$ into $\sqrt{z-\lambda_{\an,i}}\, \hp\ad$
 $\big(\sqrt{z-\lambda_{\bn,j}}\, \hp'\bd\big)$.  That is why  we
introduce the operators $\cE(l)=\diag\hs{0.03cm}\{\cE_0,\,\cE_1\}$ where
$\cE_0$ is the identity operator in $\what{\cH}_0$ if $l_0=0$,
and $\cE_0$ is the inversion $\big(\cE_0
f\big)\big(\hP\hs{0.02cm}\big)= f\big(-\hP\hs{0.02cm}\big)$ if
$l_0=\pm 1$. Analogously, $\cE_1(l)=\diag\hs{0.03cm}\{
\cE_{1,1},\,\ldots,\cE_{1,n_1}; \cE_{2,1},\,\ldots,\cE_{2,n_2};
\cE_{3,1},\,\ldots,\cE_{3,n_3}\}$ where $\cE\bjd$ is the
identity operator in $\what{\cH}^{(\bn,j)}$ if $l\bjd=0$ and
$\cE\bjd$ is the inversion $(\cE\bjd f)(\hp\bd)=f(-\hp\bd)$, if
$l\bjd=1$.  By  $\re_1(l)$ we denote the diagonal matrix
\linebreak $\re_1(l)=\diag\hs{0.03cm}\{
\re_{1,1},\,\ldots,\re_{1,n_1};\re_{2,1},\,\ldots,\re_{2,n_2};
\re_{3,1},\,\ldots,\re_{3,n_3}\}$ with the elements  $\re\bjd=1$
if $l\bjd=0$ and $\re\bjd=-1$ if $l\bjd=1$.  Let
$\re(l)=\diag\hs{0.03cm}\{ \re_0, \re_1 \}$ where $\re_0=+1$.
\begin{theorem}\label{ThS3lRepr}
If there exists a path on the surface $\Re$ such that while
moving along it from the domain $\Pi^{\rm (hol)}_{l'}$ on
$\Pi_0$ to the domain $\Pi^{\rm (hol)}_{l'}\cap\,\Pi^{\rm
(hol)}_{l'l}$ on $\Pi_{l}$ the parameter $z$ stays in
intermediate sheets $\Pi_{l''}$ always contained in the domains
$\Pi^{\rm (hol)}_{l'}\cap\Pi^{\rm (hol)}_{l''l'}$, then the
truncated scattering matrices $S_{l'}(z)$ and $\St_{l'}(z)$
admit analytic continuation in $z$ on the domain
$\Pi^{\rm(hol)}_{l'}\cap\Pi^{\rm(hol)}_{l'l}$ of the sheet
$\Pi_{l}$.  The continuation is described by
\begin{eqnarray}
\reduction{S_{l'}(z)}{\Pi_{l}} \!\!& = &\!\! \cE(l)\!\!\left[ \hbI+
\tL' \what{\cT} L' A \re(l) -
\tL' \what{\cT} L A S^{-1}_{l} \tL \hat{\cT} L' A\re(l)
\right]\! \cE(l)\,,
\label{Slfin}                               \\
\reduction{\St_{l'}(z)}{\Pi_{l}} \!\!& = & \!\!\cE(l)\!\!\left[ \hbI+
 \re(l) A  L' \what{\cT} \tL'  -
\re(l) A\, L \what{\cT} \tL\, [\St_{l}]^{-1} A  L \what{\cT} \tL'
\right]\!\! \cE(l)\,,
\label{Stlfin}
\end{eqnarray}
where
\begin{eqnarray*}
& & L'\EQ \left\{l'_0, l'_{1,1},\,\ldots,l'_{1,n_1},l'_{2,1},\,\ldots,
l'_{2,n_2}, l'_{3,1},\,\ldots,l'_{3,n_3}\right\}\,, \\ [0.1cm]
 & & \tL'\EQ \left\{|l'_0|, l'_{1,1},\,\ldots,l'_{1,n_1},
 l'_{2,1},\,\ldots,l'_{2,n_2},
  l'_{3,1},\,\ldots,l'_{3,n_3}\right\}\,.
  \end{eqnarray*}
\end{theorem}
\begin{proof}
We give the proof for the case of $S_{l'}(z)$.
Using the definition~(\ref{T3body}) of the operator $\cT(z)$
we rewrite $S_{l'}(z)$ in the form
$$
S_{l'}(z)\EQ\hbI+\tL'\!\left[        \left(\!\!\begin{array}{c}
\rJo \Om     \\ [0.2cm] \rJ_1\Psis\Y   \end{array}\!\!\right)\!
 M \!\left(   \Omt\rJot\, ,    \Y\Psi\rJt_1  \right)
 \,+\, \left(\!\!\begin{array}{cc}
 0      &      0     \\ [0.2cm]
 0      &   \rJ_1\Psis\Y\bv\Psi\rJt_1
\end{array}\!\!\right)
\right] \!L' A\,.
$$
Note that when continuing into the sheet $\Pi_{l''}$, the operators
$\rJo(z),$ $\rJot(z),$ $\rJ_1(z)$ and  $\rJt_1(z)$
turn into  $\cE_0(l'')\rJo(z),$
$\rJot(z)\cE_0(l''),$ $\cE_1(l'')\rJ_1(z)$
and  $\rJt_1(z)\cE_1(l'')$, respectively. At the same time
the matrix\,--\,function $A(z)$ turns into $A(z)\re(l'')$.
Then, using Theorem~\ref{ThMlRepr} in the domains
$\Pi^{\rm (hol)}_{l'}\cap\Pi^{\rm (hol)}_{l'l''}$
of the intermediate sheets $\Pi_{l''}$ we have
\begin{eqnarray}
\nonumber
\reduction{S_{l'}(z)}{\Pi_{l''}}\! &=&\!\hbI\,+\,
\cE(l'')\tL'\!\left[ \left(\!\!\begin{array}{c}
\rJo \Om     \\ [0.1cm] \rJ_1\Psis\Y   \end{array} \!\!\right)\!
 M^{l''}\!\! \left( \Omt\rJot\, ,\,\,\,\,     \Y\Psi\rJt_1  \right) \right.
 \nonumber \\ [-0.2cm]
& & \mbox{ } \\ [-0.2cm]
\label{S3lin}
 & & \left. +\; \left(\!\!\begin{array}{cc}
 0      &      0     \\ [0.1cm]
 0      &   \rJ_1\Psis\Y\bv\Psi\rJt_1
\end{array}\!\!\right)\right]\!L' \cE(l'') A\re(l'')\,. \nonumber
\end{eqnarray}
Substitution of $M^{l''}(z)$  from~(\ref{Ml3fin})
shows that
\begin{eqnarray}
\nonumber
\reduction{S_{l'}(z)}{\Pi_{l''}} \!\! & = & \!\! \hbI\,+\,
\cE(l'')\tL'\what{\cT}L'\cE(l'') A\re(l'')     \\ [0.1cm]
\label{S3l}
& & \!\! -\;\;\cE(l'')\tL' \!\left(\!\! \begin{array}{c}
 \rJo\Om    \\  [0.1cm]    \rJ_1\Psis\Y
\end{array} \!\!\right)
\!\!
\left( M\Omt\rJot\, , [\bv+M\Y]\Psi\rJt_1  \right)
\!\!L''A\, S_{l''}^{-1} \\ [0.1cm]
\nonumber &&  \!\! \times \;\; \tL''\!\left(\!\! \begin{array}{c}
 \rJo\Om M   \\ [0.1cm]     \rJ_1\Psis[\bv+\Y M]
\end{array} \right)\! \!
\left( \Omt\rJot\, ,\, \Y\Psi\rJt_1  \right)\!L' \cE(l'') A \re(l'')
\end{eqnarray}
where the summand immediately following $\hbI$ is generated by
the term $M(z)$ of the right\,--\,hand side of~(\ref{Ml3fin}).
The last summand of~(\ref{S3l}) is comes from the second summand
of~(\ref{Ml3fin}).

In view of~(\ref{JFJtot}) we have
$\rJ_1\Psis\bv\hs{0.02cm}\Omt\rJot=\rJ_1\Phis\rJot\Omt=0$.
Analogously, $\rJo\Om\hs{0.02cm}\bv\Psi\rJt_1$ is equal to zero,
too. Thus, taking into account (\ref{T3body}) we find
\begin{eqnarray}
\label{S3ll}
\nonumber \reduction{S_{l'}(z)}{\Pi_{l''}}\!\! & = & \!\!\hbI\,+\,
\cE(l'') \tL'\what{\cT}L' \cE(l'') A\re(l'') \\ [-0.2cm]
& &  \mbox{ } \\ [-0.2cm]
& & \!\! - \;\;
\cE(l'') \tL'\what{\cT}L'' A\, S_{l''}^{-1}
\tL''\what{\cT}L' \cE(l'') A\re(l'')\,. \nonumber
\end{eqnarray}
By the assumption, the parameter $z$ moves along  a path such
that in the sheet $\Pi_{l''}$ it is situated in the domain
$\Pi^{\rm (hol)}_{l'}\cap\Pi^{\rm (hol)}_{l'l''}. $
In this domain, the operators
$\big(\tL'\hat{\cT}L'\big)(z)$, $\big(\tL'\hat{\cT}L''\big)(z)$ and
$\big(\tL''\hat{\cT}L'\big)(z)$ are defined and
depend on $z$ analytically.
Consequently, the same may be said also
about the function $\reduction{S_{l'}(z)}{\Pi_{l''}}$.
In equal degree, this statement is related to the sheet
$\Pi_{l}$. Replacing the values of the multi\,--\,index $l''$ in the
representations  (\ref{S3lin}) -- (\ref{S3ll}) with $l$,
we come to the assertion of the theorem
for $\reduction{S_{l'}(z)}{\Pi_{l}}$.
The validity of the representations
(\ref{Stlfin}) for  $\reduction{\St_{l'}(z)}{\Pi_{l}}$
is established in the same  way.%
\end{proof}

\begin{remark}\label{NSlPil}
{\rm If $l_0=0$ then the representation (\ref{Slfin})
for the analytic continuation of $S_l(z)$
into the sheet $\Pi_l$ (its ``own'') acquires the simple form
[cf.~(\ref{S2P1})]
$$
\reduction{S_l(z)}{\Pi_{l}} \EQ \cE(l)\!\!\left[ \hbI +\re(l) -
S_l^{-1}(z)\re(l)  \right]\!\cE(l) \EQ \cE(l) S_l^{-1}(z) \cE(l)\,.
$$
In the same way
$\reduction{\St_l(z)}{\Pi_{l}}=\cE(l)\big[\St_l(z)\big]^{-1} \cE(l).$
}
\end{remark}

\newsection{Representations for the analytic continuation of the resolvent
             in the unphysical sheets}
\label{SResolvl}
The resolvent $R(z)$ of the Hamiltonian $H$ for the three\,--\,body system
is expressed by $M(z)$ according to Eq.~(\ref{RMR}).
As we have established, the
kernels of all the operators included in the
right\,--\,hand side of~(\ref{RMR}) admit, in the sense of
distributions over  $\Osix$, analytic continuation on the
domains $\Pilh$ of the unphysical sheets $\Pi_l\subset\Re$.  This
means that such a continuation is possible as well for the
kernel $R(P,P',z)$ of $R(z)$.  Moreover, there exists an
explicit representation for this continuation analogous to the
representation (\ref{res2}) for the two\,--\,body resolvent.
\begin{theorem}\label{ThResolvl}
The analytic continuation, in the sense of distributions over $\Osix$,
of the resolvent $R(z)$ on the domain $\Pilh$ of the unphysical sheet
$\Pi_l\subset\Re$ is described by
\begin{small}
\begin{eqnarray}\label{R3l}
\reduction{R(z)}{\Pi_l}\!\! \!\! & = & \!\!\!\! R+
 \bigl( [I-RV]\rJot , \Om\hs{0.02cm}[\bI-\bRo M\Y]\Psi\rJt_1 \bigr)
LA S_l^{-1} \tL \!\left(\! \!\!
\begin{array}{c}  \rJo[I-VR]   \\ [0.1cm]   \rJ_1 \Psis[\bI-\Y M\bRo]\Omt
\end{array} \!\!\! \right).
\end{eqnarray}
\end{small}
The kernels of all the operators present on the right\,--\,hand side
of Eq.~{\rm(\ref{R3l})} are taken in the physical sheet.
\end{theorem}

\begin{proof}
For the analytic continuation $R^l(z)$ of the kernel  $R(P,P',z)$ of
$R(z)$ into the sheet $\Pi_l$ we have, according to
Eq.~(\ref{RMR}),
\be
\label{RMRl}
  R^l(z) \EQ R_0^l(z)\,-\,R_0^l(z) \hs{0.02cm}\Om  \hs{0.02cm}M^l(z)
   \hs{0.02cm}\Omt R_0^l(z)\,.
\ee
For $M^l(z)$ we have found already the
representation~(\ref{Ml3fin}).  Since $R^l_0=R_0+L_0
A_0\rJot\rJo$ we can rewrite Eq.~(\ref{RMRl}) in the form
\begin{eqnarray}
\nonumber
R^l \!\!& = &\!\! R_0\,-\,R_0 \hs{0.02cm}\Om  \hs{0.02cm}M^l\Omt R_0
\,+\,A_0 L_0 \rJot\!\left( \hIo-\rJo \hs{0.02cm}\Om  \hs{0.02cm}M^l\Omt\rJot
L_0 A_0 \right)\!\rJo
 \\ [-0.2cm]
& & \mbox{ } \\ [-0.2cm]
\label{Rlini}
 & &\!\! -\;\;A_0 L_0 \rJot\rJo \hs{0.02cm}\Om  \hs{0.02cm}M^l\Omt R_0 \,-\,
 R_0 \hs{0.02cm}\Om  \hs{0.02cm}M^l\Omt\rJot\rJo L_0 A_0\,. \nonumber
\end{eqnarray}
We consider separately the
contributions of each summand of (\ref{Rlini}).
To this end we shall use the notations
$$
\bB \EQ \Big( \Om M\Omt\rJot\, ,
\Om M\Y\Psi\rJt_1+\Om\Phi\rJt_1 \Big), \quad
\bBt \EQ \left(\!\!
\begin{array}{c}   \rJo\Om M\Omt \\ [0.1cm]
                 \rJ_1\Psis\Y M\Omt +\rJ_1\Phis\Omt
\end{array} \!\!\right).
$$

It follows from~(\ref{Ml3fin}) that
$
\Om M^l\Omt =\Om M \Omt -\bB LA S_l^{-1} \tL\bBt.
$
Hence, the first two summands of~(\ref{Rlini}) give together
$$
R_0\,-\,R_0 \hs{0.02cm}\Om
\hs{0.02cm}M\Omt R_0 \,+\, R_0\bB L A  S_l^{-1} \tL\bBt R_0 \EQ
R \,+\, R_0\bB L A S_l^{-1} \tL\bBt R_0\,.
$$

To transform the third term of~(\ref{Rlini}) we again use
the representation~(\ref{Ml3fin}). We find
\begin{eqnarray*}
\rJo\Om M^l\Omt\rJot L_0 A_0 \!\!& = & \!\!
\what{\cT}_{00} L_0 A_0 \,-\,
\left(\what{\cT}_{00}\, ,\what{\cT}_{01}\right)\!
L\, A S_l^{-1}\,\tL\left(\!\!    \begin{array}{c}
   \what{\cT}_{00}   \\ [0.1cm]
   \what{\cT}_{10}    \end{array}\!\! \right) \! L_0 A_0 \\ [0.1cm]
 & = & \!\! \omo \what{\cT} L  A \omos \,-\,
 \omo \what{\cT} L A S_l^{-1} \tL\what{\cT} L A\omos\\ [0.1cm]
& = & \!\! \omo\what{\cT} L A\!\left(  \hbI-S_l^{-1} \tL\what{\cT} L A
\right)\!\!\tL\omos
\end{eqnarray*}
where  $\omo$ stands for the projection from
$\what{\cH}_0\oplus\what{\cH}_1$ on  $\what{\cH}_0$ defined
by $\omo\!\left(\!\!
\begin{array}{c} f_0 \\ f_1 \end{array} \!\!\right)=f_0$,
$f_0\in \what{\cH}_0$, $f_1\in\what{\cH}_1$. By $\omos$ we
understand, as usual, the adjoint operator of $\omo$.  Since
$S_l=\hbI+ \tL\what{\cT} L A$ we have
$
    \hbI-S_l^{-1}\tL\what{\cT} L A=
$
$
   S_l^{-1}\!\!\left( \hbI+ \tL\what{\cT} L A-\tL\hat{\cT} L A\right)
$
$
   = S_l^{-1}.
$
Taking in account that $L=L\cdot\tL$ we find
$$
L_0 A_0 \big(\hIo-\rJo\Om M^l\Omt\rJot L_0 A_0\big) \EQ
\omo A L\big(\hbI - \tL\what{\cT} LA S_l^{-1}\big) \tL \omos \EQ
\omo L A\, S_l^{-1} \tL\omos\,.
$$
This means that the third term of~(\ref{Rlini}) may be represented
as $\rJot\omo L A S_l^{-1}\tL\omos.$

When studying the fourth summand of (\ref{Rlini})
we begin by transforming the product
$A_0 L_0\rJo \Om M^l \Omt$ into a more convenient form.
It follows from (\ref{Ml3fin}) that
$$
A_0 L_0\rJo\Om M^l\Omt \EQ A_0 L_0\rJo\Om M\Omt \,-\,
A_0 L_0\!\Big(\what{\cT}_{00}\, ,\what{\cT}_{01}\Big)LA S_l^{-1}\tL\bBt\,.
$$
In view of $A_0 L_0\rJo\Om M\Omt=\omo A L\bBt$ and
$
A_0 L_0\!\Big(\what{\cT}_{00},\what{\cT}_{01}\Big)\!L =
\omo AL\what{\cT} L
$
we have
$$
A_0 L_0\rJo \Om M^l\Omt \EQ \omo \!\left( AL-
A\, L\what{\cT}L A S_l^{-1} \tL\right) \!\bBt \EQ \omo L A S_l^{-1}\,\tL\bBt\,.
$$
Analogously, in the fifth term of (\ref{Rlini}) is
$$
\Om M^l\Omt\rJot L_0 A_0  \EQ
\bB\tL \bigl[S_l^{\dagger}\bigr]^{-1}A L \hs{0.02cm}\omos  \EQ
 \bB L A S_l^{-1} \tL\hs{0.02cm}\omos\,.
$$
Thus the last two summands of (\ref{Rlini}) give toge\-ther
$$
-\,\rJot\omo L A S_l^{-1}\tL\bBt R_0 \,- \,
R_0\bB L A S_l^{-1}\tL\hs{0.02cm}\omos\rJo\,.
$$
Substituting these expressions into Eq.~(\ref{Rlini})
we find
$$
R^l \EQ R+\big(\rJot\omo-R_0\bB\big)L A S_l^{-1}\tL\!
\left(\omos\rJo-\bBt R_0\right)\,.
$$
Taking into account the definitions of $\bB$ and $\bBt$ as well
as the obvious identities $R_0\Om M\Omt=RV$,
$R_0\Om\Phi\rJ_1=-\Om\Psi\rJ_1$ and
$\rJ_1\Phis\Omt R_0=-\rJ_1\Psis\Omt$
we come finally to Eq.~(\ref{R3l}) and this completes the proof.%
\end{proof}

\newsection{On the use of the
          Faddeev differential equations for computations
              of three\,--\,body resonances}
\label{SNumerMethod}
It follows from the representations~(\ref{Ml3fin}),
(\ref{Slfin}) and (\ref{R3l}) that the matrices
$\reduction{M(z)}{\Pi_l}$, $\reduction{S_{l'}(z)}{\Pi_l}$
and the resolvent $\reduction{R(z)}{\Pi_l}$ may have poles at points
belonging to the discrete spectrum $\sigma_d(H)$ of the
Hamiltonian $H$. Nontrivial singularities of
$\reduction{M(z)}{\Pi_l}$, $\reduction{S_{l'}(z)}{\Pi_l}$
and  $\reduction{R(z)}{\Pi_l}$ correspond to those points
$z\in\Pi_0\cap\Pilh$ where the inverse truncated scattering
matrix $[S_l(z)]^{-1}$ $\big($or $\big[\St_l(z)\big]^{-1}$ and this is the
same$\big)$ does not exist or where it represents an unbounded
operator.  The points $z$ where  $[S_l(z)]^{-1}$ does not exist
are the poles of $\reduction{M(z)}{\Pi_l}$,
$\reduction{S_{l'}(z)}{\Pi_l}$ and  $\reduction{R(z)}{\Pi_l}$.
Such points are called (three\,--\,body) resonances.

A necessary and sufficient condition~\cite{ReedSimonI} for
irreversibility of the operator $S_l(z)$ is the existence of a
nontrivial solution
$\cA^{\rm(res)}\in\what{\cH}_0\oplus\what{\cH}_1$ of the
equation
\be
\label{ZeroSl}
    S_l(z)\cA^{\rm (res)}\EQ 0\,.
\ee
The investigation of this equation may be carried out on the basis of
the results obtained in Sec.$\;$4 regarding the
properties of the kernels of the operator $\what{\cT}(z)$. We should
postpone this investigation for another paper.  Here, we
restrict ourselves to the observation  that the equation
(\ref{ZeroSl}) may evidently be applied to practical
computations of resonances situated in the domains
$\Pilh\subset\Pi_l$. The resonances have to be considered as
those values of $z\in\Pi_0\cap\Pilh$ for which the operators
$S_l(z)$ and $\St_l(z)$ have zero as eigenvalue.

The elements of the scattering matrices  $S_l(z)$ and $\St_l(z)$ are
expressed in terms of the amplitudes (continued in
$z$ into the physical sheet) for different processes taking place
in the three\,--\,body system under consideration. The respective
formulas~\cite{MF} rewritten for the components of $\what{\cT}$,
are the following
$$
\begin{array}{lcl}
\what{\cT}_{\an,j;\,\bn,k}\big(\hp\ad,\hp'\bd,z\big) \!\!&=& \!\!
C_0^{(3)}(z){\cA}_{\an,j;\,\bn,k}\big(\hp\ad,\hp'\bd,z\big)\, ,\\ [0.15cm]
\what{\cT}_{\an,j;\,0}\big(\hp\ad,\hP',z\big)\!\! &=& \!\!
C_0^{(3)}(z)\cA_{\an,j;\,0}\big(\hp\ad,\hP',z\big)\,, \\ [0.15cm]
\hat{\cT}_{0;\,\bn,\,k}\big(\hP,\hp'\bd,z\big)\!\! &=& \!\!
C_0^{(6)}(z)\cA_{0;\,\bn,k}\big(\hP,\hp'\bd,z\big)\,,\\ [0.15cm]
\what{\cT}_{00}\big(\hP,\hP',z\big)\!\! &=& \!\!
C_0^{(6)}(z)\cA_{00}\big(\hP,\hP',z\big)\,,
\end{array}
$$
with
$$
C_0^{(N)}(z) \EQ -\,\Frac{ { e}^{ i\pi(N-3)/4 }  }
{ 2^{(N-1)/2} \pi^{(N+1)/2} z^{(N-3)/4} }\,,
$$
where for the function $z^{(N-3)/4}$ one chooses the main
branch.  The functions  $\cA_{\an,j;\,\bn,k}$ represent the
amplitudes of elastic ($\an=\bn;$ $j=k$) or inelastic
($\an=\bn;$ $j\neq k$) scattering and rearrangement
($\an\neq\bn$) for the process ($2\rightarrow 2,3)$ in the
initial state of which the pair subsystem $\bn$ is in the $k$\,--\,th
bound state and the complementary particle is asymptotically
free.  The function $\cA_{0;\,\bn,k}$ represents in the same
process the amplitude of the system for breakup into three
separate particles.  The amplitudes  $\cA_{\an,j;\,0}$ and
$\cA_{00}$ correspond, respectively, to the processes
($3\rightarrow 2)$ and ($3\rightarrow 3)$ beginning from the
state with all three particles asymptotically free.  Recall that
the contributions to $\cA_{00}$ from the single and double
rescattering represent singular distributions (see
Sec.$\:$4).

By describing in Sec.$\;$4  the analytical properties of the
matrix $\what{\cT}$ kernels in the variable $z$ and the smoothness
properties in the angular variables $\hP$ or $\hp\ad$ and $\hP'$ or
$\hp'\bd$, we have in effect described as well the properties of
the amplitudes  $\cA(z)$.

To search for the amplitudes $\cA(z)$ continued into the
physical sheet, one can use e.\,g., the formulation~\cite{MF},
\cite{EChAYa} of the three\,--\,body scattering problem based on the
Faddeev differential equations in coordinate space. It is
only necessary to come, in this formulation, to complex values
of the energy  $z$.  The square roots $z^{1/2}$ and
$(z-\lambda\ajd)^{1/2}$,  $\alpha =$ $= 1,2,3, \; j=1,2, \,\ldots,
n_{\alpha}$, which are present in the
formulas of \cite{MF}, \cite{EChAYa} describing
asymptotical boundary conditions for the wave function
components at infinity, have to be considered as the main
branches, i.\,e., as $\sqrt{z}$ and $\sqrt{z-\lambda\ajd}$.
Solving the Faddeev differential equations with such conditions
one finds in fact the analytical continuation of the wave
functions into the physical sheet and, thus, the continuation
of the amplitudes  $\cA(z)$.  With the known amplitudes $\cA(z)$
one can construct a necessary truncated scattering matrix
$S_l(z)$ and find then those values of $z$ for which there exists a
nontrivial solution $\cA^{\rm (res)}$ of Eq.~(\ref{ZeroSl}). As
mentioned above, these values of $z$ represent the three\,--\,body
resonances in the respective the unphysical sheet $\Pi_l$.

\pagebreak
\begin{acknowledgements}

The author is much indebted to Prof. \sc S.~Albeverio \it for
the warm hospitality at the Ruhr-Uni\-ver\-si\-t\"{a}t-Bo\-chum
and to Prof. \sc K.\,A.~Ma\-ka\-rov \it for the stimulating
support.  The author is grateful to both of them as well as to
Prof. \sc V.\,B. Be\-lya\-ev \it for the valuable discussions
and to Prof. \sc P.~Ex\-ner \it for useful remarks.  Also, the
author is thankful to Prof. \sc R.~Mennicken \it for his
interest in this work.

The work was supported in part by the International Science
Foundation {\rm (}Grants~{\rm\#~RFB000} and {\rm\#~RFB300)} and
the Russian Foundation for Basic Research
{\rm (}Project {\rm\#~96\,--\,02\,--\,17021)}. \rm
\end{acknowledgements}

\address{Bogoliubov Laboratory of Theoretical Physics\\
Joint Institute for Nuclear Research\\
{\rm 141980} Dubna, Moscow Region\\ Russia
\\e\,--\,mail:\\
motovilv@thsun1.jinr.dubna.su}

\end{document}